%% file: 0908.3307.V1.English.tex
\documentclass{amsart}
\def\EprintCalculus{}
\input{Commands}
\xRefDef{0701.238}{math/pdf/0701/0701238}
\xRefDef{0812.4763}{math/pdf/0812/0812.4763v4.pdf}
\begin{document}
\pdfbookmark[1]{The G\^ateaux Derivative of Map over Division Ring}{TitleEnglish}
\title{The G\^ateaux Derivative of Map over Division Ring}

\begin{abstract}
I consider differential of mapping $f$ of continuous division ring
as linear mapping the most close to
mapping $f$.
Different expressions which correspond to known deffinition of
derivative are supplementary.
I explore the G\^ateaux derivative of higher order and
Taylor series. The Taylor series allow solving
of simple differential equations.
As an example of solution of differential equation
I considered a model of exponent.

I considered application of obtained theorems to complex field
and quaternion algebra. In contrast to complex field
in quaternion algebra congugation is linear
function of original number
\[
\overline a=a+iai+jaj+kak
\]
In quaternion algebra this difference leads to the absence of analogue
of the Cauchy Riemann equations
that are well known in the theory of complex function.
\end{abstract}
\maketitle
\tableofcontents

This paper is written on basis of eprints
\citeBib{0701.238,0812.4763}.
I submitted brief version of this paper to be
published in transactions of
the XXI International Summer School\Hyph Seminar
on Recent Problems in Theoretical and Mathematical Physics
"Volga'21\Hyph 2009" (XXI Petrov School).

In this paper I consider linear maps over division ring.
In particular, differential of map $f$ of division ring
is linear map.
However for this version I decided to add sections
dedicated to the G\^ateaux derivative of second order,
Taylor series expansion, and solving of differential equation.
The paper may be interesting for physicists using
division ring of quaternions.

\input{Convention.English}

\input{Ring.Additive.Map.English}

\input{Differential.English}
\input{Quaternion.English}

\input{Biblio.English}
\input{Index.English}
\input{Symbol.English}
\end{document}

%% file: Commands.tex
\def\Defined{}
\ifx\UseRussian\Defined
	\usepackage[T2A,T2B]{fontenc}
	\usepackage[cp1251]{inputenc}
	\usepackage[english,russian]{babel}
\fi
\usepackage{footmisc}
\usepackage[all]{xy}
\usepackage{color}
\definecolor{UrlColor}{rgb}{.9,0,.3}
\definecolor{SymbColor}{rgb}{.4,0,.9}
\definecolor{IndexColor}{rgb}{1,.3,.6}
\definecolor{eml1}{rgb}{.8,.1,.1}
\definecolor{eml2}{rgb}{.1,.6,.6}

\usepackage{xr-hyper}
\usepackage[unicode]{hyperref}
\hypersetup{pdfdisplaydoctitle=true}
\hypersetup{colorlinks}
\hypersetup{citecolor=UrlColor}
\hypersetup{urlcolor=UrlColor}
\hypersetup{pdffitwindow=true}
\hypersetup{pdfnewwindow=true}
\hypersetup{pdfstartview={FitH}}

\def\hyph{\penalty0\hskip0pt\relax-\penalty0\hskip0pt\relax}
\def\Hyph{-\penalty0\hskip0pt\relax}

\newcommand{\Basis}[1]{\overline{\overline{#1}}{}}
\newcommand{\Vector}[1]{\overline{#1}{}}
\newcommand{\gi}[1]{\boldsymbol{\textcolor{IndexColor}{#1}}}
\newcommand{\bVector}[2]{{}_{\gi{#2}}\Vector{#1}}
\makeatletter
\newcommand{\NameDef}[1]{%
	\expandafter\gdef\csname #1\endcsname%
}%
\newcommand{\ShowSymbol}[1]{%
	\@nameuse{ViewSymbol#1}%
}%
\newcommand{\symb}[3]{%
	\@ifundefined{ViewSymbol#3}{%
		\NameDef{ViewSymbol#3}{\textcolor{SymbColor}{#1}}%
		\NameDef{RefSymbol#3}{\pageref{symbol: #3}}%
		\@namedef{LabelSymbol#3}{\label{symbol: #3}}%
	}{%
		\NameDef{RefSymbol#3}{}%
		\@namedef{LabelSymbol#3}{}%
	}%
	\ifcase#2
	\or
		$\@nameuse{ViewSymbol#3}$%
	\or
		\[\@nameuse{ViewSymbol#3}\]%
	\else%
	\fi%
	\@nameuse{LabelSymbol#3}%
}%
\newcommand{\DefEq}[2]{%
	\@ifundefined{ViewEq#2}{%
		\NameDef{ViewEq#2}{#1}%
	}{%
	}%
}%
\newcommand{\ShowEq}[1]{%
	\@ifundefined{ViewEq#1}{%
		\message {missed ShowEq #1}
  }{%
	\@nameuse{ViewEq#1}%
	}%
}%
\makeatother

\newcommand{\subs}{${}_*$\Hyph}
\newcommand{\sups}{${}^*$\Hyph}

\newcommand{\CRstar}{{}_*{}^*}
\newcommand{\RCstar}{{}^*{}_*}

\newcommand{\RC}{$\RCstar$\Hyph}
\newcommand{\CR}{$\CRstar$\Hyph}
\newcommand{\drc}{$D\RCstar$\Hyph}
\newcommand{\Drc}{$\mathcal D\RCstar$\Hyph}
\newcommand{\dcr}{$D\CRstar$\hyph}
\newcommand{\rcd}{$\RCstar D$\Hyph}
\newcommand{\crd}{$\CRstar D$\Hyph}

\newcommand\sT{$\star T$\Hyph}%
\newcommand\Ts{$T\star$\Hyph}%
\newcommand\sD{$\star D$\Hyph}%
\newcommand\Ds{$D\star$\Hyph}%
\newcommand\pC[2]{{}_{(#1)#2}}%
\newcommand\DrcPartial[1]%
{%
	\def\tempa{}%
	\def\tempb{#1}%
	\ifx\tempa\tempc%
		(\partial\RCstar)%
	\else%
		({}_{#1}\partial\RCstar)%
	\fi%
}%
\newcommand\crDPartial[1]%
{%
	\def\tempa{}%
	\def\tempb{#1}%
	\ifx\tempa\tempc%
		(\CRstar\partial)%
	\else%
		(\CRstar\partial_{#1})%
	\fi%
}%
\newcommand\StandPartial[3]%
{%
	\frac{\partial^{\gi{#3}} #1}{\partial #2}%
}%

\renewcommand{\uppercasenonmath}[1]{}

\makeatletter
\newcommand\@dotsep{4.5}
\def\@tocline#1#2#3#4#5#6#7
{\relax
		\par \addpenalty\@secpenalty\addvspace{#2}%
		\begingroup 
		\@ifempty{#4}{%
			\@tempdima\csname r@tocindent\number#1\endcsname\relax
		}{%
			\@tempdima#4\relax
		}%
		\parindent\z@ \leftskip#3\relax \advance\leftskip\@tempdima\relax
		\rightskip\@pnumwidth plus1em \parfillskip-\@pnumwidth
		#5\leavevmode\hskip-\@tempdima #6\relax
		\leaders\hbox{$\m@th
		\mkern \@dotsep mu\hbox{.}\mkern \@dotsep mu$}\hfill
		\hbox to\@pnumwidth{\@tocpagenum{#7}}\par
		\nobreak
		\endgroup
}
\makeatother 

\ifx\PrintBook\undefined
	\usepackage{fancyhdr}
	\makeatletter
	\renewcommand{\@indextitlestyle}{%
		\twocolumn[\section{\indexname}]%
		\def\IndexSpace{off}%
	}
	\makeatother 
	\thanks{\href{mailto:Aleks\_Kleyn@MailAPS.org}{Aleks\_Kleyn@MailAPS.org}}
	\thanks{\ \ \ \url{http://www.geocities.com/aleks\_kleyn}}
\else

\makeatletter
\def\@maketitle{%
  \cleardoublepage \thispagestyle{empty}%
  \begingroup \topskip\z@skip
  \null\vfil
  \begingroup
  \LARGE\bfseries \centering
  \openup\medskipamount
  \@title\par\vspace{24pt}%
  \def\and{\par\medskip}\centering
  \mdseries\authors\par\bigskip
  \endgroup
  \vfil
  \ifx\@empty\addresses \else \@setaddresses \fi
  \vfil
  \ifx\@empty\@dedicatory
  \else \begingroup
    \centering{\footnotesize\itshape\@dedicatory\@@par}%
    \endgroup
  \fi
  \vfill
  \newpage\thispagestyle{empty}
  \@setabstract
  \begin{center}
    \ifx\@empty\@subjclass\else\@setsubjclass\fi
    \ifx\@empty\@keywords\else\@setkeywords\fi
    \ifx\@empty\@translators\else\vfil\@settranslators\fi
    \ifx\@empty\thankses\else\vfil\@setthanks\fi
  \end{center}
  \vfil
  \endgroup}
\makeatother 

	\makeatletter
	\renewcommand{\@indextitlestyle}{%
		\twocolumn[\chapter{\indexname}]%
		\def\IndexSpace{off}%
		\let\@secnumber\@empty
		\chaptermark{\indexname}%
	}
	\makeatother 
	\email{\href{mailto:Aleks\_Kleyn@MailAPS.org}{Aleks\_Kleyn@MailAPS.org}}
	\urladdr{\url{http://www.geocities.com/aleks\_kleyn}}
\fi

\ifx\SelectlEnglish\undefined
	\ifx\UseRussian\undefined
		\def\SelectlEnglish{}
	\fi
\fi

\newcommand\LanguagePrefix{}%
\ifx\SelectlEnglish\undefined
	\newcommand\CurrentLanguage{Russian.}%
	\author{Александр Клейн}
	\newtheorem{theorem}{Теорема}[section]
	\newtheorem{corollary}[theorem]{Следствие}
	\theoremstyle{definition}
	\newtheorem{definition}[theorem]{Определение}
	\newtheorem{example}[theorem]{Пример}
	\newtheorem{xca}[theorem]{Exercise}
	\theoremstyle{remark}
	\newtheorem{remark}[theorem]{Замечание}
	
	\newcommand\Gbasis{$G$\Hyph базис}
	\newcommand\Gcoords{$G$\Hyph координат}
	\newcommand\Gspace{$G$\Hyph пространств}
	\newcommand\xRefDef[2]
		{
		\externaldocument[#1-Russian-]{#1.Russian}[http://arxiv.org/PS_cache/#2.pdf]
		\NameDef{xRefDef#1}{}%
		}
	\makeatletter
	\newcommand\xRef[2]%
	{%
		\@ifundefined{xRefDef#1}{%
		\ref{#2}%
		}{%
		\citeBib{#1}-\ref{#1-Russian-#2}%
		}%
	}%
	\newcommand\xEqRef[2]%
	{%
		\@ifundefined{xRefDef#1}{%
		\eqref{#2}%
		}{%
		\citeBib{#1}-\eqref{#1-Russian-#2}%
		}%
	}%
	\makeatother
	\ifx\PrintBook\undefined
		\newcommand{\BibTitle}{%
			\section{Список литературы}%
		}
	\else
		\newcommand{\BibTitle}{%
			\chapter{Список литературы}%
		}
	\fi
\else
	\newcommand\CurrentLanguage{English.}%
	\author{Aleks Kleyn}
	\newtheorem{theorem}{Theorem}[section]
	
	\theoremstyle{definition}
	\newtheorem{definition}[theorem]{Definition}
	\newtheorem{example}[theorem]{Example}
	
	\theoremstyle{remark}
	\newtheorem{remark}[theorem]{Remark}
	\newcommand\Gbasis{$G$\Hyph basis}
	\newcommand\Gcoords{$G$\Hyph coordinates}
	\newcommand\Gspace{$G$\Hyph space}
	\newcommand\xRefDef[2]
		{
		\externaldocument[#1-English-]{#1.English}[http://arxiv.org/PS_cache/#2.pdf]
		\NameDef{xRefDef#1}{}%
		}
	\makeatletter
	\newcommand\xRef[2]%
	{%
		\@ifundefined{xRefDef#1}{%
		\ref{#2}%
		}{%
		\citeBib{#1}-\ref{#1-English-#2}%
		}%
	}%
	\newcommand\xEqRef[2]%
	{%
		\@ifundefined{xRefDef#1}{%
		\eqref{#2}%
		}{%
		\citeBib{#1}-\eqref{#1-English-#2}%
		}%
	}%
	\makeatother
	\ifx\PrintBook\undefined
		\newcommand{\BibTitle}{%
			\section{References}%
		}
	\else
		\newcommand{\BibTitle}{%
			\chapter{References}%
		}
	\fi
\fi

\ifx\PrintBook\undefined
	\numberwithin{Hfootnote}{section}
\else
	\numberwithin{section}{chapter}
	\numberwithin{footnote}{chapter}
	\numberwithin{Hfootnote}{chapter}
\fi

\numberwithin{equation}{section}
\numberwithin{figure}{section}
\numberwithin{table}{section}
\numberwithin{Item}{section}

\makeatletter
\newcommand\org@maketitle{}
\let\org@maketitle\maketitle
\def\maketitle{%
	\hypersetup{pdftitle={\@title}}%
	\hypersetup{pdfauthor={\authors}}%
	\hypersetup{pdfsubject=\@keywords}%
	\org@maketitle
}
\def\make@stripped@name#1{%
	\begingroup
		\escapechar\m@ne
		\global\let\newname\@empty
		\protected@edef\Hy@tempa{\CurrentLanguage #1}%
		\edef\@tempb{%
			\noexpand\@tfor\noexpand\Hy@tempa:=%
			\expandafter\strip@prefix\meaning\Hy@tempa
		}%
		\@tempb\do{%
			\if\Hy@tempa\else
				\if\Hy@tempa\else
					\xdef\newname{\newname\Hy@tempa}%
				\fi
			\fi
		}%
	\endgroup
}%
\newenvironment{enumBib}{%
	\BibTitle
	\advance\@enumdepth \@ne
	\edef\@enumctr{enum\romannumeral\the\@enumdepth}\list
	{\csname biblabel\@enumctr\endcsname}{\usecounter
	{\@enumctr}\def\makelabel##1{\hss\llap{\upshape##1}}}
}{%
	\endlist
}

\def\Chapters#1{\ChapterList#1,LastChapter,}%
\def\LastChapter{LastChapter}%
\def\ChapterList#1,{\def\temp{#1}%
	\ifx\temp\LastChapter
	\else
		\@ifundefined{#1}{%
		}{%
			\def\Semafor{on}
		}
		\expandafter\ChapterList
	\fi
}%
\makeatletter
\newcommand{\BiblioItem}[3]
{
	\def\Semafor{off}
	\@ifundefined{ViewCite#2}{}{%
		\def\Semafor{on}%
	}%
	\ifx\Semafor\ValueOff
		\@ifundefined{xRefDef#2}{}{%
		\def\Semafor{on}%
		}%
	\fi
	\ifx\Semafor\ValueOn
		\ifx\IndexState\ValueOff
			\begin{enumBib}
			\def\IndexState{on}
		\fi
		\item \label{\LanguagePrefix bibitem: #2}#3%
	\fi
}
\makeatother
\newcommand{\OpenBiblio}
{
	\def\IndexState{off}
}
\newcommand{\CloseBiblio}
{
	\ifx\IndexState\ValueOn
		\end{enumBib}
		\def\IndexState{off}
	\fi
}

\makeatletter
\def\StartCite{[}%
\def\citeBib#1{\protect\showCiteBib#1,endCite,}%
\def\endCite{endCite}%
\def\showCiteBib#1,{\def\temp{#1}%
\ifx\temp\endCite
]%
\def\StartCite{[}%
\else
	\StartCite\LanguagePrefix \ref{\LanguagePrefix bibitem: #1}%
	\@ifundefined{ViewCite#1}{%
		\NameDef{ViewCite#1}{}%
	}{%
	}%
	\def\StartCite{, }%
\expandafter\showCiteBib%
\fi}%
\makeatother

\makeatother 
\newcommand{\arp}{\ar @{-->}}

\newcommand{\bundle}[4]%
{%
	\def\tempa{}%
	\def\tempb{#3}%
	\def\tempc{#1}%
	\ifx\tempa\tempb%
		\ifx\tempa\tempc%
			#2%
		\else%
			\xymatrix{#2:#1\arp[r]&#4}%
		\fi%
	\else%
		\ifx\tempa\tempc%
			#2[#3]%
		\else%
			\xymatrix{#2[#3]:#1\arp[r]&#4}%
		\fi%
	\fi%
}%
\newcommand{\AddIndex}[2]%
{%
	{\bf #1}%
	\label{index: #2}%
}%
\newcommand{\Index}[3]%
{%
	\def\Semafor{off}%
	\Chapters{#1}%
	\ifx\Semafor\ValueOn%
		\def\tempa{}%
		\def\tempb{#3}%
		\ifx\IndexState\ValueOff%
			\begin{theindex}%
			\def\IndexState{on}%
		\fi%
		\ifx\IndexSpace\ValueOn%
			\indexspace%
			\def\IndexSpace{off}%
		\fi%
		\item #2%
		\ifx\tempa\tempb%
		\else%
			\ \pageref{index: #3}%
		\fi%
	\fi%
}%
\newcommand{\SubIndex}[3]
{
	\def\Semafor{off}
	\Chapters{#1}
	\ifx\Semafor\ValueOn
		\subitem #2 \pageref{index: #3}
	\fi
}%

\makeatletter
\newcommand{\Symb}[3]
{
	\def\Semafor{off}
	\Chapters{#1}
	\ifx\Semafor\ValueOn
		\ifx\IndexState\ValueOff
			\begin{theindex}
			\def\IndexState{on}
		\fi
		\ifx\IndexSpace\ValueOn
			\indexspace
			\def\IndexSpace{off}
		\fi
		\item $\displaystyle\@nameuse{ViewSymbol#3}$\ \ #2
		\@nameuse{RefSymbol#3}%
	\fi
}

\makeatother

\newcommand{\SetIndexSpace}%
{%
	\def\IndexSpace{on}%
}%
\def\ValueOff{off}
\def\ValueOn{on}

\newcommand{\OpenIndex}
{
	\def\IndexState{off}
}
\newcommand{\CloseIndex}
{
	\ifx\IndexState\ValueOn
		\end{theindex}
		\def\IndexState{off}
	\fi
}

\def\LastMemo{LastMemo}%
\def\MemoList#1//{\def\temp{#1}%
	\ifx\temp\LastMemo
	\else%
		\par
		\textcolor{blue}{#1}%
		\expandafter\MemoList%
	\fi%
}%

%% file: Convention.English.tex

\section{Conventions}

\begin{enumerate}

\item Function and map are synonyms. However according to
tradition, correspondence between either rings or vector
spaces is called map and map of
either real field or quaternion algebra is called function.
I also follow this tradition.

\item
We can consider division ring $D$ as $D$\Hyph vector space
of dimension $1$. According to this statement, we can explore not only
homomorphisms of division ring $D_1$ into division ring $D_2$,
but also linear maps of division rings.
This means that map is multiplicative over
maximum possible field. In particular, linear map
of division ring $D$ is multiplicative over center $Z(D)$. This statement
does not contradict with
definition of linear map of field because for field $F$ is true
$Z(F)=F$.
When field $F$ is different from
maximum possible, I explicit tell about this in text.

\item
In spite of noncommutativity of product a lot of statements
remain to be true if we substitute, for instance, right representation by
left representation or right vector space by left
vector space.
To keep this symmetry in statements of theorems
I use symmetric notation.
For instance, I consider \Ds vector space
and \sD vector space.
We can read notation \Ds vector space
as either D\Hyph star\Hyph vector space or
left vector space.
We can read notation \Ds linear dependent vectors
as either D\Hyph star\Hyph linear dependent vectors or
vectors that are linearly dependent from left.

\end{enumerate}

%% file: Ring.Additive.Map.English.tex
\def\texRingAdditiveMap{}
\ifx\PrintBook\undefined
\else
\chapter{Linear Map of Division Ring}
\label{chapter: Linear map, Division Ring}
\fi

\input{Ring.Additive.Map.Eq}

\section{Additive Map of Ring}

\begin{definition}
\label{definition: Additive map of Ring}
Homomorphism \[f:R_1\rightarrow R_2\] of additive group of ring $R_1$
into additive group of ring $R_2$
is called
\AddIndex{additive map of ring $R_1$
into ring $R_2$}{Additive map of Ring}.
\qed
\end{definition} 

According to definition of homomorphism of additive group,
additive map $f$ of ring $R_1$
into ring $R_2$
holds
\begin{equation}
f(a+b)=f(a)+f(b)
\label{eq: additive map, ring}
\end{equation}
\ifx\EprintCalculus\undefined
We do not expect that additive map
of ring holds product.
\fi

\begin{theorem}
\label{theorem: sum of additive maps, ring}
Let us consider ring $R_1$ and ring $R_2$.
Let maps \[f:R_1\rightarrow R_2\]
\[g:R_1\rightarrow R_2\] be additive maps.
Then map $f+g$ is additive.
\end{theorem}
\begin{proof}
Statement of theorem follows from chain of equations
\begin{align*}
(f+g)(x+y)
=&f(x+y)+g(x+y)
=f(x)+f(y)+g(x)+g(y)
\\
=&(f+g)(x)+(f+g)(y)
\end{align*}
\end{proof}

\begin{theorem}
\label{theorem: additive map times constant, ring}
Let us consider ring $R_1$ and ring $R_2$.
Let map \[f:R_1\rightarrow R_2\]
be additive map.
Then maps $af$, $fb$, $a$, $b\in R_2$ are additive.
\end{theorem}
\begin{proof}
Statement of theorem follows from chain of equations
\begin{align*}
(af)(x+y)
=&a(f(x+y))
=a(f(x)+f(y))
=af(x)+af(y)
\\
=&(af)(x)+(af)(y)
\\
(fb)(x+y)
=&(f(x+y))b
=(f(x)+f(y))b
=f(x)b+f(y)b
\\
=&(fb)(x)+(fb)(y)
\end{align*}
\end{proof}

\begin{theorem}
\label{theorem: additive map, ring, morphism}
We may represent additive map of ring $R_1$
into associative ring $R_2$ as
\ShowEq{additive map, ring, morphism}
where $G_{(s)}$ is set of additive maps of ring $R_1$
into ring $R_2$.\footnote{Here and in the following text we assume
sum over index that is written in brackets and used in product few times.
Equation \eqref{eq: additive map, ring, morphism}
is recursive definition and there is hope that it is possible
to simplify it.}
\end{theorem}
\begin{proof}
The statement of theorem follows from theorems
\ref{theorem: sum of additive maps, ring} and
\ref{theorem: additive map times constant, ring}. 
\end{proof}

\begin{definition}
\label{definition: map multiplicative over commutative ring, ring}
Let commutative ring $P$ be subring of center $Z(R)$ of ring $R$.
Map
\[f:R\rightarrow R\]
of ring $R$
is called
\AddIndex{multiplicative over commutative ring $P$}
{map multiplicative over commutative ring, ring}, if
\[
f(px)=pf(x)
\]
for any $p\in P$.
\qed
\end{definition}

\begin{definition}
Let commutative ring $F$ be subring of center $Z(D)$ of ring $R$.
Additive, multiplicative over commutative ring $F$, map
\[
f:R\rightarrow R
\]
is called \AddIndex{linear map over commutative ring $F$}
{linear map over commutative ring, ring}.
\qed
\end{definition}

\begin{definition}
\label{definition: map projective over commutative ring, ring}
Let commutative ring $P$ be subring of center $Z(R)$ of ring $R$.
Map
\[f:R\rightarrow R\]
of ring $R$
is called
\AddIndex{projective over commutative ring $P$}
{map projective over commutative ring, ring}, if
\[
f(px)=f(x)
\]
for any $p\in P$. Set
\[
Px=\{px:p\in P,x\in R\}
\]
is called
\AddIndex{direction $x$ over commutative ring $P$}
{direction over commutative ring, ring}.\footnote{Direction
over commutative ring $P$ is subset of ring $R$.
However we denote direction $Px$ by element
$x\in R$ when this does not lead to ambiguity.
We tell about direction over commutative ring $Z(R)$
when we do not show commutative ring $P$ explicitly.}
\qed
\end{definition}

\begin{example}
\label{example: map projective over commutative ring, ring}
If map $f$
of ring $R$ is multiplicative over commutative ring $P$,
then map
\[
g(x)=x^{-1}f(x)
\]
is projective over commutative ring $P$.
\qed
\end{example}

\begin{definition}
Denote \symb{\mathcal A(R_1;R_2)}1{set additive maps, ring}
set of additive maps
\[
f:R_1\rightarrow R_2
\]
of ring $R_1$
into ring $R_2$.
\qed
\end{definition}

\begin{theorem}
\label{theorem: additive map, ring, integer}
Let map
\[
f:D\rightarrow D
\]
is additive map of ring $R$.
Then
\[
f(nx)=nf(x)
\]
for any integer $n$.
\end{theorem}
\begin{proof}
We prove the theorem by induction on $n$.
Statement is obvious for $n=1$ because
\[
f(1x)=f(x)=1f(x)
\]
Let statement is true for $n=k$. Then
\[
f((k+1)x)=f(kx+x)=f(kx)+f(x)=kf(x)+f(x)=(k+1)f(x)
\]
\end{proof}

\section{Additive Map of Division Ring}

\begin{theorem}
\label{theorem: additive map, division ring, rational}
Let map
\ShowEq{map, division ring D1 D2}
is additive map of division ring $D_1$ into division ring $D_2$.
Then
\[
f(ax)=af(x)
\]
for any rational $a$.
\end{theorem}
\begin{proof}
Let $a=\frac pq$.
Assume $y=\frac 1qx$. Then
\begin{equation}
f(x)=f(qy)=qf(y)=qf\left(\frac 1qx\right)
\label{eq: additive map, division ring, rational, 1}
\end{equation}
From equation \eqref{eq: additive map, division ring, rational, 1}
it follows
\begin{equation}
\frac 1qf(x)=f\left(\frac 1qx\right)
\label{eq: additive map, division ring, rational, 2}
\end{equation}
From equation \eqref{eq: additive map, division ring, rational, 2}
it follows
\[
f\left(\frac pqx\right)=pf\left(\frac 1qx\right)=\frac pqf(x)
\]
\end{proof}

\begin{theorem}
\label{theorem: additive map multiplicative over rational field, division ring}
Additive map
\ShowEq{map, division ring D1 D2}
of division ring $D_1$ into division ring $D_2$
is multiplicative over field of rational numbers.
\end{theorem}
\begin{proof}
Corollary of theorem \ref{theorem: additive map, division ring, rational}.
\end{proof}

We cannot extend the statement of theorem
\ref{theorem: additive map multiplicative over rational field, division ring}
for arbitrary subfield of center $Z(D)$ of division ring $D$.

\begin{theorem}
\label{theorem: exists additive map not multiplicative over complex field, division ring}
Let complex field $C$ be subfield of
the center of division ring $D$.
There exists additive map
\ShowEq{map, division ring D1 D2}
of division ring $D_1$ into division ring $D_2$
which is not multiplicative over field of complex numbers.
\end{theorem}
\begin{proof}
To prove the theorem it is enough to consider the complex field
$C$ because $C=Z(C)$. The map
\[
z\rightarrow \overline z
\]
is additive. However the equation
\[
\overline{az}=a\overline z
\]
is not true.
\end{proof}

The theory of complex vector spaces so well understood that the
proof of theorem
\ref{theorem: exists additive map not multiplicative over complex field, division ring}
easily leads to the following design.
Let for some division ring $D$ fields $F_1$, $F_2$ be such that
$F_1\ne F_2$, $F_1\subset F_2\subset Z(D)$.
In this case there exists map $I$ of division ring $D$
that is linear over field $F_1$,
but not linear over field $F_2$.\footnote{For instance, in case
of complex numbers map $I$ is map of complex
conjugation. The set of maps $I$ depends on the
division ring. We consider these operators when we explorer
map of division ring when structure of operation changes. For instance,
the map of complex numbers $z\rightarrow\overline z$.}
It is easy to see that this map is additive.

Let $D_1$, $D_2$ be division rings of characteristic $0$.
According to theorem
\ref{theorem: additive map, ring, morphism}
additive map
\ShowEq{additive map D1 D2}
has form
\eqref{eq: additive map, ring, morphism}.
Let us choose map $G_{(s)}(x)=G(x)$.
Additive map
\ShowEq{additive map, ring, morphism G}
is called
\AddIndex{additive map generated by map $G$}
{additive map generated by map, division ring}.
Map $G$ is called
\AddIndex{generator of additive map}
{generator of additive map, division ring}.

\begin{theorem}
\label{theorem: additive map, division ring, center}
Let $F$, $F\subset Z(D_1)$, $F\subset Z(D_2)$, be field.
Additive map
\eqref{eq: additive map, ring, morphism G}
generated by $F$\Hyph linear map $G$ is multiplicative
over field $F$.
\end{theorem}
\begin{proof}
Immediate corollary of representation
\eqref{eq: additive map, ring, morphism G}
of additive map.
For any $a\in F$
\ShowEq{additive map, division ring, center}
\end{proof}

\begin{theorem}
\label{theorem: additive map based G, standard form, division ring}
Let $D_1$, $D_2$ be division rings of characteristic $0$.
Let $F$, $F\subset Z(D_1)$, $F\subset Z(D_2)$, be field.
Let $G$ be $F$\Hyph linear map.
Let $\Basis q$ be basis of division ring $D_2$ over field $F$.
\symb{f^{\gi i\gi j}_G}0
{standard component of additive map, division ring}
\AddIndex{Standard representation of additive map
\eqref{eq: additive map, ring, morphism G}
over field $F$}
{additive map, standard representation, division ring}
has form\footnote{Representation of additive map using
components of additive map
is ambiguous. We can increase or decrease number of summands
using algebraic operations.
Since dimension of division ring $D$ over field $F$ is finite, standard representation
of additive map guarantees finiteness of set
of items in the representation of map.}
\ShowEq{additive map, division ring, G, standard representation}
Expression
$\ShowSymbol{standard component of additive map, division ring}$
in equation \eqref{eq: additive map, division ring, G, standard representation}
is called
\AddIndex{standard component of additive map $f$ over field $F$}
{standard component of additive map, division ring}.
\end{theorem}
\begin{proof}
Components of additive map $f$
have expansion
\ShowEq{additive map, division ring, components extention}
relative to bais $\Basis q$.
If we substitute \eqref{eq: additive map, division ring, components extention}
into \eqref{eq: additive map, ring, morphism G},
we get
\ShowEq{additive map, division ring, G, standard representation, 1}
If we substitute expression
\ShowEq{additive map, division ring, G, standard representation, 2}
into equation \eqref{eq: additive map, division ring, G, standard representation, 1}
we get equation \eqref{eq: additive map, division ring, G, standard representation}.
\end{proof}

\begin{theorem}
\label{theorem: additive multiplicative map over field, G, division ring}
Let $D_1$, $D_2$ be division rings of characteristic $0$.
Let $F$, $F\subset Z(D_1)$, $F\subset Z(D_2)$, be field
Let $G$ be $F$\Hyph linear map.
Let $\Basis p$ be basis of division ring $D_1$ over field $F$.
Let $\Basis q$ be basis of division ring $D_2$ over field $F$.
Let ${}_{\gi{kl}}B^{\gi p}$ be structural constants of division ring $D_2$.
Then it is possible to represent additive map
\eqref{eq: additive map, ring, morphism G}
generated by $F$\Hyph linear map $G$
as
\ShowEq{additive multiplicative map over field, G, division ring}
\end{theorem}
\begin{proof}
According to theorem \ref{theorem: additive map, division ring, center}
additive map of division ring $D$ is linear over field $F$.
Let us consider map
\ShowEq{additive map D1 D2 generator}
According to theorem \xRef{0701.238}{theorem: drc linear map of drc vector space}
additive map $f(a)$ relative to basis $\Basis e$
has form
\eqref{eq: additive multiplicative map over field, G, division ring}.
From equations \eqref{eq: additive map, division ring, G, standard representation}
and \eqref{eq: additive map D1 D2 generator}
it follows
\ShowEq{additive map, division ring, F linear map G, standard representation}
From equations \eqref{eq: additive multiplicative map over field, G, division ring}
and \eqref{eq: additive map, division ring, F linear map G, standard representation}
it follows
\ShowEq{additive multiplicative map over field, G, division ring, 1}
Since vectors ${}_{\gi r}\Vector e$ are linear independent over field $F$ and values
$a^{\gi k}$ are arbitrary, then
equation \eqref{eq: additive multiplicative map over field, G, division ring, relation}
follows from equation
\eqref{eq: additive multiplicative map over field, G, division ring, 1}.
\end{proof}

\begin{theorem}
\label{theorem: additive generator is regular map, division ring}
Let field $F$ be subring of center $Z(D)$ of division ring $D$ of characteristic $0$.
$F$\Hyph linear map generating additive map is nonsingular
map.
\end{theorem}
\begin{proof}
According to isomorphism theorem we can represent additive map
\eqref{eq: linear map, division ring} as composition
\[
f(x)=f_1(x+H)
\]
of canonical map $x\rightarrow x+H$ and isomorphism $f_1$.
$H$ is ideal of additive group of division ring $D$.
Let ideal $H$ be non-trivial. Then
there exist $x_1\ne x_2$, $f(x_1)=f(x_2)$.
Therefore, image under map $f$ contains
cyclic subgroup.
It conflict with statement that characteristic of division ring $D$ equal $0$.
Therefore, either $H=\{0\}$ and canonical map is nonsingular $F$\Hyph linear
map or $H=D$ and canonical map is singular
map.
\end{proof}

\begin{definition}
\label{definition: linear map, division ring}
Additive map that is linear over
center of division ring is called \AddIndex{linear map of division ring}
{linear map, division ring}.
\qed
\end{definition}

\begin{theorem}
\label{theorem: linear map, division ring, canonical morphism}
Let $D$ be division ring of characteristic $0$.
Linear map
\ShowEq{linear map, division ring}
has form
\ShowEq{linear map, division ring, canonical morphism}
Expression
\symb{\pC{s}{p}f}1
{component of linear map, division ring}, $p=0$, $1$,
in equation \eqref{eq: linear map, division ring, canonical morphism}
is called
\AddIndex{component of linear map $f$}
{component of linear map, division ring}.
\end{theorem}

\begin{theorem}
\label{theorem: linear map, standard form, division ring}
Let $D$ be division ring of characteristic $0$.
Let $\Basis e$ be the basis of division ring $D$ over center $Z(D)$.
\symb{f^{\gi i\gi j}}0
{standard component of linear map, division ring}
\AddIndex{Standard representation of linear map
\eqref{eq: linear map, division ring, canonical morphism}
of division ring}
{linear map, standard representation, division ring}
has form\footnote{Representation of linear map of of division ring using
components of linear map
is ambiguous. We can increase or decrease number of summands
using algebraic operations.
Since dimension of division ring $D$ over field $Z(D)$ is finite, standard representation
of linear map guarantees finiteness of set
of items in the representation of map.}
\ShowEq{linear map, division ring, standard representation}
Expression
$\ShowSymbol{standard component of linear map, division ring}$
in equation \eqref{eq: linear map, division ring, standard representation}
is called
\AddIndex{standard component of linear map $f$}
{standard component of linear map, division ring}.
\end{theorem}

\begin{theorem}
\label{theorem: linear map over field, division ring}
Let $D$ be division ring of characteristic $0$.
Let $\Basis e$ be basis of division ring $D$ over field $Z(D)$.
Then it is possible to represent linear map
\eqref{eq: linear map, division ring}
as
\ShowEq{linear map over field, division ring}
\end{theorem}
\begin{proof}
Equation
\eqref{eq: linear map, division ring, canonical morphism}
is special case of equation
\eqref{eq: additive map, ring, morphism G}
when $G(x)=x$.
Theorem
\ref{theorem: linear map, standard form, division ring}
is special case of theorem
\ref{theorem: additive map based G, standard form, division ring}
when $G(x)=x$.
Theorem
\ref{theorem: linear map over field, division ring}
is special case of theorem
\ref{theorem: additive multiplicative map over field, G, division ring}
when $G(x)=x$.
Our goal is to show that we can assume $G(x)=x$.

Equation
\eqref{eq: linear map over field, division ring, relation}
binds coordinates of linear transformation $f$ relative given
basis $\Basis e$ of division ring $D$ over center $Z(D)$
and standard components of this transformation when
we consider this transformation as linear map of division ring.
For given coordinates of linear transformation
we consider equation
\eqref{eq: linear map over field, division ring, relation}
as the system of linear equations in
standard components.
From theorem \ref{theorem: linear map over field, division ring}
it follows that if determinant of the system of linear equations
\eqref{eq: linear map over field, division ring, relation}
is different from $0$, then for any linear transformation of division ring $D$ over field $Z(D)$
there exists corresponding linear transformation of division ring.

If determinant of system of linear equations
\eqref{eq: linear map over field, division ring, relation}
equal $0$, then there exist maps different from map $G(x)=x$.
If dimension of division ring is finite then these maps form
finite dimensional algebra. I will consider the structure of this algebra in separate paper.
Without loss of generality, we will assume $G(x)=x$ also in this case.
\end{proof}

\begin{theorem}
\label{theorem: linear map, division ring}
Let field $F$ be subring of center $Z(D)$ of division ring $D$ of characteristic $0$.
Linear map of division ring
is multiplicative
over field $F$.
\end{theorem}
\begin{proof}
Immediate corollary of definition
\ref{definition: linear map, division ring}.
\end{proof}


\begin{theorem}
Expression
\[
{}_{}f^{\gi k\gi r}=f^{\gi i\gi j}
\ {}_{\gi{ik}}B^{\gi p}\ {}_{\gi{pj}}B^{\gi r}
\]
is tensor over field $F$
\begin{equation}
{}_{\gi i}f'^{\gi j}={}_{\gi i}A^{\gi k}
\ {}_{\gi k}f^{\gi l}\ {}_{\gi l}A^{-1}{}^{\gi j}
\label{eq: standard component of additive map, division ring, transformation}
\end{equation}
\end{theorem}
\begin{proof}
$D$\Hyph linear map has form
\eqref{eq: linear map over field, division ring}
relative to basis $\Basis e$.
Let $\Basis e{}'$ be another basis. Let
\begin{equation}
{}_{\gi i}\Vector e{}'={}_{\gi i}A^{\gi j}\ {}_{\gi j}\Vector e
\label{eq: standard component of additive map, division ring, change basis}
\end{equation}
be transformation mapping basis $\Basis e$ to
basis $\Basis e{}'$.
Since additive map $f$ is the same, then
\begin{equation}
f(x)=x'^{\gi k}\ {}_{\gi k}f'^{\gi l}\ {}_{\gi l}\Vector e{}'
\label{eq: standard component of additive map, division ring, transformation, 1}
\end{equation}
Let us substitute \xEqRef{0701.238}{eq: division ring over field, change basis, 3},
\eqref{eq: standard component of additive map, division ring, change basis}
into equation \eqref{eq: standard component of additive map, division ring, transformation, 1}
\begin{equation}
f(x)
=x^{\gi i}\ {}_{\gi i}A^{-1}{}^{\gi k}\ {}_{\gi k}f'^{\gi l}
\ {}_{\gi l}A^{\gi j}\ {}_{\gi j}\Vector e
\label{eq: standard component of additive map, division ring, transformation, 2}
\end{equation}
Because vectors ${}_{\gi j}\Vector e$ are linear independent and
components of vector $x^{\gi i}$ are arbitrary, the equation
\eqref{eq: standard component of additive map, division ring, transformation}
follows
from equation
\eqref{eq: standard component of additive map, division ring, transformation, 2}.
Therefore, expression
${}_kf^r$
is tensor over field $F$.
\end{proof}

\ifx\EprintCalculus\undefined
\begin{definition}
The set
\symb{\text{ker}f}
0{kernel of additive map, division ring}
\[
\ShowSymbol{kernel of additive map, division ring}=\{x\in D_1:f(x)=0\}
\]
is called
\AddIndex{kernel of additive map}{kernel of additive map, division ring}
\[f:D_1\rightarrow D_2\]
of division ring $D_1$
into division ring $D_2$.
\qed
\end{definition}

\begin{theorem}
Kernel of additive map
\[f:D_1\rightarrow D_2\]
is subgroup of additive group of division ring $D_1$.
\end{theorem}
\begin{proof}
Let $a$, $b\in\text{ker}f$. Then
\begin{align*}
f(a)&=0
\\
f(b)&=0
\\
f(a+b)=f(a)+f(b)&=0
\end{align*}
Therefore, $a+b\in\text{ker}f$.
\end{proof}

\begin{definition}
The additive map
\[f:D_1\rightarrow D_2\]
of division ring $D_1$
into division ring $D_2$
is called
\AddIndex{singular}{singular additive map, division ring},
when
\[\text{ker}f\ne\{0\}\]
\qed
\end{definition}
\fi

\begin{theorem}
\label{theorem: product of additive map, D D D}
Let $D$ be division ring of characteristic $0$.
Let $\Basis e$ be basis of division ring $D$ over center $Z(D)$ of division ring $D$.
Let
\ShowEq{product of additive map, D D D}
be linear maps of division ring $D$.
Map
\begin{equation}
h(x)=gf(x)=g(f(x))
\label{eq: D linear map, division ring, h1}
\end{equation}
is linear map
\ShowEq{D linear map, division ring, h2}
where
\ShowEq{D linear map, division ring, h gf}
\end{theorem}
\begin{proof}
Map \eqref{eq: D linear map, division ring, h1}
is additive because
\ShowEq{product of additive map, D D D, 1}
Map \eqref{eq: D linear map, division ring, h1}
is multiplicative over $Z(D)$ because
\ShowEq{product of additive map, D D D, 2}

If we substitute
\eqref{eq: D linear map, division ring, f} and
\eqref{eq: D linear map, division ring, g} into
\eqref{eq: D linear map, division ring, h1},
we get
\begin{equation}
h(x)=\pC{t}{0}g\ f(x)\ \pC{t}{1}g
=\pC{t}{0}g\ \pC{s}{0}f\ x\ \pC{s}{1}f\ \pC{t}{1}g
\label{eq: additive map, division ring, h3}
\end{equation}
Comparing \eqref{eq: additive map, division ring, h3}
and \eqref{eq: D linear map, division ring, h2},
we get \eqref{eq: D linear map, division ring, h gf, 0},
\eqref{eq: D linear map, division ring, h gf, 1}.

If we substitute
\eqref{eq: D linear map, division ring, standard representation, f} and
\eqref{eq: D linear map, division ring, standard representation, g} into
\eqref{eq: D linear map, division ring, h1},
we get
\begin{align}
h(x)=&g^{\gi{ij}}\ {}_{\gi i}\Vector e\ f(x)\ {}_{\gi j}\Vector e
\nonumber
\\
=&g^{\gi{ij}}\ {}_{\gi i}\Vector e\ f^{\gi{kl}}\ {}_{\gi k}\Vector e
\ x\ {}_{\gi l}\Vector e\ {}_{\gi j}\Vector e
\label{eq: D linear map, division ring, standard representation, h3}
\\
=&g^{\gi{ij}}\ f^{\gi{kl}}\ {}_{\gi{ik}}B^{\gi p}
\ {}_{\gi{lj}}B^{\gi r}\ {}_{\gi p}e\ x\ {}_{\gi r}\Vector e
\nonumber
\end{align}
Comparing \eqref{eq: D linear map, division ring, standard representation, h3}
and \eqref{eq: D linear map, division ring, standard representation, h2},
we get \eqref{eq: D linear map, division ring, standard representation, h gf}.
\end{proof}

\section{Polylinear Map of Division Ring}

\begin{definition}
\label{definition: polyadditive map of division rings}
Let $R_1$, ..., $R_n$ be
rings and $S$ be module.
We call map
\begin{equation}
f:R_1\times...\times R_n\rightarrow
S
\label{eq: polyadditive map of rings}
\end{equation}
\AddIndex{polyadditive map of rings}{polyadditive map of rings}
$R_1$, ..., $R_n$
into module
$S$,
if
\[
f(
p_1, ...,
p_i+ q_i, ...,
p_n)
=
f(
p_1, ...,
p_i, ...,
p_n)
+
f(
p_1, ...,
q_i, ...,
p_n)
\]
for any $1\le i\le n$ and
for any $p_i$, $q_i \in R_i$.
Let us denote
\symb{\mathcal A(R_1,...,R_n;
S)}1{set polyadditive maps, ring}
set of polyadditive maps
of rings
$R_1$, ..., $R_n$
into module
$S$.
\qed
\end{definition}

\begin{theorem}
\label{theorem: map polylinear over commutative ring, ring}
Let $R_1$, ..., $R_n$, $P$ be rings of characteristic $0$.
Let $S$ be module over ring $P$.
Let
\[
f:R_1\times...\times R_n\rightarrow
S
\]
be polyadditive map.
There exists commutative ring $F$ which is for any $i$ is subring of center
of ring $R_i$ and such that
for any $i$ and $b\in F$
\[
f(a_1, ..., b a_i, ..., a_n)=bf(a_1, ..., a_i, ..., a_n)
\]
\end{theorem}
\begin{proof}
For given $a_1$, ..., $a_{i-1}$, $a_{i+1}$, ..., $a_n$
map $f(a_1, ..., a_n)$ is additive by $a_i$.
According to theorem
\ref{theorem: additive map, ring, integer},
we can select ring of integers as ring $F$.
\end{proof}

\begin{definition}
\label{definition: map polylinear over commutative ring, ring}
Let $R_1$, ..., $R_n$, $P$ be rings of characteristic $0$.
Let $S$ be module over ring $P$.
Let $F$ be commutative ring which is for any $i$ is subring of center
of ring $R_i$.
Map
\[
f:R_1\times...\times R_n\rightarrow
S
\]
is called \AddIndex{polylinear over commutative ring $F$}
{map polylinear over commutative ring, ring},
if map $f$
is polyadditive, and for any $i$, $1\le i\le n$,
for given $a_1$, ..., $a_{i-1}$, $a_{i+1}$, ..., $a_n$
map $f(a_1, ..., a_n)$ is multiplicative by $a_i$.
If ring $F$ is maximum ring such that for any $i$, $1\le i\le n$,
for given $a_1$, ..., $a_{i-1}$, $a_{i+1}$, ..., $a_n$
map $f(a_1, ..., a_n)$ is linear by $a_i$ over ring $F$,
then map $f$ is called
\AddIndex{polylinear map of rings}{polylinear map of rings}
$R_1$, ..., $R_n$
into module
$S$.
Let us denote
\symb{\mathcal L(R_1,...,R_n;
S)}1{set polylinear maps, ring}
set of polylinear maps
of rings
$R_1$, ..., $R_n$
into module
$S$.
\qed
\end{definition}

\begin{theorem}
\label{theorem: polylinear map, division ring}
Let $D$ be division ring of characteristic $0$.
Polylinear map
\ShowEq{polylinear map, division ring}
has form
\ShowEq{polylinear map, division ring, canonical morphism}
$\sigma_s$ is a transposition of set of variables
$\{d_1,...,d_n\}$
\ShowEq{transposition of set of variables, division ring}
\end{theorem}
\begin{proof}
We prove statement by induction on $n$.

When $n=1$ the statement of theorem is corollary of theorem
\ref{theorem: linear map, division ring, canonical morphism}.
In such case we may identify\footnote{In representation
\eqref{eq: polylinear map, division ring, canonical morphism}
we will use following rules.
\begin{itemize}
\item If range of any index is set
consisting of one element, then we will omit corresponding
index.
\item If $n=1$, then $\sigma_s$ is identical transformation.
We will not show such transformation in the expression.
\end{itemize}} ($p=0$, $1$)
\ShowEq{polylinear map, division ring, 1, canonical morphism}

Let statement of theorem be true for $n=k-1$.
Then it is possible to represent map
\eqref{eq: polylinear map, division ring}
as
\ShowEq{polylinear map, induction on n, 1, division ring}
According to statement of induction polyadditive map
$h$ has form
\ShowEq{polylinear map, induction on n, 2, division ring}
According to construction $h=g(d_k)$.
Therefore, expressions $\pC{t}{p}h$
are functions of $d_k$.
Since $g(d_k)$ is additive map of $d_k$,
then only one expression $\pC{t}{p}h$
is additive map of $d_k$, and rest expressions
$\pC{t}{q}h$
do not depend on $d_k$.

Without loss of generality, assume $p=0$.
According to equation
\eqref{eq: linear map, division ring, canonical morphism}
for given $t$
\ShowEq{polylinear map, induction on n, 3, division ring}
Assume $s=tr$. Let us define transposition $\sigma_s$ according to rule
\ShowEq{polylinear map, induction on n, 4, division ring}
Suppose
\ShowEq{polylinear map, induction on n, 5, division ring}
for $q=1$, ..., $k-1$.
\ShowEq{polylinear map, induction on n, 6, division ring}
for $q=0$, $1$.
We proved step of induction.
\end{proof}

\begin{definition}
\begin{sloppypar}
Expression
\symb{\pC{s}{p}f^n}1
{component of polylinear map, division ring}
in equation \eqref{eq: polylinear map, division ring, canonical morphism}
is called
\AddIndex{component of polylinear map $f$}
{component of polylinear map, division ring}.
\qed
\end{sloppypar}
\end{definition}

\begin{theorem}
\label{theorem: standard representation of polylinear map, division ring}
Let $D$ be division ring of characteristic $0$.
Let $\Basis e$
be basis in division ring $D$ over field $Z(D)$.
\symb{\pC{t}{}f^{\gi{i_0...i_n}}}0
{standard component of polylinear map, division ring}
\AddIndex{Standard representation of polylinear map of division ring}
{polylinear map, standard representation, division ring}
has form
\ShowEq{polylinear map, division ring, standard representation}
Index $t$ enumerates every possible transpositions
$\sigma_t$ of the set of variables
$\{d_1,...,d_n\}$.
Expression
$\ShowSymbol{standard component of polylinear map, division ring}$
in equation \eqref{eq: polylinear map, division ring, standard representation}
is called
\AddIndex{standard component of polylinear map $f$}
{standard component of polylinear map, division ring}.
\end{theorem}
\begin{proof}
Components of polylinear map $f$
have expansion
\ShowEq{polylinear map, division ring, components extention}
relative to basis $\Basis e$.
If we substitute \eqref{eq: polylinear map, division ring, components extention}
into \eqref{eq: polylinear map, division ring, canonical morphism},
we get
\ShowEq{polylinear map, division ring, standard representation, 1}
Let us consider expression
\ShowEq{polylinear map, division ring, standard representation, 2}
The right-hand side is supposed to be the sum of the terms
with the index $s$, for which the transposition $\sigma_s$ is the same.
Each such sum
has a unique index $t$.
If we substitute expression
\eqref{eq: polylinear map, division ring, standard representation, 2}
into equation \eqref{eq: polylinear map, division ring, standard representation, 1}
we get equation \eqref{eq: polylinear map, division ring, standard representation}.
\end{proof}

\begin{theorem}
\label{theorem: polylinear map over field, division ring}
Let $\Basis e$ be basis of division ring $D$ over field $Z(D)$.
Polyadditive map
\eqref{eq: polylinear map, division ring}
can be represented as $D$-valued form of degree $n$
over field $Z(D)$\footnote{We proved the theorem by analogy with
theorem in \citeBib{Rashevsky}, p. 107, 108}
\ShowEq{polylinear map over field, division ring}
where
\ShowEq{polylinear map over field, coordinates, division ring}
and values ${}_{\gi{i_1...i_n}}f$
are coordinates of $D$-valued covariant tensor over field $F$.
\end{theorem}
\begin{proof}
According to theorem
\ref{theorem: map polylinear over commutative ring, ring},
the equation
\eqref{eq: polylinear map over field, division ring}
follows from the chain of equations
\ShowEq{polylinear map over field, 1, division ring}
Let $\Basis e'$ be another basis. Let
\ShowEq{polylinear map over field, division ring, change basis}
be transformation, mapping basis $\Basis e$ into
basis $\Basis e'$.
From equations \eqref{eq: polylinear map over field, division ring, change basis}
and \eqref{eq: polylinear map over field, coordinates, division ring}
it follows
\ShowEq{polylinear map over field, division ring, change coordinates}
From equation \eqref{eq: polylinear map over field, division ring, change coordinates}
the tensor law of transformation of coordinates of polylinear map follows.
From equation \eqref{eq: polylinear map over field, division ring, change coordinates}
and theorem \xRef{0701.238}{theorem: division ring over field, change basis}
it follows that value of the map $f(\Vector a_1,...,\Vector a_n)$ does not depend from
choice of basis.
\end{proof}

Polylinear map
\eqref{eq: polylinear map, division ring}
is \AddIndex{symmetric}
{polylinear map symmetric, division ring}, if
\[
f(d_1,...,d_n)=f(\sigma(d_1),...,\sigma(d_n))
\]
for any transposition $\sigma$ of set $\{d_1,...,d_n\}$.

\begin{theorem}
\label{theorem: polylinear map symmetric, division ring}
If polyadditive map $f$ is symmetric,
then
\ShowEq{polylinear map symmetric, division ring}
\end{theorem}
\begin{proof}
Equation
\eqref{eq: polylinear map symmetric, division ring}
follows from equation
\ShowEq{polylinear map symmetric, 1, division ring}
\end{proof}

Polylinear map
\eqref{eq: polylinear map, division ring}
is \AddIndex{skew symmetric}
{polylinear map skew symmetric, division ring}, if
\[
f(d_1,...,d_n)=|\sigma|f(\sigma(d_1),...,\sigma(d_n))
\]
for any transposition $\sigma$ of set $\{d_1,...,d_n\}$.
Here
\[
|\sigma|=
\left\{
\begin{matrix}
1&\textrm{transposition }\sigma\textrm{ even}
\\
-1&\textrm{transposition }\sigma\textrm{ odd}
\end{matrix}
\right.
\]

\begin{theorem}
\label{theorem: polylinear map skew symmetric, division ring}
If polylinear map $f$ is skew symmetric,
then
\ShowEq{polylinear map skew symmetric, division ring}
\end{theorem}
\begin{proof}
Equation
\eqref{eq: polylinear map skew symmetric, division ring}
follows from equation
\ShowEq{polylinear map skew symmetric, 1, division ring}
\end{proof}

\begin{theorem}
\label{theorem: coordinates of polylinear map, division ring over field}
The polylinear over field $F$ map
\eqref{eq: polylinear map, division ring}
is polylinear iff
\ShowEq{coordinates of polylinear map, division ring over field}
\end{theorem}
\begin{proof}
In equation \eqref{eq: polylinear map, division ring, standard representation},
we assume
\ShowEq{coordinates of polylinear map, 1, division ring over field}
Then equation \eqref{eq: polylinear map, division ring, standard representation}
gets form
\ShowEq{polylinear map, division ring, standard representation, 1a}
From equation \eqref{eq: polylinear map over field, division ring}
it follows that
\ShowEq{polylinear map over field, division ring, 1a}
Equation
\eqref{eq: coordinates of polylinear map, division ring over field}
follows from comparison of equations
\eqref{eq: polylinear map, division ring, standard representation, 1a}
and
\eqref{eq: polylinear map over field, division ring}.
Equation
\eqref{eq: coordinates of polylinear map, division ring over field, 1}
follows from comparison of equations
\eqref{eq: polylinear map, division ring, standard representation, 1a}
and
\eqref{eq: polylinear map over field, division ring, 1a}.
\end{proof}

%% file: Ring.Additive.Map.Eq.tex

\DefEq
{
\begin{equation}
f(x)=\pC{s}{0}f\ G_{(s)}(x)\ \pC{s}{1}f
\label{eq: additive map, ring, morphism}
\end{equation}
}
{additive map, ring, morphism}

\DefEq
{
\begin{equation}
f(x)=\pC{s}{0}f\ G(x)\ \pC{s}{1}f
\label{eq: additive map, ring, morphism G}
\end{equation}
}
{additive map, ring, morphism G}

\DefEq
{
\begin{align}
f(a)=&a^{\gi i}\ {}_{\gi i}f^{\gi j}\ {}_{\gi j}\Vector e
&{}_{\gi k}f^{\gi j}&\in Z(D)
\label{eq: linear map over field, division ring}
\\
a=&a^{\gi i}\ {}_{\gi i}\Vector e&a^{\gi i}&\in Z(D)\ \ \ a\in D
\nonumber
\\
{}_{\gi i}f^{\gi j}=&f^{\gi{kr}}\ {}_{\gi{ki}}B^{\gi p}
\ {}_{\gi{pr}}B^{\gi j}
\label{eq: linear map over field, division ring, relation}
\end{align}
}
{linear map over field, division ring}

\DefEq
{
\begin{align}
f(a)=&a^{\gi i}\ {}_{\gi i}f^{\gi j}\ {}_{\gi j}\Vector q
&{}_{\gi k}f^{\gi j}&\in F
\label{eq: additive multiplicative map over field, G, division ring}
\\
a=&a^{\gi i}\ {}_{\gi i}\Vector p&a^{\gi i}&\in F\ \ \ a\in D_1
\nonumber
\\
{}_{\gi i}f^{\gi j}=&{}_{\gi i}G^{\gi l}\ f^{\gi{kr}}_G\ {}_{\gi{kl}}B^{\gi p}
\ {}_{\gi{pr}}B^{\gi j}
\label{eq: additive multiplicative map over field, G, division ring, relation}
\end{align}
}
{additive multiplicative map over field, G, division ring}

\DefEq
{
\begin{equation}
f:D\rightarrow D
\label{eq: linear map, division ring}
\end{equation}
}
{linear map, division ring}

\DefEq
{
\begin{equation}
f(x)=\pC{s}{0}f\ x\ \pC{s}{1}f
\label{eq: linear map, division ring, canonical morphism}
\end{equation}
}
{linear map, division ring, canonical morphism}

\DefEq
{
\begin{equation}
f(x)=\pC{s}{0}f^{\gi i}\ {}_{\gi i}\Vector e
\ x\ \pC{s}{1}f^{\gi j}\ {}_{\gi j}\Vector e
\label{eq: additive map, division ring, standard representation, 1}
\end{equation}
}
{additive map, division ring, standard representation, 1}

\DefEq
{
\begin{equation}
f(x)=\pC{s}{0}f^{\gi i}\ {}_{\gi i}\Vector q
\ G(x)\ \pC{s}{1}f^{\gi j}\ {}_{\gi j}\Vector q
\label{eq: additive map, division ring, G, standard representation, 1}
\end{equation}
}
{additive map, division ring, G, standard representation, 1}

\DefEq
{
\[
f^{\gi i\gi j}=\pC{s}{0}f^{\gi i}\ \pC{s}{1}f^{\gi j}
\]
}
{additive map, division ring, standard representation, 2}

\DefEq
{
\[
f^{\gi i\gi j}_G=\pC{s}{0}f^{\gi i}\ \pC{s}{1}f^{\gi j}
\]
}
{additive map, division ring, G, standard representation, 2}

\DefEq
{
\begin{equation}
\pC{s}{p}f=\pC{s}{p}f^{\gi i}\ {}_{\gi i}\Vector q
\label{eq: additive map, division ring, components extention}
\end{equation}
}
{additive map, division ring, components extention}

\DefEq
{
\begin{equation}
f(x)=\ShowSymbol{standard component of linear map, division ring}
\ {}_{\gi i}\Vector e\ x\ {}_{\gi j}\Vector e
\label{eq: linear map, division ring, standard representation}
\end{equation}
}
{linear map, division ring, standard representation}

\DefEq
{
\begin{equation}
f(x)=\ShowSymbol{standard component of additive map, division ring}
\ {}_{\gi i}\Vector q\ G(x)\ {}_{\gi j}\Vector q
\label{eq: additive map, division ring, G, standard representation}
\end{equation}
}
{additive map, division ring, G, standard representation}

\DefEq
{
\[
f(ax)=\pC{s}{0}f\ G(ax)\ \pC{s}{1}f=\pC{s}{0}f\ aG(x)\ \pC{s}{1}f
=a\ \pC{s}{0}f\ G(x)\ \pC{s}{1}f=af(x)
\]
}
{additive map, division ring, center}

\DefEq
{
\begin{equation}
\begin{matrix}
G:D_1\rightarrow D_2
&
a=a^{\gi i}\ \bVector{p}{i}\rightarrow G(a)=a^{\gi i}
\ {}_{\gi i}G^{\gi j}\ \bVector{q}{j}
\\
&
x^{\gi i}\in F\ \ \ {}_{\gi i}G^{\gi j}\in F
\end{matrix}
\label{eq: additive map D1 D2 generator}
\end{equation}
}
{additive map D1 D2 generator}

\DefEq
{
\[
f:D_1\rightarrow D_2
\]
}
{map, division ring D1 D2}

\DefEq
{
\begin{equation}
f(a)=a^{\gi i}\ {}_{\gi i}G^{\gi l}\ f^{\gi k\gi j}_G
\ {}_{\gi k}\Vector q\ {}_{\gi l}\Vector e\ {}_{\gi j}\Vector q
\label{eq: additive map, division ring, F linear map G, standard representation}
\end{equation}
}
{additive map, division ring, F linear map G, standard representation}

\DefEq
{
\begin{equation}
a^{\gi i}\ {}_{\gi i}f^{\gi j}\ {}_{\gi j}\Vector q
=a^{\gi i}\ {}_{\gi i}G^{\gi l}\ f^{\gi{kr}}_G\ {}_{\gi k}\Vector q
\ {}_{\gi l}\Vector q\ {}_{\gi r}\Vector q
=a^{\gi i}\ {}_{\gi i}G^{\gi l}\ f^{\gi{kr}}_G\ {}_{\gi{kl}}B^{\gi p}
\ {}_{\gi{pr}}B^{\gi j}\ {}_{\gi j}\Vector q
\label{eq: additive multiplicative map over field, G, division ring, 1}
\end{equation}
}
{additive multiplicative map over field, G, division ring, 1}

\DefEq
{
\begin{equation}
f:D_1\rightarrow D_2
\label{eq: additive map D1 D2}
\end{equation}
}
{additive map D1 D2}

\DefEq
{
\[
h(x+y)=g(f(x+y))=g(f(x)+f(y))=g(f(x))+g(f(y))=h(x)+h(y)
\]
}
{product of additive map, D D D, 1}

\DefEq
{
\[
h(ax)=g(f(ax))=g(af(x))=ag(f(x))=ah(x)
\]
}
{product of additive map, D D D, 2}

\DefEq
{
\begin{align}
f&:D\rightarrow D
&f(x)&=\pC{s}{0}f\ x\ \pC{s}{1}f
\label{eq: D linear map, division ring, f}
\\
&&&=f^{\gi{ij}}\ {}_{\gi i}\Vector e\ x\ {}_{\gi j}\Vector e
\label{eq: D linear map, division ring, standard representation, f}
\\
g&:D\rightarrow D
&g(x)&=\pC{t}{0}g\ x\ \pC{t}{1}g
\label{eq: D linear map, division ring, g}
\\
&&&=g^{\gi{ij}}\ {}_{\gi i}\Vector e\ x\ {}_{\gi j}\Vector e
\label{eq: D linear map, division ring, standard representation, g}
\end{align}
}
{product of additive map, D D D}

\DefEq
{
\begin{align}
\pC{ts}{0}h&=\pC{t}{0}g\ \pC{s}{0}f
\label{eq: D linear map, division ring, h gf, 0}
\\
\pC{ts}{1}h&=\pC{s}{1}f\ \pC{t}{1}g
\label{eq: D linear map, division ring, h gf, 1}
\\
h^{\gi{pr}}&=g^{\gi{ij}}\ f^{\gi{kl}}
\ {}_{\gi{ik}}B^{\gi p}\ {}_{\gi{lj}}B^{\gi r}
\label{eq: D linear map, division ring, standard representation, h gf}
\end{align}
}
{D linear map, division ring, h gf}

\DefEq
{
\begin{align}
h(x)&=\pC{ts}{0}h\ x\ \pC{ts}{1}h
\label{eq: D linear map, division ring, h2}
\\
&=h^{\gi{pr}}\ {}_{\gi p}\Vector e\ x\ {}_{\gi r}\Vector e
\label{eq: D linear map, division ring, standard representation, h2}
\end{align}
}
{D linear map, division ring, h2}

\DefEq
{
\begin{equation}
f:D^n\rightarrow D,
d=f(d_1,...,d_n)
\label{eq: polylinear map, division ring}
\end{equation}
}
{polylinear map, division ring}

\DefEq
{
\begin{equation}
f(a_1,...,a_n)=a^{\gi{i_1}}_1...a^{\gi{i_n}}_n\ {}_{\gi{i_1...i_n}}f
\label{eq: polylinear map over field, division ring}
\end{equation}
}
{polylinear map over field, division ring}

\DefEq
{
\begin{equation}
{}_{\gi i}\Vector e'={}_{\gi i}A^{\gi j}\ {}_{\gi j}\Vector e
\label{eq: polylinear map over field, division ring, change basis}
\end{equation}
}
{polylinear map over field, division ring, change basis}

\DefEq
{
\begin{align}
{}_{\gi{i_1...i_n}}f'&=f({}_{\gi{i_1}}\Vector e',...,{}_{\gi{i_n}}\Vector e')
\nonumber
\\
&=f({}_{\gi{i_1}}A^{\gi{j_1}}\ {}_{\gi{j_1}}\Vector e,...,
{}_{\gi{i_n}}A^{\gi{j_n}}\ {}_{\gi{j_n}}\Vector e')
\label{eq: polylinear map over field, division ring, change coordinates}
\\
&={}_{\gi{i_1}}A^{\gi{j_1}}\ ...
\ {}_{\gi{i_n}}A^{\gi{j_n}}\ f({}_{\gi{j_1}}\Vector e,...,{}_{\gi{j_n}}\Vector e)
\nonumber
\\
&={}_{\gi{i_1}}A^{\gi{j_1}}\ ...\ {}_{\gi{i_n}}A^{\gi{j_n}}\ {}_{\gi{j_1...j_n}}f
\nonumber
\end{align}
}
{polylinear map over field, division ring, change coordinates}

\DefEq
{
\begin{equation}
f(\Vector a_1,...,\Vector a_n)=a^{\gi{i_1}}_1...a^{\gi{i_n}}_n\ {}_{\gi{i_1...i_n}}f^{\gi p}
\ {}_{\gi p}\Vector e
\label{eq: polylinear map over field, division ring, 1a}
\end{equation}
}
{polylinear map over field, division ring, 1a}

\DefEq
{
\begin{align}
a_j&=a^{\gi i}_j\ {}_{\gi i}\Vector e
\nonumber
\\
{}_{\gi{i_1...i_n}}f&=f({}_{\gi{i_1}}\Vector e,...,{}_{\gi{i_n}}\Vector e)
\label{eq: polylinear map over field, coordinates, division ring}
\end{align}
}
{polylinear map over field, coordinates, division ring}

\DefEq
{
\[
f(a_1,...,a_n)=
f(a^{\gi{i_1}}_1\ {}_{\gi{i_1}}\Vector e,...,a^{\gi{i_n}}_n\ {}_{\gi{i_n}}\Vector e)
=a^{\gi{i_1}}_1...a^{\gi{i_n}}_nf({}_{\gi{i_1}}e,...,{}_{\gi{i_n}}\Vector e)
\]
}
{polylinear map over field, 1, division ring}

\DefEq
{
\begin{equation}
{}_{\gi{i_1,...,i_n}}f={}_{\sigma(\gi{i_1}),...,\sigma(\gi{i_n})}f
\label{eq: polylinear map symmetric, division ring}
\end{equation}
}
{polylinear map symmetric, division ring}

\DefEq
{
\begin{align*}
a^{\gi{i_1}}_1\ ...\ a^{\gi{i_n}}_n\ {}_{\gi{i_1...i_n}}f
=&f(a_1,...,a_n)
\\
=&|\sigma|f(\sigma(a_1),...,\sigma(a_n))
\\
=&a_1^{\gi{i_1}}\ ...\ a_n^{\gi{i_n}}\ |\sigma|\ {}_{\sigma(\gi{i_1})...\sigma(\gi{i_n})}f
\end{align*}
}
{polylinear map skew symmetric, 1, division ring}

\DefEq
{
\begin{align}
{}_{\gi{j_1...j_n}}f
=&\pC{t}{}f^{\gi{i_0...i_n}}
\ {}_{\gi{i_0}\sigma_t(\gi{j_1})}B^{\gi{k_1}}
\ {}_{\gi{k_1i_1}}B^{\gi{l_1}}
\ ...\ {}_{\gi{l_{n-1}}\sigma_t(\gi{j_n})}B^{\gi{k_n}}
\ {}_{\gi{k_ni_n}}B^{\gi{l_n}}\ {}_{\gi{l_n}}\Vector e
\label{eq: coordinates of polylinear map, division ring over field}
\\
{}_{\gi{j_1...j_n}}f^{\gi p}
=&\pC{t}{}f^{\gi{i_0...i_n}}
\ {}_{\gi{i_0}\sigma_t(\gi{j_1})}B^{\gi{k_1}}
\ {}_{\gi{k_1i_1}}B^{\gi{l_1}}
\ ...\ {}_{\gi{l_{n-1}}\sigma_t(\gi{j_n})}B^{\gi{k_n}}
\ {}_{\gi{k_ni_n}}B^{\gi p}
\label{eq: coordinates of polylinear map, division ring over field, 1}
\end{align}
}
{coordinates of polylinear map, division ring over field}

\DefEq
{
\[
d_i=d_i^{\gi{j_i}}\ {}_{\gi{j_i}}\Vector e
\]
}
{coordinates of polylinear map, 1, division ring over field}

\DefEq
{
\begin{align}
f(d_1,...,d_n)
=&\pC{t}{}f^{\gi{i_0...i_n}}
\ {}_{\gi{i_0}}\Vector e\ \sigma_t(d_1^{\gi{j_1}}\ {}_{\gi{j_1}}\Vector e)
\ {}_{\gi{i_1}}\Vector e\ ...\ \sigma_t(d_n^{\gi{j_n}}
\ {}_{\gi{j_n}}\Vector e)\ {}_{\gi{i_n}}\Vector e
\nonumber\\
=&d_1^{\gi{j_1}}\ ...d_n^{\gi{j_n}}\ \pC{t}{}f^{\gi{i_0...i_n}}
\ {}_{\gi{i_0}}\Vector e\ \sigma_t({}_{\gi{j_1}}\Vector e)
\ {}_{\gi{i_1}}\Vector e\ ...\ \sigma_t({}_{\gi{j_n}}\Vector e)\ {}_{\gi{i_n}}\Vector e
\label{eq: polylinear map, division ring, standard representation, 1a}
\\
=&d_1^{\gi{j_1}}\ ...d_n^{\gi{j_n}}\ \pC{t}{}f^{\gi{i_0...i_n}}
\ {}_{\gi{i_0}\sigma_t(\gi{j_1})}B^{\gi{k_1}}
\ {}_{\gi{k_1i_1}}B^{\gi{l_1}}
\nonumber\\
&...\ {}_{\gi{l_{n-1}}\sigma_t(\gi{j_n})}B^{\gi{k_n}}
\ {}_{\gi{k_ni_n}}B^{\gi{l_n}}\ {}_{\gi{l_n}}\Vector e
\nonumber
\end{align}
}
{polylinear map, division ring, standard representation, 1a}

\DefEq
{
\begin{equation}
d=\pC{s}{0}f^{n\gi{j_1}}\ {}_{\gi{j_1}}\Vector e\ \sigma_s(d_1)
\ \pC{s}{1}f^{n\gi{j_2}}\ {}_{\gi{j_2}}\Vector e\ ...
\ \sigma_s(d_n)\ \pC{s}{n}f^{n\gi{j_n}}\ {}_{\gi{j_n}}\Vector e
\label{eq: polylinear map, division ring, standard representation, 1}
\end{equation}
}
{polylinear map, division ring, standard representation, 1}

\DefEq
{
\begin{equation}
\pC{t}{}f^{\gi{j_0...j_n}}=
\pC{s}{0}f^{n\gi{j_1}}\ ...\pC{s}{n}f^{n\gi{j_n}}
\label{eq: polylinear map, division ring, standard representation, 2}
\end{equation}
}
{polylinear map, division ring, standard representation, 2}

\DefEq
{
\begin{equation}
\pC{s}{p}f^n=\pC{s}{p}f^{n\gi i}\ {}_{\gi i}\Vector e
\label{eq: polylinear map, division ring, components extention}
\end{equation}
}
{polylinear map, division ring, components extention}

\DefEq
{
\begin{equation}
{}_{\gi{i_1,...,i_n}}f=|\sigma|\ {}_{\sigma(\gi{i_1}),...,\sigma(\gi{i_n})}f
\label{eq: polylinear map skew symmetric, division ring}
\end{equation}
}
{polylinear map skew symmetric, division ring}

\DefEq
{
\begin{align*}
a^{\gi{i_1}}_1\ ...\ a^{\gi{i_n}}_n\ {}_{\gi{i_1...i_n}}f
=&f(a_1,...,a_n)
\\
=&f(\sigma(a_1),...,\sigma(a_n))
\\
=&a_1^{\gi{i_1}}\ ...\ a_n^{\gi{i_n}}\ {}_{\sigma(\gi{i_1})...\sigma(\gi{i_n})}f
\end{align*}
}
{polylinear map symmetric, 1, division ring}

\DefEq
{
\begin{equation}
\label{eq: polylinear map, division ring, canonical morphism}
d=\pC{s}{0}f^n\ \sigma_s(d_1)
\ \pC{s}{1}f^n\ ...\ \sigma_s(d_n)\ \pC{s}{n}f^n
\end{equation}
}
{polylinear map, division ring, canonical morphism}

\DefEq
{
\[
\sigma_s=
\begin{pmatrix}
d_1&...&d_n
\\
\sigma_s(d_1)&...&\sigma_s(d_n)
\end{pmatrix}
\]
}
{transposition of set of variables, division ring}

\DefEq
{
\[
\pC{s}{p}f^1=\pC{s}{p}f
\]
}
{polylinear map, division ring, 1, canonical morphism}

\DefEq
{
\[
\xymatrix{
D^k\ar[rr]^f\ar@{=>}[drr]_{g(d_k)}
& & D
\\
& &
\\
&&D^{k-1}\ar[uu]_h &
}
\]
\[
d=f(d_1,...,d_k)=g(d_k)(d_1,...,d_{k-1})
\]
}
{polylinear map, induction on n, 1, division ring}

\DefEq
{
\[
d=
\pC{t}{0}h^{k-1}\ \sigma_t(d_1)
\ \pC{t}{1}h^{k-1}\ ...
\ \sigma_t(d_{k-1})\ \pC{t}{k-1}h^{k-1}
\]
}
{polylinear map, induction on n, 2, division ring}

\DefEq
{
\[
\pC{t}{0}h^{k-1}
=\pC{tr}{0}g
\ d_k\ \pC{tr}{1}g
\]
}
{polylinear map, induction on n, 3, division ring}

\DefEq
{
\[
\sigma_s=\sigma_{tr}=
\left(
\begin{array}{cccc}
d_k&d_1&...&d_{k-1}
\\
d_k&\sigma_t(d_1)&...&\sigma_t(d_{k-1})
\end{array}
\right)
\]
}
{polylinear map, induction on n, 4, division ring}

\DefEq
{
\[
\pC{tr}{q+1}f^k=\pC{t}{q}h^{k-1}
\]
}
{polylinear map, induction on n, 5, division ring}

\DefEq
{
\[
\pC{tr}{q}f^k=\pC{tr}{q}g
\]
}
{polylinear map, induction on n, 6, division ring}

\DefEq
{
\begin{equation}
f(d_1,...,d_n)=
\ShowSymbol{standard component of polylinear map, division ring}
\ {}_{\gi{i_0}}\Vector e\ \sigma_t(d_1)\ {}_{\gi{i_1}}\Vector e
\ ...\ \sigma_t(d_n)\ {}_{\gi{i_n}}\Vector e
\label{eq: polylinear map, division ring, standard representation}
\end{equation}
}
{polylinear map, division ring, standard representation}

%% file: Differential.English.tex
\def\texDifferential{}

\input{Differential.Eq}

\ifx\PrintBook\undefined
\else
	\chapter{Differentiable Maps}
	\label{chapter: Differentiable maps}
\fi
\section{Topological Division Ring}
\label{section: Topological Division Ring}

\begin{definition}
Division ring $D$ is called 
\AddIndex{topological division ring}{topological division ring}\footnote{I
made definition according to definition
from \citeBib{Pontryagin: Topological Group},
chapter 4}
if $D$ is topological space and the algebraic operations
defined in $D$ are continuous in the topological space $D$.
\qed
\end{definition}

According to definition, for arbitrary elements $a, b\in D$
and for arbitrary neighborhoods $W_{a-b}$ of the element $a-b$,
$W_{ab}$ of the element $ab$ there exists neighborhoods $W_a$
of the element $a$ and $W_b$ of the element $b$ such
that $W_a-W_b\subset W_{a-b}$,
$W_aW_b\subset W_{ab}$. For any $a\ne 0$ and for arbitrary
neighborhood $W_{a^{-1}}$ there exists neighborhood $W_a$ of the element $a$,
satisfying the condition $W_a^{-1}\subset W_{a^{-1}}$.

\begin{definition}
\label{definition: absolute value on division ring}
\AddIndex{Absolute value on division ring}
{absolute value on division ring} $D$\footnote{I
made definition according to definition
from \citeBib{Bourbaki: General Topology: Chapter 5 - 10},
IX, \S 3.2} is a map
\[d\in D\rightarrow |d|\in R\]
which satisfies the following axioms
\begin{itemize}
\item $|a|\ge 0$
\item $|a|=0$ if, and only if, $a=0$
\item $|ab|=|a|\ |b|$
\item $|a+b|\le|a|+|b|$
\end{itemize}

Division ring $D$, endowed with the structure defined by a given absolute value on
$D$, is called \AddIndex{valued division ring}{valued division ring}.
\qed
\end{definition}

Invariant distance on additive group of division ring $D$
\[d(a,b)=|a-b|\]
defines topology of metric space,
compatible with division ring structure of $D$.

\begin{definition}
\label{definition: limit of sequence, valued division ring}
Let $D$ be valued division ring.
Element $a\in D$ is called 
\AddIndex{limit of a sequence}{limit of sequence, valued division ring}
$\{a_n\}$
\symb{\lim_{n\rightarrow\infty}a_n}0{limit of sequence, valued division ring}
\[
a=\ShowSymbol{limit of sequence, valued division ring}
\]
if for every $\epsilon\in R$, $\epsilon>0$
there exists positive integer $n_0$ depending on $\epsilon$ and such,
that $|a_n-a|<\epsilon$ for every $n>n_0$.
\qed
\end{definition}

\ifx\EprintCalculus\undefined
\begin{theorem}
\label{theorem: limit of sequence, valued division ring, product on scalar}
Let $D$ be valued division ring of characteristic $0$ and let $d\in D$.
Let $a\in D$ be limit of a sequence $\{a_n\}$.
Then
\[
\lim_{n\rightarrow\infty}(a_nd)=ad
\]
\[
\lim_{n\rightarrow\infty}(da_n)=da
\]
\end{theorem}
\begin{proof}
Statement of the theorem is trivial, however I give this proof
for completeness sake. 
Since $a\in D$ is limit of the sequence $\{a_n\}$,
then according to definition \ref{definition: limit of sequence, valued division ring}
for given $\epsilon\in R$, $\epsilon>0$,
there exists positive integer $n_0$ such,
that $|a_n-a|<\epsilon/|d|$ for every $n>n_0$.
According to definition \ref{definition: absolute value on division ring}
the statement of theorem follows from inequations
\[
|a_nd-ad|=|(a_n-a)d|=|a_n-a||d|<\epsilon/|d||d|=\epsilon
\]
\[
|da_n-da|=|d(a_n-a)|=|d||a_n-a|<|d|\epsilon/|d|=\epsilon
\]
for any $n>n_0$.
\end{proof}
\fi

\begin{definition}
Let $D$ be valued division ring.
The sequence $\{a_n\}$, $a_n\in D$ is called 
\AddIndex{fundamental}{fundamental sequence, valued division ring}
or \AddIndex{Cauchy sequence}{Cauchy sequence, valued division ring},
if for every $\epsilon\in R$, $\epsilon>0$
there exists positive integer $n_0$ depending on $\epsilon$ and such,
that $|a_p-a_q|<\epsilon$ for every $p$, $q>n_0$.
\qed
\end{definition}

\begin{definition}
Valued division ring $D$ is called
\AddIndex{complete}{complete division ring}
if any fundamental sequence of elements
of division ring $D$ converges, i.e.
has limit in division ring $D$.
\qed
\end{definition}

\ifx\EprintCalculus\undefined
So far everything was good. Ring of characteristic $0$ contains ring of integers.
Division ring of characteristic $0$
contains field of rational numbers.
However this does not mean that complete division ring $D$ of characteristic $0$
contains field of real numbers.
For instance, absolute value $|x|=1$
determines discrete topology on valued
division ring and is not interesting for us,
because any fundamental sequence is
constant map.
\fi

Later on, speaking about valued division ring of characteristic $0$,
we will assume that homeomorphism of field of rational numbers $Q$
into division ring $D$ is defined.

\begin{theorem}
\label{theorem: division ring contains real number}
Complete division ring $D$ of characteristic $0$
contains as subfield an isomorphic image of the field $R$ of
real numbers.
It is customary to identify it with $R$.
\end{theorem}
\begin{proof}
Let us consider fundamental sequence of rational numbers $\{p_n\}$.
Let $p'$ be limit of this sequence in division ring $D$.
Let $p$ be limit of this sequence in field $R$.
Since immersion of field $Q$ into division ring $D$ is homeomorphism,
then we may identify $p'\in D$ and $p\in R$.
\end{proof}

\begin{theorem}
\label{theorem: division ring, real number}
Let $D$ be complete division ring of characteristic $0$ and let $d\in D$.
Then any real number $p\in R$ commute with $d$.
\end{theorem}
\begin{proof}
Let us represent real number $p\in R$ as
fundamental sequence of rational numbers $\{p_n\}$.
Statement of theorem follows from chain of equations
\[
pd=\lim_{n\rightarrow\infty}(p_nd)=\lim_{n\rightarrow\infty}(dp_n)=dp
\]
based on statement of theorem
\xRef{0812.4763}{theorem: limit of sequence, valued division ring, product on scalar}.
\end{proof}

\begin{theorem}
\label{theorem: real subring of center}
Let $D$ be complete division ring of characteristic $0$.
Then field of real numbers $R$ is subfield of center $Z(D)$ of division ring $D$.
\end{theorem}
\begin{proof}
Corollary of theorem \ref{theorem: division ring, real number}.
\end{proof}

\begin{definition}
Let $D$ be complete division ring of characteristic $0$.
Set of elements $d\in D$, $|d|=1$ is called
\AddIndex{unit sphere in division ring}{unit sphere in division ring} $D$.
\qed
\end{definition}

\begin{definition}
\label{definition: continuous function, division ring}
Let $D_1$ be complete division ring of characteristic $0$ with absolute value $|x|_1$.
Let $D_2$ be complete division ring of characteristic $0$ with absolute value $|x|_2$.
Function \[f:D_1\rightarrow D_2\]
is called \AddIndex{continuous}{continuous function, division ring}, if
for every as small as we please $\epsilon>0$
there exist such $\delta>0$, that
\[
|x'-x|_1<\delta
\]
implies
\[
|f(x')-f(x)|_2<\epsilon
\]
\qed
\end{definition}

\begin{theorem}
\label{theorem: continuous map of division ring}
Let $D$ be complete division ring of characteristic $0$.
Since index $s$ in expansion
\ifx\IntroCalculusBook\Defined
\eqref{eq: linear map, division ring, canonical morphism}
\else
\ifx\EprintCalculus\Defined
\eqref{eq: linear map, division ring, canonical morphism}
\else
\xEqRef{0701.238}{eq: linear map, division ring, canonical morphism}
\fi
\fi
of additive map \[f:D\rightarrow D\]
belongs to finite range,
then additive map $f$
is continuous.
\end{theorem}
\begin{proof}
Suppose $x'=x+a$. Then
\[
f(x')-f(x)=f(x+a)-f(x)=f(a)=\pC{s}{0}f\ a\ \pC{s}{1}f
\]
\[
|f(x')-f(x)|=|\pC{s}{0}f\ a\ \pC{s}{1}f|<(|\pC{s}{0}f|\ |\pC{s}{1}f|)|a|
\]
Let us denote $F=|\pC{s}{0}f|\ |\pC{s}{1}f|$. Then
\[
|f(x')-f(x)|<F|a|
\]
Suppose $\epsilon>0$ and let us assume $\displaystyle a=\frac{\epsilon}Fe$.
Then $\displaystyle\delta=|a|=\frac{\epsilon}F$.
According to definition \ref{definition: continuous function, division ring}
additive map $f$ is continuous.
\end{proof}

\ifx\EprintCalculus\undefined
Similarly, since index $s$ belongs to countable range,
then for continuity of additive map $f$ it is necessary that series
$|\pC{s}{0}f|\ |\pC{s}{1}f|$ converges.
Since index $s$ belongs to continues range,
then for continuity of additive map $f$ it is necessary that integral
$\int|\pC{s}{0}f|\ |\pC{s}{1}f|ds$ exist.
\fi

\begin{definition}
\label{definition: norm of map, division ring}
Let
\[f:D_1\rightarrow D_2\]
map of complete division ring $D_1$ of characteristic $0$ with absolute value $|x|_1$
into complete division ring $D_2$ of characteristic $0$ with absolute value $|y|_2$.
Value \symb{\|f\|}0{norm of map, division ring}
\begin{equation}
\ShowSymbol{norm of map, division ring}=
\text{sup}\frac{|f(x)|_2}{|x|_1}
\label{eq: norm of map, division ring}
\end{equation}
is called
\AddIndex{norm of map $f$}
{norm of map, division ring}.
\qed
\end{definition}

\begin{theorem}
\label{theorem: norm of multiplicative map, division ring}
Let $D_1$ be complete division ring of characteristic $0$ with absolute value $|x|_1$.
Let $D_2$ be complete division ring of characteristic $0$ with absolute value $|x|_2$.
Let
\[f:D_1\rightarrow D_2\]
be map which is multiplicative over field $R$.
Then
\begin{equation}
\|f\|=\text{sup}\{|f(x)|_2:|x|_1=1\}
\label{eq: norm of multiplicative map, division ring}
\end{equation}
\end{theorem}
\begin{proof}
According to definition
\ifx\IntroCalculusBook\Defined
\ref{definition: map multiplicative over commutative ring, ring}
\else
\ifx\EprintCalculus\Defined
\ref{definition: map multiplicative over commutative ring, ring}
\else
\xRef{0701.238}{definition: map multiplicative over commutative ring, ring}
\fi
\fi
\[
\frac{|f(x)|_2}{|x|_1}=\frac{|f(rx)|_2}{|rx|_1}
\]
Assuming $\displaystyle r=\frac 1{|x|_1}$, we get
\begin{equation}
\frac{|f(x)|_2}{|x|_1}=\left|f\left(\frac x{|x|_1}\right)\right|_2
\label{eq: norm of multiplicative map, division ring, 1}
\end{equation}
Equation \eqref{eq: norm of multiplicative map, division ring}
follows from equations \eqref{eq: norm of multiplicative map, division ring, 1}
and \eqref{eq: norm of map, division ring}.
\end{proof}

\begin{theorem}
\label{theorem: continuity of additive map, finite norm, division ring}
Let
\[f:D_1\rightarrow D_2\]
additive map of complete division ring $D_1$
into complete division ring $D_2$.
Since $\|f\|<\infty$, then map $f$ is continuous.
\end{theorem}
\begin{proof}
Because map $f$ is additive, then
according to definition \ref{definition: norm of map, division ring}
\[
|f(x)-f(y)|_2=|f(x-y)|_2\le \|f\|\ |x-y|_1
\]
Let us assume arbitrary $\epsilon>0$. Assume $\displaystyle\delta=\frac\epsilon{\|f\|}$.
Then
\[
|f(x)-f(y)|_2\le \|f\|\ \delta=\epsilon
\]
follows from inequation
\[
|x-y|_1<\delta
\]
According to definition \ref{definition: continuous function, division ring}
map $f$ is continuous.
\end{proof}

\ifx\texFuture\Defined
\begin{theorem}
\label{theorem: continuity of map, norm, division ring}
Let
\[f:D_1\rightarrow D_2\]
map of complete division ring $D_1$
в complete division ring $D_2$.
map $f$ is continuous, if $\|f\|<\infty$.
\end{theorem}
\begin{proof}
Since map $f$ is additive, then
according to definition \ref{definition: norm of map, division ring}
\[
|f(x)-f(y)|_2=|f(x-y)|_2\le \|f\|\ |x-y|_1
\]
Let $\epsilon>0$ be arbitrary. Assume $\displaystyle\delta=\frac\epsilon{\|f\|}$.
Then from inequation
\[
|x-y|_1<\delta
\]
follows
\[
|f(x)-f(y)|_2\le \|f\|\ \delta=\epsilon
\]
According to definition \ref{definition: continuous function, division ring}
map $f$ is continuous.
\end{proof}
\fi

\begin{theorem}
\label{theorem: projective map in 0, division ring}
Let $D$ be complete division ring of characteristic $0$.
Either continuous map $f$ of division ring which is projective over field $P$,
does not depend on direction over field $P$,
or value $f(0)$ is not defined.
\end{theorem}
\begin{proof}
According to definition
\ifx\IntroCalculusBook\Defined
\ref{definition: map projective over commutative ring, ring},
\else
\xRef{0701.238}{definition: map projective over commutative ring, ring},
\fi
map $f$
is constant on direction $Pa$. Since $0\in Pa$,
then we may assume
\[
f(0)=f(a)
\]
based on continuity.
However this leads to uncertainty of value
of map $f$ in direction $0$,
when map $f$ has different values for different directions $a$.
\end{proof}

If projective over field $R$ map $f$ is continuous,
then we say that function $f$ is \AddIndex{continuous
in direction over field $R$}
{projective function is continuous in direction over field R, division ring}.
Since for any $a\in D$, $a\ne 0$ we may choose
$a_1=|a|^{-1}a$, $f(a_1)=f(a)$, then it is possible to make definition more accurate.

\begin{definition}
\label{definition: projective function is continuous in direction over field R, division ring}
Let $D$ be complete division ring of characteristic $0$.
Projective over field $R$ function $f$ is continuous
in direction over field $R$, if
for every as small as we please $\epsilon>0$
there exists such $\delta>0$, that
\[
|x'-x|_1<\delta\ \ \ |x'|_1=|x|_1=1
\]
implies
\[
|f(x')-f(x)|_2<\epsilon
\]
\qed
\end{definition}

\begin{theorem}
\label{theorem: projective function is continuous in direction over field R, division ring}
Let $D$ be complete division ring of characteristic $0$.
Projective over field $R$ function $f$ is continuous
in direction over field $R$ iif
this function is continuous on unit sphere of division ring $D$.
\end{theorem}
\begin{proof}
Corollary of definitions \ref{definition: continuous function, division ring},
\ifx\IntroCalculusBook\Defined
\ref{definition: map projective over commutative ring, ring},
\else
\xRef{0701.238}{definition: map projective over commutative ring, ring},
\fi
\ref{definition: projective function is continuous in direction over field R, division ring}.
\end{proof}

\section{Differentiable Map of Division Ring}
\label{section: Differentiable maps of Division Ring}

\ifx\IntroCalculusBook\Defined
\def\texFrechet{}
\begin{definition}
\label{definition: differentiable functions, division ring}
Let $D$ be valued division ring.\footnote{I
made definition according to definition
\citeBib{Lebedev Vorovich}-3.1.1,
page 177.}
Function \[f:D\rightarrow D\] is called
\AddIndex{\Ds differentiable in the Fr\'echet sense}
{function Dstar differentiable in Frechet sense, division ring}
on the set $U\subset D$\footnote{Like in remark
\xRef{0701.238}{remark: drc linear map}
we may define \Ds differentiability of map \[f:S\rightarrow D\]
of valued division ring $S$
into valued division ring $D$ according to rule
\[
f(x+h)-f(x)
=F(h)\frac{d f(x)}{d_*x}
+o(h)
\]
where \[F:S\rightarrow D\] is homomorphism of division rings.
However, based on isomorphism theorem, we may confine ourselves by case of
maps of division ring $D$ into $D$.},
if at every point $x\in U$ the increment of the function $f$ can be represented as
\symb{\frac{d f(x)}{d_*x}}0{Frechet Dstar derivative of map, division ring}
\begin{equation}
\label{eq: Dstar derivative, division ring}
f(x+h)-f(x)
=h\ShowSymbol{Frechet Dstar derivative of map, division ring}
+o(h)
\end{equation}
where $o$ is such continuous map
\[o:D\rightarrow D\]
that
\begin{equation}
\label{eq: infinitesimal, division ring}
\lim_{h\rightarrow 0}\frac{|o(h)|}{|h|}=0
\end{equation}
\qed
\end{definition}

According to definition \ref{definition: differentiable functions, division ring},
\AddIndex{the Fr\'echet \Ds derivative of map $f$ at point $x$}
{Frechet Dstar derivative of map, division ring}
\[
\ShowSymbol{Frechet Dstar derivative of map, division ring}\in {}_*\mathcal L(D,D)
\]
generates a homomorphism
\[
\Delta f
=\Delta x\frac{d f(x)}{d_*x}
\]
which maps increment of argument into increment of function.

If we multiply both sides of equation
\eqref{eq: Dstar derivative, division ring}
by $h^{-1}$, then we get equation
\begin{equation}
\label{eq: Dstar derivative, division ring, 1}
h^{-1}(f(x+h)-f(x))
=\frac{d f(x)}{d_*x}
+h^{-1}o(h)
\end{equation}
Alternative definition of the Fr\'echet \Ds derivative
\begin{equation}
\label{eq: Dstar derivative, division ring, 2}
\frac{d f(x)}{d_*x}
=\lim_{h\rightarrow 0}(h^{-1}(f(x+h)-f(x)))
\end{equation}
follows from equations \eqref{eq: Dstar derivative, division ring, 1} and
\eqref{eq: infinitesimal, division ring}.

In contrast to derivative of map over field $F$
the Fr\'echet \Ds derivative of map of division ring $D$ is not
linear map. Indeed, the Fr\'echet \Ds derivative is \sD
linear map; however, in general, it is not \Ds
linear map. In fact,
\[
f(x+h)a-f(x)a=h\frac{df}{d_*x}a+o(h)
\]
\begin{equation}
\label{eq: Dstar derivative, division ring, 3}
af(x+h)-af(x)=ah\frac{df}{d_*x}+o(h)
\end{equation}
However, in general,
\[
ah\frac{df}{d_*x}\ne ha\frac{df}{d_*x}
\]

Similar problem arises when we differentiate product of maps.
According to definition \ref{definition: differentiable functions, division ring}
\begin{equation}
f(x+h)g(x+h)-f(x)g(x)
=h\frac{df(x)g(x)}{d_*x}+o(h)
\label{eq: product of maps, derivative, 0}
\end{equation}
It is possible to represent expression in left side as
\begin{align}
&f(x+h)g(x+h)-f(x)g(x)\nonumber\\
\label{eq: product of maps, derivative, 1}
=&f(x+h)g(x+h)-f(x)g(x+h)+f(x)g(x+h)-f(x)g(x)\\
=&(f(x+h)-f(x))g(x+h)+f(x)(g(x+h)-g(x))\nonumber
\end{align}
According to definition \ref{definition: differentiable functions, division ring},
it is possible to represent \eqref{eq: product of maps, derivative, 1} as
\begin{align}
&f(x+h)g(x+h)-f(x)g(x)\nonumber\\
\label{eq: product of maps, derivative, 2}
=&\left(h\frac{df(x)}{d_*x}+o(h)\right)
\left(g(x)+h\frac{dg(x)}{d_*x}+o(h)\right)\\
+&f(x)\left(h\frac{dg(x)}{d_*x}+o(h)\right)\nonumber
\end{align}
Since, in general,
\[
f(x)h\frac{dg}{d_*x}\ne hf(x)\frac{dg}{d_*x}
\]
then it is natural to anticipate that map $f(x)g(x)$
is not \Ds differentiable in the Fr\'echet sense.

Thus definition the Fr\'echet \Ds derivative is extremely restricted
and does not satisfy standard definition of differentiation.
To find solution of problem of differentiation,
let us consider this problem on the other hand.
\fi

\begin{example}
\label{example: differential dxx}
Let us consider increment of map $f(x)=x^2$.
\begin{align*}
f(x+h)-f(x)
&=(x+h)^2-x^2
\\
&=xh+hx+h^2
\\
&=xh+hx+o(h)
\end{align*}
As can be easily seen, the component of the increment of the function $f(x)=x^2$
that is linearly dependent on the increment of the argument, is of the form
\[
xh+hx
\]
Since product is non commutative,
we cannot represent increment of map $f(x+h)-f(x)$
as $Ah$ or $hA$
where $A$ does not depend on $h$.
It results in the unpredictable behavior of
the increment of the function $f(x)=x^2$ when
the increment of the argument converges to $0$.
However, since infinitesimal $h$ is infinitesimal like
$h=ta$, $a\in D$, $t\in R$, $t\rightarrow 0$, the answer becomes
more definite
\[
(xa+ax)t
\]
\qed
\end{example}

\begin{definition}
\label{definition: function differentiable in Gateaux sense, division ring}
Let $D$ be complete division ring of characteristic $0$.\footnote{I
made definition according to definition
\citeBib{Lebedev Vorovich}-3.1.2,
page 177.}
The function \[f:D\rightarrow D\]
is called \AddIndex{differentiable in the G\^ateaux sense}
{function differentiable in Gateaux sense, division ring} on the set $U\subset D$,
if at every point $x\in U$
the increment of the function $f$ can be represented as
\symb{\partial f(x)}
0{Gateaux derivative of map, division ring}
\symb{\frac{\partial f(x)}{\partial x}}
0{Gateaux derivative of map, fraction, division ring}
\begin{equation}
\label{eq: Gateaux differential of map, division ring}
f(x+a)-f(x)
=\ShowSymbol{Gateaux derivative of map, division ring}(a)
+o(a)
=\ShowSymbol{Gateaux derivative of map, fraction, division ring}(a)
+o(a)
\end{equation}
where
\AddIndex{the G\^ateaux derivative
$\ShowSymbol{Gateaux derivative of map, division ring}$
of map $f$}
{Gateaux derivative of map, division ring}
is linear map of increment $a$ and
$o:D\rightarrow D$ is such continuous map that
\[
\lim_{a\rightarrow 0}\frac{|o(a)|}{|a|}=0
\]
\qed
\end{definition}

\begin{remark}
\label{remark: differential L(D,D)}
According to definition
\ref{definition: function differentiable in Gateaux sense, division ring}
for given $x$, the G\^ateaux derivative $\partial f(x)\in\mathcal L(D;D)$.
Therefore, the G\^ateaux derivative of map $f$ is map
\ShowEq{differential L(D,D)}
Expressions $\partial f(x)$ and $\displaystyle\frac{\partial f(x)}{\partial x}$
are different notations for the same function.
We will use notation $\displaystyle\frac{\partial f(x)}{\partial x}$
to underline that this is  the G\^ateaux derivative with respect to variable $x$.
\qed
\end{remark}

\begin{theorem}
\label{theorem: function differentiable in Gateaux sense, division ring}
\symb{\partial f(x)(a)}
0{Gateaux differential of map, division ring}
It is possible to represent \AddIndex{the G\^ateaux differential
$\ShowSymbol{Gateaux differential of map, division ring}$
of map $f$}
{Gateaux differential of map, division ring}
as
\begin{equation}
\ShowSymbol{Gateaux differential of map, division ring}
=\frac{\pC{s}{0}\partial f(x)}{\partial x}
\ a
\ \frac{\pC{s}{1}\partial f(x)}{\partial x}
\label{eq: Gateaux differential, representation}
\end{equation}
\end{theorem}
\begin{proof}
Corollary of definitions \ref{definition: function differentiable in Gateaux sense, division ring}
and theorem
\ifx\IntroCalculusBook\Defined
\ref{theorem: linear map, division ring, canonical morphism}.
\else
\ifx\EprintCalculus\Defined
\ref{theorem: linear map, division ring, canonical morphism}.
\else
\xRef{0701.238}{theorem: linear map, division ring, canonical morphism}.
\fi
\fi
\end{proof}

\begin{definition}
\begin{sloppypar}
Expression
\symb{\frac{\pC{s}{p}\partial f(x)}{\partial x}}0
{component of Gateaux derivative of map, division ring}
$\displaystyle\ShowSymbol{component of Gateaux derivative of map, division ring}$,
$p=0$, $1$,
is called
\AddIndex{component of the G\^ateaux derivative of map $f(x)$}
{component of Gateaux derivative of map, division ring}.
\qed
\end{sloppypar}
\end{definition}

\begin{theorem}
\label{theorem: Gateaux derivative is multiplicative over field R, division ring}
Let $D$ be division ring of characteristic $0$.
The G\^ateaux derivative
of function \[f:D\rightarrow D\]
is multiplicative over field $R$.
\end{theorem}
\begin{proof}
Corollary of theorems
\ref{theorem: real subring of center},
\ifx\IntroCalculusBook\Defined
\ref{theorem: additive map, division ring, center},
\else
\ifx\EprintCalculus\Defined
\ref{theorem: additive map, division ring, center},
\else
\xRef{0701.238}{theorem: additive map, division ring, center},
\fi
\fi
and
definition \ref{definition: function differentiable in Gateaux sense, division ring}.
\end{proof}

From theorem \ref{theorem: Gateaux derivative is multiplicative over field R, division ring}
it follows
\begin{equation}
\label{eq: Gateaux differential is multiplicative over field R, division ring}
\partial f(x)(ra)
=
r\partial f(x)(a)
\end{equation}
for any $r\in R$, $r\ne 0$ and $a\in D$, $a\ne 0$.
Combining equation \eqref{eq: Gateaux differential is multiplicative over field R, division ring}
and definition
\ref{definition: function differentiable in Gateaux sense, division ring},
we get known definition of the G\^ateaux
differential
\begin{equation}
\partial f(x)(a)=\lim_{t\rightarrow 0,\ t\in R}(t^{-1}(f(x+ta)-f(x)))
\label{eq: Gateaux differential, representation, 2}
\end{equation}
Definitions of the G\^ateaux derivative \eqref{eq: Gateaux differential of map, division ring}
and \eqref{eq: Gateaux differential, representation, 2}
are equivalent. Using this equivalence we tell
that map $f$
is called differentiable in the G\^ateaux sense on the set $U\subset D$,
if at every point $x\in U$
the increment of the function $f$ can be represented as
\begin{equation}
\label{eq: Gateaux differential of map, t, division ring}
f(x+ta)-f(x)
=t\partial f(x)(a)
+o(t)
\end{equation}
where
$o:R\rightarrow D$ is such continuous map that
\[
\lim_{t\rightarrow 0}\frac{|o(t)|}{|t|}=0
\]

Since infinitesimal $ta$ is differential $dx$,
then equation
\eqref{eq: Gateaux differential, representation}
gets form
\begin{equation}
\partial f(x)(dx)
=\frac{\pC{s}{0}\partial f(x)}{\partial x}
dx
\frac{\pC{s}{1}\partial f(x)}{\partial x}
\label{eq: Gateaux differential, form}
\end{equation}

\begin{theorem}
\label{theorem: Gateaux differential, standard form, division ring}
Let $D$ be division ring of characteristic $0$.
Let $\Basis e$ be basis of division ring $D$ over center $Z(D)$ of division ring $D$.
\symb{\StandPartial{f(x)}{x}{ij}}0
{standard component of Gateaux differential, division ring}
\AddIndex{Standard representation of the G\^ateaux differential
\eqref{eq: Gateaux differential, representation}
of map\[f:D\rightarrow D\]}
{Gateaux differential, standard representation, division ring}
has form
\begin{equation}
\partial f(x)(a)
=\ShowSymbol{standard component of Gateaux differential, division ring}
\ {}_{\gi i}\Vector e\ a\ {}_{\gi j}\Vector e
\label{eq: Gateaux differential, division ring, standard representation}
\end{equation}
Expression
$\displaystyle\ShowSymbol{standard component of Gateaux differential, division ring}$
in equation \eqref{eq: Gateaux differential, division ring, standard representation}
is called
\AddIndex{standard component of the G\^ateaux differential of map $f$}
{standard component of Gateaux differential, division ring}.
\end{theorem}
\begin{proof}
Statement of theorem is corollary of theorem
\ifx\IntroCalculusBook\Defined
\ref{theorem: linear map, standard form, division ring}.
\else
\ifx\EprintCalculus\Defined
\ref{theorem: linear map, standard form, division ring}.
\else
\xRef{0701.238}{theorem: linear map, standard form, division ring}.
\fi
\fi
\end{proof}

\begin{theorem}
Let $D$ be division ring of characteristic $0$.
Let $\Basis e$ be basis of division ring $D$ over center $Z(D)$ of division ring $D$.
Then it is possible to represent the G\^ateaux differential of map
\[
f:D\rightarrow D
\]
as
\begin{equation}
\partial f(x)(a)=a^{\gi i}\ \frac{\partial f^{\gi j}}{\partial x^{\gi i}}
\ {}_{\gi j}\Vector e
\label{eq: Gateaux differential and Jacobian, division ring}
\end{equation}
where $a\in D$ has expansion
\[
\begin{array}{ccc}
a=a^{\gi i}\ {}_{\gi i}\Vector e&&a^{\gi i}\in F
\end{array}
\]
relative to basis $\Basis e$ and Jacobian of map $f$ has form
\begin{equation}
\label{eq: standard components and Jacobian, division ring}
\frac{\partial f^{\gi j}}{\partial x^{\gi i}}=
\StandPartial{f(x)}{x}{kr}\ {}_{\gi{ki}}B^{\gi p}\ {}_{\gi{pr}}B^{\gi j}
\end{equation}
\end{theorem}
\begin{proof}
Statement of theorem is corollary of theorem
\ifx\IntroCalculusBook\Defined
\ref{theorem: linear map over field, division ring}.
\else
\ifx\EprintCalculus\Defined
\ref{theorem: linear map over field, division ring}.
\else
\xRef{0701.238}{theorem: linear map over field, division ring}.
\fi
\fi
\end{proof}

\ifx\EprintCalculus\undefined
To find structure similar to derivative in commutative calculus.
we need extract derivative from differential
of map. For this purpose we need
take the increment of argument outside brackets.\footnote{It is possible when all
$\displaystyle \frac{\pC{s}{0}\partial f(x)}{\partial x}=e$
or all
$\displaystyle \frac{\pC{s}{1}\partial f(x)}{\partial x}=e$.
Therefore definition
\citeBib{Lebedev Vorovich}-3.1.2,
page 177, leads us either to definition of
the Fr\'echet \Ds derivative,
or to definition of
the Fr\'echet \sD derivative.}
It is possible to take $a$ outside brackets
using equations
\[
\frac{\pC{s}{0}\partial f(x)}{\partial x}a
\frac{\pC{s}{1}\partial f(x)}{\partial x}
=\frac{\pC{s}{0}\partial f(x)}{\partial x}a
\frac{\pC{s}{1}\partial f(x)}{\partial x}a^{-1}a
\]
\[
\frac{\pC{s}{0}\partial f(x)}{\partial x}a
\frac{\pC{s}{1}\partial f(x)}{\partial x}
=aa^{-1}\frac{\pC{s}{0}\partial f(x)}{\partial x}a
\frac{\pC{s}{1}\partial f(x)}{\partial x}
\]
\fi

\begin{definition}
\label{definition: Gateaux Dstar derivative, division ring}
Let $D$ be complete division ring of characteristic $0$
and $a\in D$.
\symb{\frac{\partial f(x)(a)}{\partial_* x}}
0{Gateaux Dstar derivative of map, division ring}
We define \AddIndex{the G\^ateaux \Ds derivative
$\displaystyle\ShowSymbol{Gateaux Dstar derivative of map, division ring}$
of map $f:D\rightarrow D$}
{Gateaux Dstar derivative of map, division ring}
using equation
\begin{equation}
\label{eq: Gateaux Dstar derivative, division ring}
\partial f(x)(a)
=a\ShowSymbol{Gateaux Dstar derivative of map, division ring}
\end{equation}
\symb{\frac{\partial f(x)(a)}{{}_*\partial x}}
0{Gateaux starD derivative of map, division ring}
We define \AddIndex{the G\^ateaux \sD derivative
$\displaystyle\ShowSymbol{Gateaux starD derivative of map, division ring}$
of map $f:D\rightarrow D$}
{Gateaux starD derivative of map, division ring}
using equation
\begin{equation}
\label{eq: Gateaux starD derivative, division ring}
\partial f(x)(a)
=\ShowSymbol{Gateaux starD derivative of map, division ring}a
\end{equation}
\qed
\end{definition}

Let us consider the basis ${}_1e=1$, ${}_2e=i$, ${}_3e=j$, ${}_4e=k$
of division ring of quaternions over real field.
From straightforward calculation, it follows that
standard \Ds representation of the G\^ateaux differential of map $x^2$
has form
\[
\partial x^2(a)=(x+x_1)a+x_2ai+x_3aj+x_4ak
\]

\ifx\EprintCalculus\undefined
For representation of the G\^ateaux \Ds derivative
we also use notation
\[
\frac{\partial f(x)(a)}{\partial_* x}=\frac{\partial}{\partial_* x} f(x)(a)=\partial_* f(x)(a)
\]

For representation of the G\^ateaux \sD derivative
we also use notation
\[
\frac{\partial f(x)(a)}{{}_*\partial x}=\frac{\partial}{{}_*\partial x} f(x)(a)={}_*\partial f(x)(a)
\]

Based on duality principle
\xRef{0701.238}{theorem: duality principle, vector space},
we will study \Ds derivative, keeping in mind that dual statement
is true for \sD derivative.

Because product is not commutative, concept of fraction is not defined in division ring.
We need explicitly show from which side denominator has an effect on numerator.
Symbol $*$ in denominator of fraction implements this function.
According to definition
\[
\frac{\partial f}{\partial_* x}=(\partial_* x)^{-1}\partial f
\]
\[
\frac{\partial f}{{}_*\partial x}=\partial f({}_*\partial x)^{-1}
\]
Thus, we may represent
the G\^ateaux \Ds derivative of function $f$
as ratio of change of function to change of argument.
We do not extend this remark to
components of the G\^ateaux differential.
In this case we consider notation
$\displaystyle\frac{\partial}{\partial_* x}$
and $\displaystyle\frac{\partial}{{}_*\partial x}$ not
as fraction, but as symbol of operator.

Then we can write equation \eqref{eq: Gateaux Dstar derivative, division ring}
as
\begin{equation}
\label{eq: Gateaux differential, division ring, 1}
\partial f(x)(dx)
=dx(\partial_* x)^{-1}\partial f(x)(dx)
\end{equation}
It is easy to see that
we wrote \Ds differential of variable $x$
in denominator of fraction in equation
\eqref{eq: Gateaux differential, division ring, 1},
however we multiply fraction by
differential of variable $x$.
\fi

\begin{theorem}
\label{theorem: Dstar derivative is multiplicative over field R, division ring}
Let $D$ be complete division ring of characteristic $0$.
The G\^ateaux \Ds derivative is projective
over field of real numbers $R$.
\end{theorem}
\begin{proof}
Corollary of theorems \ref{theorem: Gateaux derivative is multiplicative over field R, division ring}
and example
\ifx\IntroCalculusBook\Defined
\ref{example: map projective over commutative ring, ring}.
\else
\ifx\EprintCalculus\Defined
\ref{example: map projective over commutative ring, ring}.
\else
\xRef{0701.238}{example: map projective over commutative ring, ring}.
\fi
\fi
\end{proof}

From theorem \ref{theorem: Dstar derivative is multiplicative over field R, division ring}
it follows
\begin{equation}
\label{eq: Dstar derivative is multiplicative over field R, division ring}
\frac{\partial f(x)(ra)}{\partial_* x}
=\frac{\partial f(x)(a)}{\partial_* x}
\end{equation}
for every $r\in R$, $r\ne 0$ and  $a\in D$, $a\ne 0$.
Therefore the G\^ateaux \Ds derivative
is well defined in direction $a$ over field $R$, $a\in D$, $a\ne 0$,
and does not depend on the choice of value in this direction.

\begin{theorem}
Let $D$ be complete division ring of characteristic $0$
and $a\ne 0$.
The G\^ateaux \Ds derivative and the G\^ateaux \sD derivative of map $f$ of division ring $D$
are bounded by relationship
\begin{equation}
\label{eq: Dstar derivative and starD derivative in direction, division ring}
\frac{\partial f(x)(a)}{{}_*\partial x} 
=a\frac{\partial f(x)(a)}{\partial_* x} a^{-1}
\end{equation}
\end{theorem}
\begin{proof}
From equations \eqref{eq: Gateaux Dstar derivative, division ring}
and \eqref{eq: Gateaux starD derivative, division ring}
it follows
\[
\frac{\partial f(x)(a)}{{}_*\partial x}
=\partial f(x)(a)a^{-1}
=a\frac{\partial f(x)(a)}{\partial_* x}a^{-1}
\]
\end{proof}

\begin{theorem}
Let $D$ be complete division ring of characteristic $0$.
The G\^ateaux differential satisfies to relationship
\begin{equation}
\label{eq: Gateaux differential of product, division ring}
\partial (f(x)g(x))(a)
=\partial f(x)(a)\ g(x)+f(x)\ \partial g(x)(a)
\end{equation}
\end{theorem}
\begin{proof}
Equation \eqref{eq: Gateaux differential of product, division ring}
follows from chain of equations
\begin{align*}
\partial (f(x)g(x))(a)
&=\lim_{t\rightarrow 0}(t^{-1}(f(x+ta)g(x+ta)-f(x)g(x)))
\\
&=\lim_{t\rightarrow 0}(t^{-1}(f(x+ta)g(x+ta)-f(x)g(x+ta)))
\\
&+\lim_{t\rightarrow 0}(t^{-1}(f(x)g(x+ta)-f(x)g(x)))
\\
&=\lim_{t\rightarrow 0}(t^{-1}(f(x+ta)-f(x)))g(x)
\\
&+f(x)\lim_{t\rightarrow 0}(t^{-1}(g(x+ta)-g(x)))
\end{align*}
based on definition \eqref{eq: Gateaux differential, representation, 2}.
\end{proof}

\begin{theorem}
Let $D$ be complete division ring of characteristic $0$.
Suppose the G\^ateaux derivative
of map $f:D\rightarrow D$
has expansion
\begin{equation}
\label{eq: Gateaux differential of f, division ring}
\partial f(x)(a)
=\frac{\pC{s}{0}\partial f(x)}{\partial x}
a
\frac{\pC{s}{1}\partial f(x)}{\partial x}
\end{equation}
Suppose the G\^ateaux differential
of map $g:D\rightarrow D$
has expansion
\begin{equation}
\label{eq: Gateaux differential of g, division ring}
\partial g(x)(a)
=\frac{\pC{t}{0}\partial g(x)}{\partial x}
a
\frac{\pC{t}{1}\partial g(x)}{\partial x}
\end{equation}
Components of the G\^ateaux differential of map $f(x)g(x)$
have form
\begin{align}
\label{eq: 0 component of Gateaux derivative, fg, division ring}
\frac{\pC{s}{0}\partial f(x)g(x)}{\partial x}
&
=\frac{\pC{s}{0}\partial f(x)}{\partial x}
&
\frac{\pC{t}{0}\partial f(x)g(x)}{\partial x}
&
=f(x)\frac{\pC{t}{0}\partial g(x)}{\partial x}
\\
\label{eq: 1 component of Gateaux derivative, fg, division ring}
\frac{\pC{s}{1}\partial f(x)g(x)}{\partial x}
&
=\frac{\pC{s}{1}\partial f(x)}{\partial x}
g(x)
&
\frac{\pC{t}{1}\partial f(x)g(x)}{\partial x}
&
=\frac{\pC{t}{1}\partial g(x)}{\partial x}
\end{align}
\end{theorem}
\begin{proof}
Let us substitute \eqref{eq: Gateaux differential of f, division ring}
and \eqref{eq: Gateaux differential of g, division ring}
into equation \eqref{eq: Gateaux differential of product, division ring}
\begin{align}
\label{eq: Gateaux differential of fg, division ring}
\partial (f(x) g(x))(a)
&=\partial f(x)(a)\ g(x)+f(x)\ \partial g(x)(a)
\\ \nonumber
&=\frac{\pC{s}{0}\partial f(x)}{\partial x}a
\frac{\pC{s}{1}\partial f(x)}{\partial x}
g(x)
+f(x)\frac{\pC{t}{0}\partial g(x)}{\partial x}a
\frac{\pC{t}{1}\partial g(x)}{\partial x}
\end{align}
Based \eqref{eq: Gateaux differential of fg, division ring},
we define equations
\eqref{eq: 0 component of Gateaux derivative, fg, division ring},
\eqref{eq: 1 component of Gateaux derivative, fg, division ring}.
\end{proof}

\begin{theorem}
Let $D$ be complete division ring of characteristic $0$.
The G\^ateaux \Ds derivative satisfy to relationship
\begin{equation}
\label{eq: Gateaux Dstar derivative of product, division ring}
\frac{\partial f(x)g(x)}{\partial_* x}(a)
=\frac{\partial f(x)(a)}{\partial_* x}g(x)+a^{-1}f(x)a\frac{\partial g(x)(a)}{\partial_* x}
\end{equation}
\end{theorem}
\begin{proof}
Equation \eqref{eq: Gateaux Dstar derivative of product, division ring}
follows from chain of equations
\begin{align*}
\frac{\partial f(x)g(x)}{\partial_* x}(a)
&=a^{-1}\partial f(x)g(x)(a)
\\
&=a^{-1}(\partial f(x)(a)g(x)+f(x)\partial g(x)(a))
\\
&=a^{-1}\partial f(x)(a)g(x)+a^{-1}f(x)aa^{-1}\partial g(x)(a)
\\
&=\frac{\partial f(x)(a)}{\partial_* x}g(x)+a^{-1}f(x)a\frac{\partial g(x)(a)}{\partial_* x}
\end{align*}
\end{proof}

\begin{theorem}
Let $D$ be complete division ring of characteristic $0$.
Either the G\^ateaux \Ds derivative  does not depend on direction,
or the G\^ateaux \Ds derivative in direction $0$ is not defined.
\end{theorem}
\begin{proof}
Statement of theorem is corollary
of theorem \ref{theorem: Dstar derivative is multiplicative over field R, division ring}
and theorem \ref{theorem: projective map in 0, division ring}.
\end{proof}

\begin{theorem}
\label{theorem: norm of Gateaux differential, division ring}
Let $D$ be complete division ring of characteristic $0$.
Let unit sphere of division ring $D$ be compact.
If the G\^ateaux \Ds derivative
$\displaystyle\frac{\partial f(x)(a)}{\partial_* x}$
exists in point $x$ and is continuous in direction over field $R$,
then there exist norm $\|\partial f(x)\|$ of the G\^ateaux \Ds differential.
\end{theorem}
\begin{proof}
From definition \ref{definition: Gateaux Dstar derivative, division ring}
it follows
\begin{equation}
\label{eq: norm of Gateaux differential, division ring, 1}
|\partial f(x)(a)|
=|a|\ \left|\frac{\partial f(x)(a)}{\partial_* x}\right|
\end{equation}
From theorems
\xRef{0812.4763}{theorem: projective function is continuous in direction over field R, division ring},
\ref{theorem: Dstar derivative is multiplicative over field R, division ring}
it follows, that the G\^ateaux \Ds derivative is continuous
on unit sphere.
Since unit sphere is compact, then range
the G\^ateaux \Ds derivative of function $f$ at point $x$ is bounded
\[
\left|
\frac{\partial f(x)(a)}{\partial_* x}
\right|
<F=\text{sup}\left|
\frac{\partial f(x)(a)}{\partial_* x}
\right|
\]
According to definition \ref{definition: norm of map, division ring}
\[
\|\partial f(x)\|=F
\]
\end{proof}

\begin{theorem}
\label{theorem: Dstar derivative and continuos function, division ring}
Let $D$ be complete division ring of characteristic $0$.
Let unit sphere of division ring $D$ be compact.
If the G\^ateaux \Ds derivative
$\displaystyle\frac{\partial f(x)(a)}{\partial_* x}$
exists in point $x$ and is continuous in direction over field $R$,
then function $f$ is continuous at point $x$.
\end{theorem}
\begin{proof}
From theorem \ref{theorem: norm of Gateaux differential, division ring}
it follows
\begin{equation}
\label{eq: Dstar derivative and continuos function, division ring, 1}
|\partial f(x)(a)|
\le\|\partial f(x)\| |a|
\end{equation}
From \eqref{eq: Gateaux differential of map, division ring},
\eqref{eq: Dstar derivative and continuos function, division ring, 1}
it follows
\begin{equation}
\label{eq: Dstar derivative and continuos function, division ring, 2}
|f(x+a)-f(x)|<|a|\ \|\partial f(x)\|
\end{equation}
Let us assume arbitrary $\epsilon>0$. Assume
\[\delta=\frac\epsilon{\|\partial f(x)\|}\]
Then from inequation
\[
|a|<\delta
\]
it follows
\[
|f(x+a)-f(x)|\le \|\partial f(x)\|\ \delta=\epsilon
\]
According to definition \ref{definition: continuous function, division ring}
map $f$ is continuous at point $x$.
\end{proof}

\ifx\EprintCalculus\undefined
Theorem \ref{theorem: Dstar derivative and continuos function, division ring}
has interesting generalization. If unit sphere of division ring $D$ is not compact,
then we may consider compact set of directions instead of unit sphere.
In this case we may speak about continuity of function $f$ along selected set
of directions.
\fi

\section{Table of Derivatives of Map of Division Ring}

\begin{theorem}
\label{theorem: Gateaux differential of const, division ring}
Let $D$ be complete division ring of characteristic $0$.
Then for any $b\in D$
\begin{equation}
\label{eq: Gateaux differential of const, division ring}
\partial (b)(a)
=0
\end{equation}
\end{theorem}
\begin{proof}
Immediate corollary of definition
\ref{definition: function differentiable in Gateaux sense, division ring}.
\end{proof}

\begin{theorem}
\label{theorem: Gateaux differential, product over constant, division ring}
Let $D$ be complete division ring of characteristic $0$.
Then for any $b$, $c\in D$
\begin{align}
\label{eq: Gateaux differential, product over constant, division ring}
\partial (bf(x)c)(a)
&=b\partial f(x)(a)c
\\
\label{eq: 0 component of Gateaux differential, product over constant, division ring}
\frac{\pC{s}{0}\partial bf(x)c}{\partial x}
&=b\frac{\pC{s}{0}\partial f(x)}{\partial x}
\\
\label{eq: 1 component of Gateaux differential, product over constant, division ring}
\frac{\pC{s}{1}\partial bf(x)c}{\partial x}
&=\frac{\pC{s}{1}\partial f(x)}{\partial x}c
\\
\label{eq: Gateaux derivative, product over constant, division ring}
\frac{\partial bf(x)c}{\partial_* x}(a)
&=a^{-1}ba\frac{\partial f(x)(a)}{\partial_* x}c
\end{align}
\end{theorem}
\begin{proof}
Immediate corollary of equations
\eqref{eq: Gateaux differential of product, division ring},
\eqref{eq: 0 component of Gateaux derivative, fg, division ring},
\eqref{eq: 1 component of Gateaux derivative, fg, division ring},
\eqref{eq: Gateaux Dstar derivative of product, division ring}
because $\partial b=\partial c=0$.
\end{proof}

\begin{theorem}
\label{theorem: Gateaux differential, product over constant, fx=x, division ring}
Let $D$ be complete division ring of characteristic $0$.
Then for any $b$, $c\in D$
\begin{align}
\label{eq: Gateaux differential, product over constant, fx=x, division ring}
\partial (bxc)(h)
&=bhc
\\
\label{eq: 0 component of Gateaux differential, product over constant, fx=x, division ring}
\frac{\pC{1}{0}\partial bxc}{\partial x}
&=b
\\
\label{eq: 1 component of Gateaux differential, product over constant, fx=x, division ring}
\frac{\pC{1}{1}\partial bxc}{\partial x}
&=c
\\
\label{eq: Gateaux derivative, product over constant, fx=x, division ring}
\frac{\partial bxc}{\partial_* x}(h)
&=h^{-1}bhc
\end{align}
\end{theorem}
\begin{proof}
Corollary of theorem
\ref{theorem: Gateaux differential, product over constant, division ring},
when $f(x)=x$.
\end{proof}

\begin{theorem}
Let $D$ be complete division ring of characteristic $0$.
Then for any $b\in D$
\begin{align}
\label{eq: Gateaux differential, xb-bx, division ring}
\partial (xb-bx)(h)
&=hb-bh
\\
\nonumber
\frac{\pC{1}{0}\partial (xb-bx)}{\partial x}
&=1
&
\frac{\pC{1}{1}\partial (xb-bx)}{\partial x}
&=b
\\
\nonumber
\frac{\pC{2}{0}\partial (xb-bx)}{\partial x}
&=-b
&
\frac{\pC{2}{1}\partial (xb-bx)}{\partial x}
&=1
\\
\nonumber
\frac{\partial (xb-bx)}{\partial_* x}(h)
&=h^{-1}bhc
\end{align}
\end{theorem}
\begin{proof}
Corollary of theorem
\ref{theorem: Gateaux differential, product over constant, division ring},
when $f(x)=x$.
\end{proof}

\begin{theorem}
Let $D$ be complete division ring of characteristic $0$.
Then
\begin{equation}
 \begin{array}{rl}
\partial (x^2)(a)&=xa+ax
\\
\displaystyle\frac{\partial x^2}{\partial_*x}(a)
&=a^{-1}xa+x
\end{array}
\label{eq: derivative x2, division ring}
\end{equation}
\begin{equation}
 \begin{array}{ccc}
\displaystyle\frac{\pC{1}{0}\partial x^2}{\partial x}=x
&
\displaystyle\frac{\pC{1}{1}\partial x^2}{\partial x}=e
\\ \\
\displaystyle\frac{\pC{2}{0}\partial x^2}{\partial x}=e
&
\displaystyle\frac{\pC{2}{1}\partial x^2}{\partial x}=x
\end{array}
\label{eq: derivative x2, component of Gateaux differential, division ring}
\end{equation}
\end{theorem}
\begin{proof}
\eqref{eq: derivative x2, division ring} follows from example \ref{example: differential dxx}
and definition \ref{definition: Gateaux Dstar derivative, division ring}.
\eqref{eq: derivative x2, component of Gateaux differential, division ring}
follows from example \ref{example: differential dxx}
and equation \eqref{eq: Gateaux differential, form}.
\end{proof}

\begin{theorem}
Let $D$ be complete division ring of characteristic $0$.
Then
\begin{align}
\partial (x^{-1})(h)&=-x^{-1}hx^{-1}
\label{eq: derivative x power -1, division ring}
\\
\displaystyle\frac{\partial x^{-1}}{\partial_*x}(h)
&=-h^{-1}x^{-1}hx^{-1}
\nonumber
\end{align}
\[
 \begin{array}{ccc}
\displaystyle\frac{\pC{1}{0}\partial x^{-1}}{\partial x}=-x^{-1}
&
\displaystyle\frac{\pC{1}{1}\partial x^{-1}}{\partial x}=x^{-1}
\end{array}
\]
\end{theorem}
\begin{proof}
Let us substitute $f(x)=x^{-1}$ in definition \eqref{eq: Gateaux differential, representation, 2}.
\begin{align}
\partial f(x)(h)
&=\lim_{t\rightarrow 0,\ t\in R}(t^{-1}((x+th)^{-1}-x^{-1}))
\nonumber
\\
&=\lim_{t\rightarrow 0,\ t\in R}(t^{-1}((x+th)^{-1}-x^{-1}(x+th)(x+th)^{-1}))
\nonumber
\\
&=\lim_{t\rightarrow 0,\ t\in R}(t^{-1}(1-x^{-1}(x+th))(x+th)^{-1})
\label{eq: derivative x power -1, division ring, 1}
\\
&=\lim_{t\rightarrow 0,\ t\in R}(t^{-1}(1-1-x^{-1}th)(x+th)^{-1})
\nonumber
\\
&=\lim_{t\rightarrow 0,\ t\in R}(-x^{-1}h(x+th)^{-1})
\nonumber
\end{align}
Equation \eqref{eq: derivative x power -1, division ring}
follows from chain of equations
\eqref{eq: derivative x power -1, division ring, 1}.
\end{proof}

\begin{theorem}
Let $D$ be complete division ring of characteristic $0$.
Then
\begin{align}
\partial (xax^{-1})(h)&=hax^{-1}-xax^{-1}hx^{-1}
\label{eq: derivative xax power -1, division ring}
\\
\displaystyle\frac{\partial xax^{-1}}{\partial_*x}(h)
&=ax^{-1}-h^{-1}xax^{-1}hx^{-1}
\nonumber
\end{align}
\[
 \begin{array}{lll}
\displaystyle\frac{\pC{1}{0}\partial x^{-1}}{\partial x}=1
&
\displaystyle\frac{\pC{1}{1}\partial x^{-1}}{\partial x}=ax^{-1}
\\
\displaystyle\frac{\pC{2}{0}\partial x^{-1}}{\partial x}=-xax^{-1}
&
\displaystyle\frac{\pC{2}{1}\partial x^{-1}}{\partial x}=x^{-1}
\end{array}
\]
\end{theorem}
\begin{proof}
Equation \eqref{eq: derivative xax power -1, division ring}
is corollary of equations \eqref{eq: Gateaux differential of product, division ring},
\eqref{eq: Gateaux differential, product over constant, fx=x, division ring},
\eqref{eq: derivative xax power -1, division ring}.
\end{proof}

\section{Derivative of Second Order of Map of Division Ring}

Let $D$ be valued division ring of characteristic $0$.
Let
\[f:D\rightarrow D\]
function differentiable in the G\^ateaux sense.
According to remark \ref{remark: differential L(D,D)}
the G\^ateaux derivative is map
\ShowEq{differential L(D,D)}
According to theorems
\ifx\IntroCalculusBook\Defined
\ref{theorem: sum of additive maps, ring},
\ref{theorem: additive map times constant, ring}
\else
\ifx\EprintCalculus\Defined
\ref{theorem: sum of additive maps, ring},
\ref{theorem: additive map times constant, ring}
\else
\xRef{0701.238}{theorem: sum of additive maps, ring},
\xRef{0701.238}{theorem: additive map times constant, ring}
\fi
\fi
and definition \ref{definition: norm of map, division ring}
set $\mathcal L(D;D)$ is
normed $D$\Hyph vector space.
Therefore, we may consider the question,
if map $\partial f$ is differentiable in the G\^ateaux sense.

According to definition
\xRef{0812.4763}{definition: function differentiable in Gateaux sense, D vector space}
\ShowEq{Gateaux differential of map df, division ring}
where
\ShowEq{Gateaux differential of map df, division ring, o2}
is such continuous map, that
\ShowEq{Gateaux differential of map df, division ring, o2 lim}
According to definition
\xRef{0812.4763}{definition: function differentiable in Gateaux sense, D vector space}
map $\partial(\partial f(x)(a_1))(a_2)$
is linear map of variable $a_2$. From equation
\eqref{eq: Gateaux differential of map df, division ring}
it follows that map
$\partial(\partial f(x)(a_1))(a_2)$
is linear map of variable $a_1$.

\begin{definition}
\label{definition: Gateaux Differential of Second Order, division ring}
Polylinear map
\symb{\partial^2 f(x)}
0{Gateaux derivative of Second Order, division ring}
\symb{\frac{\partial^2 f(x)}{\partial x^2}}
0{Gateaux derivative of Second Order, fraction, division ring}
\begin{equation}
\label{eq: Gateaux derivative of Second Order, division ring}
\ShowSymbol{Gateaux derivative of Second Order, division ring}(a_1;a_2)
=\ShowSymbol{Gateaux derivative of Second Order, fraction, division ring}(a_1;a_2)
=\partial(\partial f(x)(a_1))(a_2)
\end{equation}
is called
\AddIndex{the G\^ateaux derivative of second order of map $f$}
{Gateaux derivative of Second Order, division ring}
\qed
\end{definition}

\begin{remark}
\label{remark: differential L(D,D;D)}
According to definition
\ref{definition: Gateaux Differential of Second Order, division ring}
for given $x$ the G\^ateaux differential of second order
\ShowEq{differential L(D,D;D), 1}.
Therefore, the G\^ateaux differential of second order of map $f$ is map
\ShowEq{differential L(D,D;D)}
\end{remark}

\begin{theorem}
\label{theorem: Gateaux Differential of Second Order, division ring}
\symb{\partial^2 f(x)(a_1;a_2)}
0{Gateaux differential of Second Order, division ring}
It is possible to represent \AddIndex{the G\^ateaux differential of second order of map $f$}
{Gateaux differential of Second Order, division ring}
as
\ShowEq{Differential of Second Order, division ring, representation}
\end{theorem}
\begin{proof}
Corollary of definition \ref{definition: Gateaux Differential of Second Order, division ring}
and theorem
\ifx\IntroCalculusBook\Defined
\ref{theorem: polyadditive map, D vector space}.
\else
\xRef{0701.238}{theorem: polyadditive map, D vector space}.
\fi
\end{proof}

\begin{definition}
\begin{sloppypar}
Expression\footnote{We suppose
\[
\frac{\pC{s}{p}\partial^2 f(x)}{\partial x^2}=\frac{\pC{s}{p}\partial^2 f(x)}{\partial x\partial x}
\]}
\symb{\frac{\pC{s}{p}\partial^2 f(x)}{\partial x^2}}0
{component of Gateaux derivative of Second Order, division ring}
$\displaystyle\ShowSymbol{component of Gateaux derivative of Second Order, division ring}$,
$p=0$, $1$,
is called
\AddIndex{component of the G\^ateaux derivative of map $f(x)$}
{component of Gateaux derivative of Second Order, division ring}.
\qed
\end{sloppypar}
\end{definition}

By induction, assuming that the G\^ateaux derivative is defined
$\partial^{n-1} f(x)$ up to order $n-1$, we define
\symb{\partial^n f(x)}
0{Gateaux derivative of Order n, division ring}
\symb{\frac{\partial^n f(x)}{\partial x^n}}
0{Gateaux derivative of Order n, fraction, division ring}
\begin{equation}
\label{eq: Gateaux derivative of Order n, division ring}
\ShowSymbol{Gateaux derivative of Order n, division ring}(a_1;...;a_n)
=\ShowSymbol{Gateaux derivative of Order n, fraction, division ring}(a_1;...;a_n)
=\partial(\partial^{n-1} f(x)(a_1;...;a_{n-1}))(a_n)
\end{equation}
\AddIndex{the G\^ateaux derivative of order $n$ of map $f$}
{Gateaux derivative of Order n, division ring}.
We also assume $\partial^0 f(x)=f(x)$.

\section{Taylor Series}
\label{section: Taylor Series}

Let us consider polynomial in one variable over division ring $D$ of power $n$, $n>0$. We
want to explore the structure of monomial $p_k(x)$ of polynomial of power $k$.

It is evident that monomial of power $0$ has form $a_0$, $a_0\in D$.
Let $k>0$. Let us prove that
\[
p_k(x)=p_{k-1}(x)xa_k
\]
where $a_k\in D$.
Actually, last factor of monomial $p_k(x)$ is either $a_k\in D$,
or has form $x^l$, $l\ge 1$.
In the later case we assume $a_k=1$.
Factor preceding $a_k$ has form $x^l$, $l\ge 1$.
We can represent this factor as $x^{l-1}x$.
Therefore, we proved the statement.

In particular, monomial of power $1$ has form $p_1(x)=a_0xa_1$.

Without loss of generality, we assume $k=n$.

\begin{theorem}
\label{theorem: Gateaux derivative of f(x)x, division ring}
For any $m>0$ the following equation is true
\begin{align}
&\partial^m (f(x)x)(h_1;...;h_m)
\nonumber \\
=&\partial^m f(x)(h_1;...;h_m)x
+\partial^{m-1} f(x)(h_1;...;h_{m-1})h_m
\label{eq: Gateaux derivative of f(x)x, division ring}
\\ \nonumber
+&\partial^{m-1} f(x)(\widehat{h_1};...;h_{m-1};h_m)h_1
+...
\\ \nonumber
+&\partial^{m-1} f(x)(h_1;...;\widehat{h_{m-1}};h_m)h_{m-1}
\end{align}
where symbol $\widehat{h^i}$ means absense of variable $h^i$ in the list.
\end{theorem}
\begin{proof}
For $m=1$, this is corollary  of equation \eqref{eq: Gateaux differential of product, division ring}
\[
\partial (f(x)x)(h_1)=\partial f(x)(h_1)x+
f(x)h_1
\]
Assume, \eqref{eq: Gateaux derivative of f(x)x, division ring} is true
for $m-1$. Than
\begin{align}
&\partial^{m-1} (f(x)x)(h_1;...;h_{m-1})
\nonumber \\
=&\partial^{m-1} f(x)(h_1;...;h_{m-1})x
+\partial^{m-2} f(x)(h_1;...;h_{m-2})h_{m-1}
\label{eq: Gateaux derivative of f(x)x, m-1, division ring}
\\ \nonumber
+&\partial^{m-2} f(x)(\widehat{h_1},...,h_{m-2},h_{m-1})h_1
+...
\\ \nonumber
+&\partial^{m-2} f(x)(h_1,...,\widehat{h_{m-2}},h_{m-1})h_{m-2}
\end{align}
Using equations \eqref{eq: Gateaux differential of product, division ring}
and \eqref{eq: Gateaux differential, product over constant, division ring}
we get
\begin{align}
&\partial^m (f(x)x)(h_1;...;h_{m-1};h_m)
\nonumber \\
=&\partial^m f(x)(h_1;...;h_{m-1};h_m)x
+\partial^{m-1} f(x)(h_1;...;h_{m-2};h_{m-1})h_m
\nonumber \\
+&\partial^{m-1} f(x)(h_1;...;h_{m-2};\widehat{h_{m-1}};h_m)h_{m-1}
\label{eq: Gateaux derivative of f(x)x, m, division ring}
\\ \nonumber
+&\partial^{m-2} f(x)(\widehat{h_1},...,h_{m-2},h_{m-1};h_m)h_1
+...
\\ \nonumber
+&\partial^{m-2} f(x)(h_1,...,\widehat{h_{m-2}},h_{m-1};h_m)h_{m-2}
\end{align}
The difference between
equations \eqref{eq: Gateaux derivative of f(x)x, division ring} and
\eqref{eq: Gateaux derivative of f(x)x, m, division ring}
is only in form of presentation.
We proved the theorem.
\end{proof}

\begin{theorem}
\label{theorem: Gateaux derivative of pn is symmetric, m less than n, division ring}
The G\^ateaux derivative $\partial^m p_n(x)(h_1,...,h_m)$
is symmetric polynomial with respect to variables $h_1$, ..., $h_m$.
\end{theorem}
\begin{proof}
To prove the theorem we consider algebraic properties
of the G\^ateaux derivative and give equivalent definition.
We start from construction of monomial.
For any monomial $p_n(x)$ we build symmetric
polynomial $r_n(x)$ according to following rules
\begin{itemize}
\item If $p_1(x)=a_0xa_1$, then $r_1(x_1)=a_0x_1a_1$
\item If $p_n(x)=p_{n-1}(x)a_n$, then
\[
r_n(x_1,...,x_n)=r_{n-1}(x_{[1},...,x_{n-1})x_{n]}a_n
\]
where square brackets express symmetrization of expression
with respect to variables $x_1$, ..., $x_n$.
\end{itemize}
It is evident that
\[
p_n(x)=r_n(x_1,...,x_n)\ \ \ x_1=...=x_n=x
\]
We define the G\^ateaux derivative of power $k$ according to rule
\begin{equation}
\partial^k p_n(x)(h_1,...,h_k)=r_n(h_1,...,h_k,x_{k+1},...,x_n)\ \ \ x_{k+1}=x_n=x
\label{eq: Gateaux derivative of monomial, algebraic definition, division ring}
\end{equation}
According to construction, polynomial $r_n(h_1,...,h_k,x_{k+1},...,x_n)$
is symmetric with respect to variables $h_1$, ..., $h_k$, $x_{k+1}$, ..., $x_n$.
Therefore, polynomial
\eqref{eq: Gateaux derivative of monomial, algebraic definition, division ring}
is symmetric with respect to variables $h_1$, ..., $h_k$.

For $k=1$, we will prove that definition
\eqref{eq: Gateaux derivative of monomial, algebraic definition, division ring}
of the G\^ateaux derivative coincides with definition
\eqref{eq: Gateaux Dstar derivative, division ring}.

For $n=1$, $r_1(h_1)=a_0h_1a_1$. This expression coincides
with expression of the G\^ateaux derivative in theorem
\ref{theorem: Gateaux differential, product over constant, fx=x, division ring}.

Let the statement be true for $n-1$.
The following equation is true
\begin{equation}
r_n(h_1,x_2,...,x_n)=r_{n-1}(h_1,x_{[2},...,x_{n-1})x_{n]}a_n+r_{n-1}(x_2,...,x_n)h_1a_n
\label{eq: Gateaux derivative of monomial, n, division ring}
\end{equation}
Assume $x_2=...=x_n=x$. 
According to suggestion of induction, from equations
\eqref{eq: Gateaux derivative of monomial, algebraic definition, division ring},
\eqref{eq: Gateaux derivative of monomial, n, division ring}
it follows that
\[
r_n(h_1,x_2,...,x_n)=\partial p_{n-1}(x)(h_1)xa_n+p_{n-1}(x)h_1a_n
\]

According to theorem
\ref{theorem: Gateaux derivative of f(x)x, division ring}
\[
r_n(h_1,x_2,...,x_n)=\partial p_n(x)(h_1)
\]
This proves the equation
\eqref{eq: Gateaux derivative of monomial, algebraic definition, division ring}
for $k=1$.

Let us prove now that definition
\eqref{eq: Gateaux derivative of monomial, algebraic definition, division ring}
of the G\^ateaux derivative coincides with definition
\eqref{eq: Gateaux derivative of Order n, division ring} for $k>1$.

Let equation
\eqref{eq: Gateaux derivative of monomial, algebraic definition, division ring}
be true for $k-1$.
Let us consider arbitrary monomial of polynomial $r_n(h_1,...,h_{k-1},x_k,...,x_n)$.
Identifying variables $h_1$, ..., $h_{k-1}$
with elements of division ring $D$, we consider polynomial
\begin{equation}
R_{n-k}(x_k,...,x_n)=r_n(h_1,...,h_{k-1},x_k,...,x_n)
\label{eq: reduced polynom, division ring}
\end{equation}
Assume $P_{n-k}(x)=R_{n-k}(x_k,...,x_n)$, $x_k=...=x_n=x$.
Therefore
\[
P_{n-k}(x)=\partial^{k-1} p_n(x)(h_1;...;h_{k-1})
\]
According to definition of the G\^ateaux derivative
\eqref{eq: Gateaux derivative of Order n, division ring}
\begin{align}
\partial P_{n-k}(x)(h_k)
=&\partial(\partial^{k-1} p_n(x)(h_1;...;h_{k-1}))(h_k)
\nonumber\\
=&\partial^k p_n(x)(h_1;...;h_{k-1};h_k)
\label{eq: Gateaux derivative of Order n, 1, division ring}
\end{align}
According to definition
\eqref{eq: Gateaux derivative of monomial, algebraic definition, division ring}
of the G\^ateaux derivative
\begin{equation}
\partial P_{n-k}(x)(h_k)=R_{n-k}(h_k,x_{k+1},...,x_n)\ \ \ x_{k+1}=x_n=x
\label{eq: Gateaux derivative of monomial, algebraic definition, 1, division ring}
\end{equation}
According to definition \eqref{eq: reduced polynom, division ring},
from equation
\eqref{eq: Gateaux derivative of monomial, algebraic definition, 1, division ring}
it follows that
\begin{equation}
\partial P_{n-k}(x)(h_k)=r_n(h_1,...,h_k,x_{k+1},...,x_n)\ \ \ x_{k+1}=x_n=x
\label{eq: Gateaux derivative of monomial, algebraic definition, 2, division ring}
\end{equation}
From comparison of equations \eqref{eq: Gateaux derivative of Order n, 1, division ring}
and \eqref{eq: Gateaux derivative of monomial, algebraic definition, 2, division ring}
it follows that
\[
\partial^k p_n(x)(h_1;...;h_k)=r_n(h_1,...,h_k,x_{k+1},...,x_n)\ \ \ x_{k+1}=x_n=x
\]
Therefore equation
\eqref{eq: Gateaux derivative of monomial, algebraic definition, division ring}
is true for any $k$ and $n$.

We proved the statement of theorem.
\end{proof}

\begin{theorem}
\label{theorem: Gateaux derivative of pn equal 0, division ring}
For any $n\ge 0$ following equation is true
\begin{equation}
\label{eq: Gateaux derivative n1 of pn, division ring}
\partial^{n+1} p_n(x)(h_1;...;h_{n+1})=0
\end{equation}
\end{theorem}
\begin{proof}
Since $p_0(x)=a_0$, $a_0\in D$, then for $n=0$ theorem is corollary
of theorem \ref{theorem: Gateaux differential of const, division ring}.
Let statement of theorem is true for $n-1$. According to theorem
\ref{theorem: Gateaux derivative of f(x)x, division ring}
when $f(x)=p_{n-1}(x)$ we get
\begin{align*}
\partial^{n+1} p_n(x)(h_1;...;h_{n+1})
=&\partial^{n+1}(p_{n-1}(x)xa_n)(h_1;...;h_{n+1})
\\
=&\partial^{n+1} p_{n-1}(x)(h_1;...;h_m)xa_n
\\
+&\partial^n p_{n-1}(x)(h_1;...;h_{m-1})h_ma_n
\\
+&\partial^n p_{n-1}(x)(\widehat{h_1};...;h_{m-1};h_m)h_1a_n
+...
\\
+&\partial^{m-1} p_{n-1}(x)(h_1;...;\widehat{h_{m-1}};h_m)h_{m-1}a_n
\end{align*}
According to suggestion of induction all monomials are equal $0$.
\end{proof}

\begin{theorem}
\label{theorem: Gateaux derivative of pn equal 0, m less than n, division ring}
If $m<n$, then following equation is true
\begin{equation}
\partial^m p_n(0)(h)=0
\label{eq: m Gateaux derivative of polinom pn, division ring}
\end{equation}
\end{theorem}
\begin{proof}
For $n=1$ following equation is true
\[
\partial^0 p_1(0)=a_0xa_1=0
\]
Assume that statement is true
for $n-1$. Then according to theorem \ref{theorem: Gateaux derivative of f(x)x, division ring}
\begin{align*}
&\partial^m (p_{n-1}(x)xa_n)(h_1;...;h_m)
\\
=&\partial^m p_{n-1}(x)(h_1;...;h_m)xa_n
+\partial^{m-1} p_{n-1}(x)(h_1;...;h_{m-1})h_ma_n
\\
+&\partial^{m-1} p_{n-1}(x)(\widehat{h_1};...;h_{m-1};h_m)h_1a_n
+...
\\
+&\partial^{m-1} p_{n-1}(x)(h_1;...;\widehat{h_{m-1}};h_m)h_{m-1}a_n
\end{align*}
First term equal $0$ because $x=0$.
Because $m-1<n-1$, then
rest terms equal $0$ according to suggestion of induction.
We proved the statement of theorem.
\end{proof}

When $h_1=...=h_n=h$, we assume
\[
\partial^nf(x)(h)=\partial^nf(x)(h_1;...;h_n)
\]
This notation does not create ambiguity, because we can determine function
according to number of arguments.

\begin{theorem}
\label{theorem: n Gateaux derivative of polinom pn, division ring}
For any $n>0$ following equation is true
\begin{equation}
\partial^n p_n(x)(h)=n!p_n(h)
\label{eq: n Gateaux derivative of polinom pn, division ring}
\end{equation}
\end{theorem}
\begin{proof}
For $n=1$ following equation is true
\[
\partial p_1(x)(h)=\partial (a_0xa_1)(h)=a_0ha_1=1!p_1(h)
\]
Assume the statement is true
for $n-1$. Then according to theorem \ref{theorem: Gateaux derivative of f(x)x, division ring}
\begin{align}
\label{eq: n Gateaux derivative of polinom pn, 1, division ring}
\partial^n p_n(x)(h)
=&\partial^n p_{n-1}(x)(h)xa_n
+\partial^{n-1} p_{n-1}(x)(h)ha_n
\\ \nonumber
+&...
+\partial^{n-1} p_{n-1}(x)(h)ha_n
\end{align}
First term equal $0$ according to theorem
\ref{theorem: Gateaux derivative of pn equal 0, division ring}.
The rest $n$ terms equal, and according to suggestion of induction
from equation \eqref{eq: n Gateaux derivative of polinom pn, 1, division ring}
it follows
\[
\partial^n p_n(x)(h)
=n\partial^{n-1} p_{n-1}(x)(h)ha_n=n(n-1)!p_{n-1}(h)ha_n=n!p_n(h)
\]
Therefore, statement of theorem is true for any $n$.
\end{proof}

Let $p(x)$ be polynomial of power $n$.\footnote{I consider
Taylor polynomial for polynomials by analogy with
construction of Taylor polynomial in \citeBib{Fihtengolts: Calculus volume 1}, p. 246.}
\[
p(x)=p_0+p_{1i_1}(x)+...+p_{ni_n}(x)
\]
We assume sum by index $i_k$ which enumerates terms
of power $k$.
According to theorem \ref{theorem: Gateaux derivative of pn equal 0, division ring},
\ref{theorem: Gateaux derivative of pn equal 0, m less than n, division ring},
\ref{theorem: n Gateaux derivative of polinom pn, division ring}
\[
\partial^k p(x)(h_1;...;h_k)=k!p_{ki_k}(x)
\]
Therefore, we can write
\[
p(x)=p_0+(1!)^{-1}\partial p(0)(x)+(2!)^{-1}\partial^2 p(0)(x)+...+(n!)^{-1}\partial^n p(0)(x)
\]
This representation of polynomial is called
\AddIndex{Taylor polynomial}{Taylor polynomial, division ring}.
If we consider substitution of variable $x=y-y_0$, then considered above construction
remain true for polynomial
\[
p(y)=p_0+p_{1i_1}(y-y_0)+...+p_{ni_n}(y-y_0)
\]
Therefore
\[
p(y)=p_0+(1!)^{-1}\partial p(y_0)(y-y_0)+(2!)^{-1}\partial^2 p(y_0)(y-y_0)+...+(n!)^{-1}\partial^n p(y_0)(y-y_0)
\]

Assume that function $f(x))$ is differentiable in the G\^ateaux sense at point $x_0$
up to any order.\footnote{I explore
construction of Taylor series by analogy with
construction of Taylor series in \citeBib{Fihtengolts: Calculus volume 1}, p. 248, 249.}

\begin{theorem}
\label{theorem: n Gateaux derivative equal 0, division ring}
If function $f(x)$ holds
\begin{equation}
f(x_0)=\partial f(x_0)(h)=...=\partial^n f(x_0)(h)=0
\label{eq: n Gateaux derivatives of function, division ring}
\end{equation}
then for $t\rightarrow 0$ expression $f(x+th)$ is infinitesimal of order
higher then $n$ with respect to $t$
\[
f(x_0+th)=o(t^n)
\]
\end{theorem}
\begin{proof}
When $n=1$ this statement follows from equation
\eqref{eq: Gateaux differential of map, t, division ring}.

Let statement be true for $n-1$.
Map
\[
f_1(x)=\partial f(x)(h)
\]
satisfies to condition
\[
f_1(x_0)=\partial f_1(x_0)(h)=...=\partial^{n-1} f_1(x_0)(h)=0
\]
According to suggestion of induction
\[
f_1(x_0+th)=o(t^{n-1})
\]
Then equation \eqref{eq: Gateaux differential, representation, 2}
gets form
\[
o(t^{n-1})=\lim_{t\rightarrow 0,\ t\in R}(t^{-1}f(x+th))
\]
Therefore,
\[
f(x+th)=o(t^n)
\]
\end{proof}

Let us form polynomial
\begin{equation}
\label{eq: Taylor polynomial, f(x), division ring}
p(x)=f(x_0)+(1!)^{-1}\partial f(x_0)(x-x_0)+...+(n!)^{-1}\partial^n f(x_0)(x-x_0)
\end{equation}
According to theorem \ref{theorem: n Gateaux derivative equal 0, division ring}
\[
f(x_0+t(x-x_0))-p(x_0+t(x-x_0))=o(t^n)
\]
Therefore, polynomial $p(x)$ is good approximation of map $f(x)$.

If map $f(x)$ has the G\^ateaux derivative of any order, the passing to the limit
$n\rightarrow\infty$, we get expansion into series
\[
f(x)=\sum_{n=0}^{\infty}(n!)^{-1}\partial^n f(x_0)(x-x_0)
\]
which is called \AddIndex{Taylor series}{Taylor series, division ring}.

\section{Integral}

Concept of integral has different aspect. In this section we consider
integration as operation inverse to differentiation.
As a matter of fact, we consider procedure of solution of ordinary differential
equation
\[
\partial f(x)(h)=F(x;h)
\]

\begin{example}
I start from example of differential equation over real field.
\begin{equation}
\label{eq: differential equation y=xx, 1, real number}
y'=3x^2
\end{equation}
\begin{equation}
\label{eq: differential equation y=xx, initial, real number}
\begin{matrix}
x_0=0&y_0=0
\end{matrix}
\end{equation}
Differentiating one after another equation \eqref{eq: differential equation y=xx, 1, real number},
we get the chain of equations
\begin{align}
\label{eq: differential equation y=xx, 2, real number}
y''&=6x
\\
\label{eq: differential equation y=xx, 3, real number}
y'''&=6
\\
\label{eq: differential equation y=xx, 4, real number}
y^{(n)}&=0&n&>3
\end{align}
From equations
\eqref{eq: differential equation y=xx, 1, real number},
\eqref{eq: differential equation y=xx, initial, real number},
\eqref{eq: differential equation y=xx, 2, real number},
\eqref{eq: differential equation y=xx, 3, real number},
\eqref{eq: differential equation y=xx, 4, real number}
it follows expansion into Taylor series
\[
y=x^3
\]
\qed
\end{example}

\begin{example}
Let us consider similar equation over division ring
\begin{equation}
\label{eq: differential equation y=xx, 1, division ring}
\partial (y)(h)=hx^2+xhx+x^2h
\end{equation}
\begin{equation}
\label{eq: differential equation y=xx, initial, division ring}
\begin{matrix}
x_0=0&y_0=0
\end{matrix}
\end{equation}
Differentiating one after another equation
\eqref{eq: differential equation y=xx, 1, division ring},
we get the chain of equations
\ShowEq{differential equation y=xx, division ring}
From equations
\eqref{eq: differential equation y=xx, 1, division ring},
\eqref{eq: differential equation y=xx, initial, division ring},
\eqref{eq: differential equation y=xx, 2, division ring},
\eqref{eq: differential equation y=xx, 3, division ring},
\eqref{eq: differential equation y=xx, 4, division ring}
expansion into Taylor series follows
\[
y=x^3
\]
\qed
\end{example}

\begin{remark}
Differential equation
\begin{equation}
\label{eq: differential equation y=xx, 1, a, division ring}
\partial (y)(h)=3hx^2
\end{equation}
\begin{equation}
\label{eq: differential equation y=xx, initial, a, division ring}
\begin{matrix}
x_0=0&y_0=0
\end{matrix}
\end{equation}
also leads to answer $y=x^3$. it is evident that this
map does not satisfies differential equation.
However, contrary to theorem
\ref{theorem: Gateaux derivative of pn is symmetric, m less than n, division ring}
second derivative is
not symmetric polynomial. This means that
equation \eqref{eq: differential equation y=xx, 1, a, division ring}
does not possess a solution.
\qed
\end{remark}

\begin{example}
\label{example: differential equation, additive function, division ring}
It is evident that, if function satisfies to differential equation
\ShowEq{differential equation, additive function, division ring}
then The G\^ateaux derivative of second order
\[
\partial^2 f(x)(h_1;h_2)=0
\]
Than, if initial condition is $y(0)=0$,
then differential equation
\eqref{eq: differential equation, additive function, division ring} has solution
\ShowEq{differential equation, additive function, solution, division ring}
\qed
\end{example}

\section{Exponent}

In this section we consider one of possible models
of exponent.

In a field we can define exponent as solution of differential equation
\ShowEq{exponent over field}
It is evident that we cannot write such equation for division ring.
However we can use equation
\ShowEq{derivative over field}
From equations
\eqref{eq: exponent over field}, \eqref{eq: derivative over field}
it follows
\ShowEq{exponent derivative over field}
This equation is closer to our goal,
however there is the question: in which order
we should multiply $y$ and $h$?
To answer this question we change equation
\ShowEq{exponent derivative over division ring}
Hence, our goal is to solve differential equation
\eqref{eq: exponent derivative over division ring}
with initial condition $y(0)=1$.

For the statement and proof of the theorem
\ref{theorem: exponent derivative n over division ring}
I introduce following notation.
Let
\ShowEq{exponent derivative, transposition n}
be transposition of the tuple of variables
\ShowEq{exponent transposition tuple 1}
Let $p_{\sigma}(h_i)$ be position that
variable $h_i$ gets in the tuple
\ShowEq{exponent transposition tuple 2}
For instance, if transposition $\sigma$ has form
\ShowEq{exponent transposition tuple 3}
then following tuples equal
\ShowEq{exponent transposition tuple 4}

\begin{theorem}
\label{theorem: exponent derivative n over division ring}
If function $y$ is solution of differential equation
\eqref{eq: exponent derivative over division ring}
then the G\^ateaux derivative of order $n$ of function $y$
has form
\ShowEq{exponent derivative n over division ring}
where sum is over transpositions
\ShowEq{exponent derivative, transposition n}
of the set of variables $y$, $h_1$, ..., $h_n$.
Transposition $\sigma$ has following properties
\begin{enumerate}
\item If there exist $i$, $j$, $i\ne j$, such that
$p_{\sigma}(h_i)$ is situated in product
\eqref{eq: exponent derivative n over division ring}
on the left side of $p_{\sigma}(h_j)$ and $p_{\sigma}(h_j)$ is situated
on the left side of $p_{\sigma}(y)$, then $i<j$.
\label{enumerate: exponent derivative, transposition n, left}
\item If there exist $i$, $j$, $i\ne j$, such that
$p_{\sigma}(h_i)$ is situated in product
\eqref{eq: exponent derivative n over division ring}
on the right side of $p_{\sigma}(h_j)$ and $p_{\sigma}(h_j)$ is situated
on the right side of $p_{\sigma}(y)$, then $i>j$.
\label{enumerate: exponent derivative, transposition n, right}
\end{enumerate}
\end{theorem}
\begin{proof}
We prove this statement by induction.
For $n=1$ the statement is true because this is differential equation
\eqref{eq: exponent derivative over division ring}.
Let the statement be true for $n=k-1$.
Hence
\ShowEq{exponent derivative n=k-1 over division ring}
where the sum is over transposition
\ShowEq{exponent derivative, transposition n=k-1}
of the set of variables $y$, $h_1$, ..., $h_{k-1}$.
Transposition $\sigma$ satisfies to conditions
\eqref{enumerate: exponent derivative, transposition n, left},
\eqref{enumerate: exponent derivative, transposition n, right}
in theorem.
According to definition
\eqref{eq: Gateaux derivative of Order n, division ring}
the G\^ateaux derivative of order $k$ has form
\ShowEq{exponent derivative n=k over division ring, 1}
From equations
\eqref{eq: exponent derivative over division ring},
\eqref{eq: exponent derivative n=k over division ring, 1}
it follows that
\ShowEq{exponent derivative n=k over division ring, 2}
It is easy to see that arbitrary transposition $\sigma$
from sum
\eqref{eq: exponent derivative n=k over division ring, 2}
forms two transpositions
\ShowEq{exponent derivative, transposition n=k}
From \eqref{eq: exponent derivative n=k over division ring, 2}
and \eqref{eq: exponent derivative, transposition n=k}
it follows that
\ShowEq{exponent derivative n=k over division ring, 3}
In expression
\eqref{eq: exponent derivative n=k over division ring, 3}
$p_{\tau_1}(h_k)$ is written immediately before $p_{\tau_1}(y)$.
Since $k$ is smallest value of index
then transposition $\tau_1$ satisfies to conditions
\eqref{enumerate: exponent derivative, transposition n, left},
\eqref{enumerate: exponent derivative, transposition n, right}
in the theorem.
In expression
\eqref{eq: exponent derivative n=k over division ring, 3}
$p_{\tau_2}(h_k)$ is written immediately after $p_{\tau_2}(y)$.
Since $k$ is largest value of index
then transposition $\tau_2$ satisfies to conditions
\eqref{enumerate: exponent derivative, transposition n, left},
\eqref{enumerate: exponent derivative, transposition n, right}
in the theorem.

It remains to show that in the expression
\eqref{eq: exponent derivative n=k over division ring, 3}
we get all transpositions $\tau$ that satisfy to conditions
\eqref{enumerate: exponent derivative, transposition n, left},
\eqref{enumerate: exponent derivative, transposition n, right}
in the theorem.
Since $k$ is largest index then according to conditions
\eqref{enumerate: exponent derivative, transposition n, left},
\eqref{enumerate: exponent derivative, transposition n, right}
in the theorem $\tau(h_k)$
is written either immediately before or immediately after $\tau(y)$.
Therefore, any transposition $\tau$ has
either form $\tau_1$ or form $\tau_2$.
Using equation
\eqref{eq: exponent derivative, transposition n=k},
we can find corresponding
transposition $\sigma$ for given transposition $\tau$.
Therefore, the statement of theorem is true for $n=k$.
We proved the theorem.
\end{proof}

\begin{theorem}
The solution of differential equation
\eqref{eq: exponent derivative over division ring}
with initial condition $y(0)=1$
is exponent
\ShowEq{exponent over division ring}
that has following Taylor series expansion
\ShowEq{exponent Taylor series over division ring}
\end{theorem}
\begin{proof}
The G\^ateaux derivative of order $n$ has $2^n$ items.
In fact, the G\^ateaux derivative of order $1$ has $2$ items,
and each differentiation increase number of items twice.
From initial condition $y(0)=1$ and theorem
\ref{theorem: exponent derivative n over division ring}
it follows that the G\^ateaux derivative of order $n$
of required solution has form
\ShowEq{exponent derivative n over division ring, x=0}
Taylor series expansion
\eqref{eq: exponent Taylor series over division ring}
follows
from equation
\eqref{eq: exponent derivative n over division ring, x=0}.
\end{proof}

\begin{theorem}
\label{theorem: exponent of sum}
The equation
\ShowEq{exponent of sum}
is true iff
\ShowEq{exponent of sum, 1}
\end{theorem}
\begin{proof}
To prove the theorem it is enough to consider Taylor series
\ShowEq{exponent a b a+b}
Let us multiply expressions
\eqref{eq: exponent a} and \eqref{eq: exponent b}.
The sum of monomials of order $3$ has form
\ShowEq{exponent ab 3}
and in general does not equal expression
\ShowEq{exponent a+b 3}
The proof of statement that \eqref{eq: exponent of sum} follows from
\eqref{eq: exponent of sum, 1}
is trivial.
\end{proof}

The meaning of the theorem \ref{theorem: exponent of sum}
becomes more clear if we recall that there exist two
models of design of exponent. First model is
the solution of differential equation
\eqref{eq: exponent derivative over division ring}.
Second model is exploring of one parameter group of transformations.
For field both models lead to the same function.
I cannot state this now for general case. This is the subject of separate
research. However if we recall that quaternion is analogue
of transformation of three dimensional space then the statement of the theorem
becomes evident.

%% file: Differential.Eq.tex

\DefEq
{
\[
\partial f:D\rightarrow \mathcal L(D;D)
\]
}
{differential L(D,D)}

\DefEq
{
\begin{equation}
\partial (y)(h)=\pC{s}{0}f\ h\ \pC{s}{1}f
\label{eq: differential equation, additive function, division ring}
\end{equation}
}
{differential equation, additive function, division ring}

\DefEq
{
\begin{equation}
\label{eq: Gateaux differential of map df, division ring}
\partial f(x+a_2)(a_1)-\partial f(x)(a_1)
=\partial(\partial f(x)(a_1))(a_2)
+\Vector o_2(a_2)
\end{equation}
}
{Gateaux differential of map df, division ring}

\DefEq
{
$\Vector o_2:D\rightarrow \mathcal L(D;D)$
}
{Gateaux differential of map df, division ring, o2}

\DefEq
{
\[
\partial^2 f:D\rightarrow \mathcal L(D,D;D)
\]
}
{differential L(D,D;D)}

\DefEq
{
$\partial^2 f(x)\in\mathcal L(D,D;D)$%
}
{differential L(D,D;D), 1}

\DefEq
{
\[
\lim_{a_2\rightarrow 0}\frac{\|\Vector o_2(a_2)\|}{|a_2|}=0
\]
}
{Gateaux differential of map df, division ring, o2 lim}

\DefEq
{
\[
y=\pC{s}{0}f\ x\ \pC{s}{1}f
\]
}
{differential equation, additive function, solution, division ring}

\DefEq
{
\begin{equation}
\ShowSymbol{Gateaux differential of Second Order, division ring}
=\frac{\pC{s}{0}\partial^2 f(x)}{\partial x^2}
\sigma_s(a_1)
\frac{\pC{s}{1}\partial^2 f(x)}{\partial x^2}
\sigma_s(a_2)
\frac{\pC{s}{1}\partial^2 f(x)}{\partial x^2}
\label{eq: Differential of Second Order, division ring, representation}
\end{equation}
}
{Differential of Second Order, division ring, representation}

\DefEq
{
\begin{align}
\label{eq: differential equation y=xx, 2, division ring}
&\partial^2 (y)(h_1;h_2)
\\ \nonumber
=&h_1h_2x+h_1xh_2
+h_2h_1x+xh_1h_2+h_2xh_1+xh_2h_1
\\
\label{eq: differential equation y=xx, 3, division ring}
&\partial^3 (y)(h_1;h_2;h_3)
\\ \nonumber
=&h_1h_2h_3+h_1h_3h_2
+h_2h_1h_3+h_3h_1h_2+h_2h_3h_1+h_3h_2h_1
\\
\label{eq: differential equation y=xx, 4, division ring}
&\partial^n (y)(h_1;...;h_n)=0\ \ \ n>3
\end{align}
}
{differential equation y=xx, division ring}

\DefEq
{
\begin{equation}
y'=y
\label{eq: exponent over field}
\end{equation}
}
{exponent over field}

\DefEq
{
\begin{equation}
\partial(y)(h)=y'h
\label{eq: derivative over field}
\end{equation}
}
{derivative over field}

\DefEq
{
\begin{equation}
\partial(y)(h)=yh
\label{eq: exponent derivative over field}
\end{equation}
}
{exponent derivative over field}

\DefEq
{
\begin{equation}
\partial(y)(h)=\frac 12(yh+hy)
\label{eq: exponent derivative over division ring}
\end{equation}
}
{exponent derivative over division ring}

\DefEq
{
\begin{equation}
\partial^n(y)(h_1,...,h_n)=\frac 1{2^n}
\sum_{\sigma} \sigma(y)\sigma(h_1) ... \sigma(h_n)
\label{eq: exponent derivative n over division ring}
\end{equation}
}
{exponent derivative n over division ring}

\DefEq
{
\begin{equation}
\partial^n(0)(h,...,h)=1
\label{eq: exponent derivative n over division ring, x=0}
\end{equation}
}
{exponent derivative n over division ring, x=0}

\DefEq
{
\[
\sigma=
\begin{pmatrix}
y&h_1&...&h_n
\\
\sigma(y)&\sigma(h_1)&...&\sigma(h_n)
\end{pmatrix}
\]
}
{exponent derivative, transposition n}

\DefEq
{
\[
\begin{pmatrix}
y&h_1&...&h_n
\end{pmatrix}
\]
}
{exponent transposition tuple 1}

\DefEq
{
\[
\begin{pmatrix}
y&h_1&h_2&h_3
\\
h_2&y&h_3&h_1
\end{pmatrix}
\]
}
{exponent transposition tuple 3}

\DefEq
{
\begin{align*}
\begin{pmatrix}
\sigma(y)&\sigma(h_1)&\sigma(h_2)&\sigma(h_3)
\end{pmatrix}
=&
\begin{pmatrix}
h_2&y&h_3&h_1
\end{pmatrix}
\\
=&
\begin{pmatrix}
p_{\sigma}(h_2)&p_{\sigma}(y)&p_{\sigma}(h_3)&p_{\sigma}(h_1)
\end{pmatrix}
\end{align*}
}
{exponent transposition tuple 4}

\DefEq
{
\[
\begin{pmatrix}
\sigma(y)&\sigma(h_1)&...&\sigma(h_n)
\end{pmatrix}
\]
}
{exponent transposition tuple 2}

\DefEq
{
\[
\sigma=
\begin{pmatrix}
y&h_1&...&h_{k-1}
\\
\sigma(y)&\sigma(h_1)&...&\sigma(h_{k-1})
\end{pmatrix}
\]
}
{exponent derivative, transposition n=k-1}

\DefEq
{
\begin{equation}
\partial^{k-1}(y)(h_1,...,h_{k-1})=\frac 1{2^{k-1}}
\sum_{\sigma} \sigma(y)\sigma(h_1) ... \sigma(h_{k-1})
\label{eq: exponent derivative n=k-1 over division ring}
\end{equation}
}
{exponent derivative n=k-1 over division ring}

\DefEq
{
\begin{align}
\partial^k(y)(h_1,...,h_k)
=&\partial(\partial^{k-1}(y)(h_1,...,h_{k-1}))(h_k)
\nonumber
\\
=&\frac 1{2^{k-1}}\partial
\left(
\sum_{\sigma} \sigma(y)\sigma(h_1) ... \sigma(h_{k-1})
\right)
(h_k)
\label{eq: exponent derivative n=k over division ring, 1}
\end{align}
}
{exponent derivative n=k over division ring, 1}

\DefEq
{
\begin{align}
&\partial^k(y)(h_1,...,h_k)
\nonumber
\\
=&\frac 1{2^{k-1}}\frac 12
\left(
\sum_{\sigma} \sigma(yh_k)\sigma(h_1) ... \sigma(h_{k-1})
+
\sum_{\sigma} \sigma(h_ky)\sigma(h_1) ... \sigma(h_{k-1})
\right)
\label{eq: exponent derivative n=k over division ring, 2}
\end{align}
}
{exponent derivative n=k over division ring, 2}

\DefEq
{
\begin{equation}
\begin{array}{rl}
\tau_1
&=
\begin{pmatrix}
y&h_1&...&h_{k-1}&h_k
\\
\tau_1(y)&\tau_1(h_1)&...&\tau_1(h_{k-1})&\tau_1(h_k)
\end{pmatrix}
\\
&=
\begin{pmatrix}
h_ky&h_1&...&h_{k-1}&
\\
\sigma(h_ky)&\sigma(h_1)&...&\sigma(h_{k-1})
\end{pmatrix}
\\
\tau_2
&=
\begin{pmatrix}
y&h_1&...&h_{k-1}&h_k
\\
\tau_2(y)&\tau_2(h_1)&...&\tau_2(h_{k-1})&\tau_2(h_k)
\end{pmatrix}
\\
&=
\begin{pmatrix}
yh_k&h_1&...&h_{k-1}&
\\
\sigma(yh_k)&\sigma(h_1)&...&\sigma(h_{k-1})
\end{pmatrix}
\end{array}
\label{eq: exponent derivative, transposition n=k}
\end{equation}
}
{exponent derivative, transposition n=k}

\DefEq
{
\begin{align}
&\partial^k(y)(h_1,...,h_k)
\nonumber
\\
=&\frac 1{2^k}
\Bigg(
\sum_{\tau_1} \tau_1(y)\tau_1(h_1) ... \tau_1(h_{k-1})\tau_1(h_k)
\label{eq: exponent derivative n=k over division ring, 3}
\\
+&
\sum_{\tau_2} \tau_2(y)\tau_2(h_1) ... \tau_2(h_{k-1})\tau_2(h_k)
\Bigg)
\nonumber
\end{align}
}
{exponent derivative n=k over division ring, 3}

\DefEq
{
\begin{equation}
e^x=\sum_{n=0}^{\infty}
\frac 1{n!}x^n
\label{eq: exponent Taylor series over division ring}
\end{equation}
}
{exponent Taylor series over division ring}

\DefEq
{
$y=e^x$
}
{exponent over division ring}

\DefEq
{
\begin{equation}
e^{a+b}=e^ae^b
\label{eq: exponent of sum}
\end{equation}
}
{exponent of sum}

\DefEq
{
\begin{equation}
ab=ba
\label{eq: exponent of sum, 1}
\end{equation}
}
{exponent of sum, 1}

\DefEq
{
\begin{align}
e^a&=\sum_{n=0}^{\infty}\frac 1{n!}a^n
\label{eq: exponent a}
\\
e^b&=\sum_{n=0}^{\infty}\frac 1{n!}b^n
\label{eq: exponent b}
\\
e^{a+b}&=\sum_{n=0}^{\infty}\frac 1{n!}(a+b)^n
\label{eq: exponent a+b}
\end{align}
}
{exponent a b a+b}

\DefEq
{
\begin{equation}
\frac 16a^3+\frac 12a^2b+\frac 12ab^2+\frac 16b^3
\label{eq: exponent ab 3}
\end{equation}
}
{exponent ab 3}

\DefEq
{
\begin{equation}
\frac 16(a+b)^3=\frac 16a^3+\frac 16a^2b+\frac 16aba+\frac 16ba^2
+\frac 16ab^2+\frac 16bab+\frac 16b^2a+\frac 16b^3
\label{eq: exponent a+b 3}
\end{equation}
}
{exponent a+b 3}

%% file: Quaternion.English.tex
\def\texQuaternion{}
\ifx\PrintBook\undefined
\else
\chapter{Quaternion Algebra}
\label{chapter: Quaternion Algebra}
\fi

\input{Quaternion.Eq}

\section{Linear Function of Complex Field}
\label{section: Function of Complex Field}

\begin{theorem}[the Cauchy Riemann equations]
\label{theorem: complex field over real field}
Let us consider complex field $C$ as two-dimensional algebra over real field.
Let ${}_{\gi 0}\Vector e=1$, ${}_{\gi 1}\Vector e=i$ be the basis of algebra $C$.
Then in this basis structural constants have form
\[
\begin{matrix}
{}_{\gi{00}}B^{\gi 0}=1&{}_{\gi{01}}B^{\gi 1}=1\ \ 
\\
{}_{\gi{10}}B^{\gi 1}=1&{}_{\gi{\gi{11}}}B^{\gi 0}=-1
\end{matrix}
\]
Matrix of linear function
\[
y^{\gi i}=x^{\gi j}\ {}_{\gi j}f^{\gi i}
\]
of complex field over real field
satisfies relationship
\begin{align}
\label{eq: complex field over real field, 0}
{}_{\gi 0}f^{\gi 0}&={}_{\gi 1}f^{\gi 1}
\\
\label{eq: complex field over real field, 1}
{}_{\gi 0}f^{\gi 1}&=-{}_{\gi 1}f^{\gi 0}
\end{align}
\end{theorem}
\begin{proof}
Value of structural constants follows from equation $i^2=-1$.
Using equation \eqref{eq: linear map over field, division ring, relation}
we get relationships
\begin{equation}
\label{eq: complex field over real field, 0, 0}
{}_{\gi 0}f^{\gi 0}=f^{\gi{kr}}\ {}_{\gi{k0}}B^{\gi p}\ {}_{\gi{pr}}B^{\gi 0}
=f^{\gi{0r}}\ {}_{\gi{00}}B^{\gi 0}\ {}_{\gi{0r}}B^{\gi 0}+f^{\gi{1r}}\ {}_{\gi{10}}B^{\gi 1}\ {}_{\gi{1r}}B^{\gi 0}
=f^{\gi{00}}-f^{\gi{11}}
\end{equation}
\begin{equation}
\label{eq: complex field over real field, 1, 0}
{}_{\gi 0}f^{\gi 1}=f^{\gi{kr}}\ {}_{\gi{k0}}B^{\gi p}\ {}_{\gi{pr}}B^{\gi 1}
=f^{\gi{0r}}\ {}_{\gi{00}}B^{\gi 0}\ {}_{\gi{0r}}B^{\gi 1}+f^{\gi{1r}}\ {}_{\gi{10}}B^{\gi 1}\ {}_{\gi{1r}}B^{\gi 1}
=f^{\gi{01}}+f^{\gi{10}}
\end{equation}
\begin{equation}
\label{eq: complex field over real field, 1, 1}
{}_{\gi 1}f^{\gi 0}=f^{\gi{kr}}\ {}_{\gi{k1}}B^{\gi p}\ {}_{\gi{pr}}B^{\gi 0}
=f^{\gi{0r}}\ {}_{\gi{01}}B^{\gi 1}\ {}_{\gi{1r}}B^{\gi 0}+f^{\gi{1r}}\ {}_{\gi{11}}B^{\gi 0}\ {}_{\gi{0r}}B^{\gi 0}
=-f^{\gi{01}}-f^{\gi{10}}
\end{equation}
\begin{equation}
\label{eq: complex field over real field, 0, 1}
{}_{\gi 1}f^{\gi 1}=f^{\gi{kr}}\ {}_{\gi{k1}}B^{\gi p}\ {}_{\gi{pr}}B^{\gi 1}
=f^{\gi{0r}}\ {}_{\gi{01}}B^{\gi 1}\ {}_{\gi{1r}}B^{\gi 1}+f^{\gi{1r}}\ {}_{\gi{11}}B^{\gi 0}\ {}_{\gi{0r}}B^{\gi 1}
=f^{\gi{00}}-f^{\gi{11}}
\end{equation}
\eqref{eq: complex field over real field, 0}
follows from equations
\eqref{eq: complex field over real field, 0, 0} and
\eqref{eq: complex field over real field, 0, 1}.
\eqref{eq: complex field over real field, 1}
follows from equations
\eqref{eq: complex field over real field, 1, 0} and
\eqref{eq: complex field over real field, 1, 1}.
\end{proof}

\begin{remark}
\label{remark: congugate complex and additive map}
In order to show how it is hard to find generators
of additive function, let us consider equation
\eqref{eq: additive multiplicative map over field, G, division ring, relation}
for complex field.
\begin{align*}
{}_{\gi 0}f^{\gi 0}
&
={}_{\gi 0}G^{\gi l}\ f^{\gi{kr}}_G\ {}_{\gi{kl}}B^{\gi p}\ {}_{\gi{pr}}B^{\gi 0}
\\
&
={}_{\gi 0}G^{\gi 0}\ f^{\gi{00}}_G\ {}_{\gi{00}}B^{\gi 0}\ {}_{\gi{00}}B^{\gi 0}
+{}_{\gi 0}G^{\gi 1}\ f^{\gi{10}}_G\ {}_{\gi{11}}B^{\gi 0}\ {}_{\gi{00}}B^{\gi 0}
\\
&
+{}_{\gi 0}G^{\gi 1}\ f^{\gi{01}}_G\ {}_{\gi{01}}B^{\gi 1}\ {}_{\gi{11}}B^{\gi 0}
+{}_{\gi 0}G^{\gi 0}\ f^{\gi{11}}_G\ {}_{\gi{10}}B^{\gi 1}\ {}_{\gi{11}}B^{\gi 0}
\\
&
={}_{\gi 0}G^{\gi 0}\ f^{\gi{00}}_G
-{}_{\gi 0}G^{\gi 1}\ f^{\gi{10}}_G
-{}_{\gi 0}G^{\gi 1}\ f^{\gi{01}}_G
-{}_{\gi 0}G^{\gi 0}\ f^{\gi{11}}_G
\\
&
={}_{\gi 0}G^{\gi 0}(f^{\gi{00}}_G-f^{\gi{11}}_G)
-{}_{\gi 0}G^{\gi 1}(f^{\gi{10}}_G+f^{\gi{01}}_G)
\end{align*}

\begin{align*}
{}_{\gi 0}f^{\gi 1}
&
={}_{\gi 0}G^{\gi l}\ f^{\gi{kr}}_G\ {}_{\gi{kl}}B^{\gi p}\ {}_{\gi{pr}}B^{\gi 1}
\\
&
={}_{\gi 0}G^{\gi 0}\ f^{\gi{01}}_G\ {}_{\gi{00}}B^{\gi 0}\ {}_{\gi{01}}B^{\gi 1}
+{}_{\gi 0}G^{\gi 1}\ f^{\gi{11}}_G\ {}_{\gi{11}}B^{\gi 0}\ {}_{\gi{01}}B^{\gi 1}
\\
&
+{}_{\gi 0}G^{\gi 1}\ f^{\gi{00}}_G\ {}_{\gi{01}}B^{\gi 1}\ {}_{\gi{10}}B^{\gi 1}
+{}_{\gi 0}G^{\gi 0}\ f^{\gi{10}}_G\ {}_{\gi{10}}B^{\gi 1}\ {}_{\gi{10}}B^{\gi 1}
\\
&
={}_{\gi 0}G^{\gi 0}\ f^{\gi{01}}_G
-{}_{\gi 0}G^{\gi 1}\ f^{\gi{11}}_G
+{}_{\gi 0}G^{\gi 1}\ f^{\gi{00}}_G
+{}_{\gi 0}G^{\gi 0}\ f^{\gi{10}}_G
\\
&
={}_{\gi 0}G^{\gi 0}(f^{\gi{01}}_G+f^{\gi{10}}_G)
+{}_{\gi 0}G^{\gi 1}(f^{\gi{00}}_G-f^{\gi{11}}_G)
\end{align*}

\begin{align*}
{}_{\gi 1}f^{\gi 0}
&
={}_{\gi 1}G^{\gi l}\ f^{\gi{kr}}_G\ {}_{\gi{kl}}B^{\gi p}\ {}_{\gi{pr}}B^{\gi 0}
\\
&
={}_{\gi 1}G^{\gi 0}\ f^{\gi{00}}_G\ {}_{\gi{00}}B^{\gi 0}\ {}_{\gi{00}}B^{\gi 0}
+{}_{\gi 1}G^{\gi 1}\ f^{\gi{10}}_G\ {}_{\gi{11}}B^{\gi 0}\ {}_{\gi{00}}B^{\gi 0}
\\
&
+{}_{\gi 1}G^{\gi 1}\ f^{\gi{01}}_G\ {}_{\gi{01}}B^{\gi 1}\ {}_{\gi{11}}B^{\gi 0}
+{}_{\gi 1}G^{\gi 0}\ f^{\gi{11}}_G\ {}_{\gi{10}}B^{\gi 1}\ {}_{\gi{11}}B^{\gi 0}
\\
&
={}_{\gi 1}G^{\gi 0}\ f^{\gi{00}}_G
-{}_{\gi 1}G^{\gi 1}\ f^{\gi{10}}_G
-{}_{\gi 1}G^{\gi 1}\ f^{\gi{01}}_G
-{}_{\gi 1}G^{\gi 0}\ f^{\gi{11}}_G
\\
&
={}_{\gi 1}G^{\gi 0}(f^{\gi{00}}_G-f^{\gi{11}}_G)
-{}_{\gi 1}G^{\gi 1}(f^{\gi{10}}_G+f^{\gi{01}}_G)
\end{align*}

\begin{align*}
{}_{\gi 1}f^{\gi 1}
&
={}_{\gi 1}G^{\gi l}\ f^{\gi{kr}}_G\ {}_{\gi{kl}}B^{\gi p}\ {}_{\gi{pr}}B^{\gi 1}
\\
&
={}_{\gi 1}G^{\gi 0}\ f^{\gi{01}}_G\ {}_{\gi{00}}B^{\gi 0}\ {}_{\gi{01}}B^{\gi 1}
+{}_{\gi 1}G^{\gi 1}\ f^{\gi{11}}_G\ {}_{\gi{11}}B^{\gi 0}\ {}_{\gi{01}}B^{\gi 1}
\\
&
+{}_{\gi 1}G^{\gi 1}\ f^{\gi{00}}_G\ {}_{\gi{01}}B^{\gi 1}\ {}_{\gi{10}}B^{\gi 1}
+{}_{\gi 1}G^{\gi 0}\ f^{\gi{10}}_G\ {}_{\gi{10}}B^{\gi 1}\ {}_{\gi{10}}B^{\gi 1}
\\
&
={}_{\gi 1}G^{\gi 0}\ f^{\gi{01}}_G
-{}_{\gi 1}G^{\gi 1}\ f^{\gi{11}}_G
+{}_{\gi 1}G^{\gi 1}\ f^{\gi{00}}_G
+{}_{\gi 1}G^{\gi 0}\ f^{\gi{10}}_G
\\
&
={}_{\gi 1}G^{\gi 0}(f^{\gi{01}}_G+f^{\gi{10}}_G)
+{}_{\gi 1}G^{\gi 1}(f^{\gi{00}}_G-f^{\gi{11}}_G)
\end{align*}
For complex field our task becomes easier because we know
that matrix of operator $G$ has form either
\begin{equation}
G=
\begin{pmatrix}
1&0\\0&1
\end{pmatrix}
\label{eq: additive map generator, complex field 1}
\end{equation}
or
\begin{equation}
G=
\begin{pmatrix}
1&0\\0&-1
\end{pmatrix}
\label{eq: additive map generator, complex field 2}
\end{equation}
For generator
\eqref{eq: additive map generator, complex field 2}
we get equation for coordinates of transformation
\[
{}_{\gi 0}f^{\gi 0}=-{}_{\gi 1}f^{\gi 1}
\ \ \ {}_{\gi 0}f^{\gi 1}={}_{\gi 1}f^{\gi 0}
\]
\qed
\end{remark}

\ifx\texFuture\Defined
\begin{theorem}[the Cauchy Riemann equations]
Since matrix
\[
\begin{pmatrix}
\displaystyle\frac{\partial y^{\gi 0}}{\partial x^{\gi 0}}
&
\displaystyle\frac{\partial y^{\gi 0}}{\partial x^{\gi 1}} 
\\ \\
\displaystyle\frac{\partial y^{\gi 1}}{\partial x^{\gi 0}}
&
\displaystyle\frac{\partial y^{\gi 1}}{\partial x^{\gi 1}}
\end{pmatrix}
\]
is Jacobian of map of complex variable
\[
x=x^{\gi 0}+x^{\gi 1}i\rightarrow y=y^{\gi 0}(x^{\gi 0},x^{\gi 1})+y^{\gi 1}(x^{\gi 0},x^{\gi 1})i
\]
over real field,
then
\begin{align*}
\frac{\partial y^{\gi 1}}{\partial x^{\gi 0}}&=-\frac{\partial y^{\gi 0}}{\partial x^{\gi 1}}
\\
\frac{\partial y^{\gi 0}}{\partial x^{\gi 0}}&=\frac{\partial y^{\gi 1}}{\partial x^{\gi 1}}
\end{align*}
\end{theorem}
\begin{proof}
The statement of theorem is
corollary of theorem
\ifx\IntroCalculusBook\Defined
\ref{theorem: complex field over real field}.
\else
\ifx\EprintCalculus\Defined
\ref{theorem: complex field over real field}.
\else
\xRef{0701.238}{theorem: complex field over real field}.
\fi
\fi
\end{proof}

Let us consider map $y=x^2$.
To see relationship of derivative of map
\[
\begin{matrix}
\displaystyle\frac{dy}{dx}=2x
& &
\displaystyle\frac{d(y^{\gi 0}+y^{\gi 1}i)}{d(x^{\gi 0}+x^{\gi 1}i)}=2x^{\gi 0}+2x^{\gi 1}i
\end{matrix}
\]
and Jacobian of map over real field
\[
\begin{pmatrix}
2x^{\gi 0}
&
-2x^{\gi 1} 
\\
2x^{\gi 1}
&
2x^{\gi 0}
\end{pmatrix}
\]
let us consider differential
\begin{align*}
dy=&2xdx
\\
d(y^{\gi 0}+y^{\gi 1}i)=&(2x^{\gi 0}+2x^{\gi 1}i)d(x^{\gi 0}+x^{\gi 1}i)
\\
=&2x^{\gi 0}dx^{\gi 0}-2x^{\gi 1}dx^{\gi 1}+(2x^{\gi 1}dx^{\gi 0}+2x^{\gi 0}dx^{\gi 1})i
\end{align*}
\fi

\section{Linear Function of Division Ring of Quaternions}
\label{section: Linear Function of Division Ring of Quaternions}

		\begin{definition}
Let $F$ be field.
Extension field $F(i,j,k)$ is called
\AddIndex{the quaternion algebra
\symb{E(F)}1{quaternion algebra over the field}
over the field $F$}{quaternion algebra over the field}\footnote{I
follow definition from \citeBib{Izrail M. Gelfand: Quaternion Groups}.}
if multiplication in algebra $E$ is defined according to rule
\begin{equation}
\begin{array}{c|ccc}
&i&j&k\\
\hline
i&-1&k&-j\\
j&-k&-1&i\\
k&j&-i&-1\\
\end{array}
\label{eq: product of quaternions}
\end{equation}
		\qed
		\end{definition}

Elements of the algebra $E(F)$ have form
	\[
x=x^{\gi 0}+x^{\gi 1}i+x^{\gi 2}j+x^{\gi 3}k
	\]
where $x^{\gi i}\in F$, $\gi i=\gi 0$, $\gi 1$, $\gi 2$, $\gi 3$.
Quaternion
\[
\overline x=x^{\gi 0}-x^{\gi 1}i-x^{\gi 2}j-x^{\gi 3}k
\]
is called conjugate to the quaternion $x$.
We define \AddIndex{the norm of the quaternion}{norm of quaternion} $x$ using equation
	\begin{equation}
	\label{eq: norm of quaternion}
|x|^2=x\overline x=(x^{\gi 0})^2+(x^{\gi 1})^2+(x^{\gi 2})^2+(x^{\gi 3})^2
	\end{equation}
From equation \eqref{eq: norm of quaternion}, it follows that
$E(F)$ is algebra with division.\footnote{In \citeBib{Izrail M. Gelfand: Quaternion Groups}, Gelfand 
gives more general definition considering quaternion algebra $E(F,a,b)$ with
product
\[
\begin{array}{c|ccc}
&i&j&k\\
\hline
i&a&k&aj\\
j&-k&b&-bi\\
k&-aj&bi&-ab\\
\end{array}
\]
where $a$, $b\in F$, $ab\ne 0$.
However this algebra becomes division ring only when $a<0$, $b<0$.
It follows from equation
\[
x\overline x=(x^{\gi 0})^2-a(x^{\gi 1})^2-b(x^{\gi 2})^2+ab(x^{\gi 3})^2
\]
In this case we can renorm basis such that $a=-1$, $b=-1$.}
In this case inverse element has form
	\begin{equation}
	\label{eq: inverce quaternion}
x^{-1}=|x|^{-2}\overline x
	\end{equation}

\begin{theorem}
\label{theorem: Quaternion over real field}
Let us consider division ring of quaternions $E(F)$ as four\Hyph dimensional algebra over field $F$.
Let ${}_{\gi 0}\Vector e=1$, ${}_{\gi 1}\Vector e=i$, ${}_{\gi 2}\Vector e=j$, ${}_{\gi 3}\Vector e=k$
be basis of algebra $E(F)$.
Then in this basis structural constants have form
\ShowEq{structural constants, quaternions}
Standard components of additive function over field $F$
and coordinats of corresponding linear map over field $F$
satisfy relationship
$\\$
\begin{tabular}{l|r}
\begin{minipage}{170pt}
\begin{equation}
\label{eq: quaternion over real field, 1, 0}
\left\{
\begin{matrix}
{}_{\gi 0}f^{\gi 0}=f^{\gi{00}}-f^{\gi{11}}-f^{\gi{22}}-f^{\gi{33}}
\\
{}_{\gi 1}f^{\gi 1}=f^{\gi{00}}-f^{\gi{11}}+f^{\gi{22}}+f^{\gi{33}}
\\
{}_{\gi 2}f^{\gi 2}=f^{\gi{00}}+f^{\gi{11}}-f^{\gi{22}}+f^{\gi{33}}
\\
{}_{\gi 3}f^{\gi 3}=f^{\gi{00}}+f^{\gi{11}}+f^{\gi{22}}-f^{\gi{33}}
\end{matrix}
\right.
\end{equation}
\end{minipage}
&
\begin{minipage}{170pt}
\begin{equation}
\label{eq: quaternion over real field, 2, 0}
\left\{
\begin{matrix}
4\ f^{\gi{00}}={}_{\gi 0}f^{\gi 0}+{}_{\gi 1}f^{\gi 1}+{}_{\gi 2}f^{\gi 2}+{}_{\gi 3}f^{\gi 3}
\\
4\ f^{\gi{11}}=-{}_{\gi 0}f^{\gi 0}-{}_{\gi 1}f^{\gi 1}+{}_{\gi 2}f^{\gi 2}+{}_{\gi 3}f^{\gi 3}
\\
4\ f^{\gi{22}}=-{}_{\gi 0}f^{\gi 0}+{}_{\gi 1}f^{\gi 1}-{}_{\gi 2}f^{\gi 2}+{}_{\gi 3}f^{\gi 3}
\\
4\ f^{\gi{33}}=-{}_{\gi 0}f^{\gi 0}+{}_{\gi 1}f^{\gi 1}+{}_{\gi 2}f^{\gi 2}-{}_{\gi 3}f^{\gi 3}
\end{matrix}
\right.
\end{equation}
\end{minipage}
\end{tabular}
$\\$
\begin{tabular}{l|r}
\begin{minipage}{170pt}
\begin{equation}
\label{eq: quaternion over real field, 1, 1}
\left\{
\begin{matrix}
{}_{\gi 0}f^{\gi 1}=f^{\gi{01}}+f^{\gi{10}}+f^{\gi{23}}-f^{\gi{32}}
\\
{}_{\gi 1}f^{\gi 0}=-f^{\gi{01}}-f^{\gi{10}}+f^{\gi{23}}-f^{\gi{32}}
\\
{}_{\gi 2}f^{\gi 3}=-f^{\gi{01}}+f^{\gi{10}}-f^{\gi{23}}-f^{\gi{32}}
\\
{}_{\gi 3}f^{\gi 2}=f^{\gi{01}}-f^{\gi{10}}-f^{\gi{23}}-f^{\gi{32}}
\end{matrix}
\right.
\end{equation}
\end{minipage}
&
\begin{minipage}{170pt}
\begin{equation}
\label{eq: quaternion over real field, 2, 1}
\left\{
\begin{matrix}
4\ f^{\gi{10}}=-{}_{\gi 1}f^{\gi 0}+{}_{\gi 0}f^{\gi 1}-{}_{\gi 3}f^{\gi 2}+{}_{\gi 2}f^{\gi 3}
\\
4\ f^{\gi{01}}=-{}_{\gi 1}f^{\gi 0}+{}_{\gi 0}f^{\gi 1}+{}_{\gi 3}f^{\gi 2}-{}_{\gi 2}f^{\gi 3}
\\
4\ f^{\gi{32}}=-{}_{\gi 1}f^{\gi 0}-{}_{\gi 0}f^{\gi 1}-{}_{\gi 3}f^{\gi 2}-{}_{\gi 2}f^{\gi 3}
\\
4\ f^{\gi{23}}={}_{\gi 1}f^{\gi 0}+{}_{\gi 0}f^{\gi 1}-{}_{\gi 3}f^{\gi 2}-{}_{\gi 2}f^{\gi 3}
\end{matrix}
\right.
\end{equation}
\end{minipage}
\end{tabular}
$\\$
\begin{tabular}{l|r}
\begin{minipage}{170pt}
\begin{equation}
\label{eq: quaternion over real field, 1, 2}
\left\{
\begin{matrix}
{}_{\gi 0}f^{\gi 2}=f^{\gi{02}}-f^{\gi{13}}+f^{\gi{20}}+f^{\gi{31}}
\\
{}_{\gi 1}f^{\gi 3}=f^{\gi{02}}-f^{\gi{13}}-f^{\gi{20}}-f^{\gi{31}}
\\
{}_{\gi 2}f^{\gi 0}=-f^{\gi{02}}-f^{\gi{13}}-f^{\gi{20}}+f^{\gi{31}}
\\
{}_{\gi 3}f^{\gi 1}=-f^{\gi{02}}-f^{\gi{13}}+f^{\gi{20}}-f^{\gi{31}}
\end{matrix}
\right.
\end{equation}
\end{minipage}
&
\begin{minipage}{170pt}
\begin{equation}
\label{eq: quaternion over real field, 2, 2}
\left\{
\begin{matrix}
4\ f^{\gi{20}}=-{}_{\gi 2}f^{\gi 0}+{}_{\gi 3}f^{\gi 1}+{}_{\gi 0}f^{\gi 2}-{}_{\gi 1}f^{\gi 3}
\\
4\ f^{\gi{31}}=+{}_{\gi 2}f^{\gi 0}-{}_{\gi 3}f^{\gi 1}+{}_{\gi 0}f^{\gi 2}-{}_{\gi 1}f^{\gi 3}
\\
4\ f^{\gi{02}}=-{}_{\gi 2}f^{\gi 0}-{}_{\gi 3}f^{\gi 1}+{}_{\gi 0}f^{\gi 2}+{}_{\gi 1}f^{\gi 3}
\\
4\ f^{\gi{13}}=-{}_{\gi 2}f^{\gi 0}-{}_{\gi 3}f^{\gi 1}-{}_{\gi 0}f^{\gi 2}-{}_{\gi 1}f^{\gi 3}
\end{matrix}
\right.
\end{equation}
\end{minipage}
\end{tabular}
$\\$
\begin{tabular}{l|r}
\begin{minipage}{170pt}
\begin{equation}
\label{eq: quaternion over real field, 1, 3}
\left\{
\begin{matrix}
{}_{\gi 0}f^{\gi 3}=f^{\gi{03}}+f^{\gi{12}}-f^{\gi{21}}+f^{\gi{30}}
\\
{}_{\gi 1}f^{\gi 2}=-f^{\gi{03}}-f^{\gi{12}}-f^{\gi{21}}+f^{\gi{30}}
\\
{}_{\gi 2}f^{\gi 1}=f^{\gi{03}}-f^{\gi{12}}-f^{\gi{21}}-f^{\gi{30}}
\\
{}_{\gi 3}f^{\gi 0}=-f^{\gi{03}}+f^{\gi{12}}-f^{\gi{21}}-f^{\gi{30}}
\end{matrix}
\right.
\end{equation}
\end{minipage}
&
\begin{minipage}{170pt}
\begin{equation}
\label{eq: quaternion over real field, 2, 3}
\left\{
\begin{matrix}
4\ f^{\gi{30}}=-{}_{\gi 3}f^{\gi 0}-{}_{\gi 2}f^{\gi 1}+{}_{\gi 1}f^{\gi 2}+{}_{\gi 0}f^{\gi 3}
\\
4\ f^{\gi{21}}=-{}_{\gi 3}f^{\gi 0}-{}_{\gi 2}f^{\gi 1}-{}_{\gi 1}f^{\gi 2}-{}_{\gi 0}f^{\gi 3}
\\
4\ f^{\gi{12}}={}_{\gi 3}f^{\gi 0}-{}_{\gi 2}f^{\gi 1}-{}_{\gi 1}f^{\gi 2}+{}_{\gi 0}f^{\gi 3}
\\
4\ f^{\gi{03}}=-{}_{\gi 3}f^{\gi 0}+{}_{\gi 2}f^{\gi 1}-{}_{\gi 1}f^{\gi 2}+{}_{\gi 0}f^{\gi 3}
\end{matrix}
\right.
\end{equation}
\end{minipage}
\end{tabular}
\end{theorem}
\begin{proof}
Value of structural constants follows from multiplication table
\eqref{eq: product of quaternions}.
Using equation \eqref{eq: linear map over field, division ring, relation}
we get relationships
\begin{align*}
{}_{\gi 0}f^{\gi 0}&=f^{\gi{kr}}\ {}_{\gi{k0}}B^{\gi p}\ {}_{\gi{pr}}B^{\gi 0}
\\
&=f^{\gi{00}}\ {}_{\gi{00}}B^{\gi 0}\ {}_{\gi{00}}B^{\gi 0}+f^{\gi{11}}\ {}_{\gi{10}}B^{\gi 1}\ {}_{\gi{11}}B^{\gi 0}
+f^{\gi{22}}\ {}_{\gi{20}}B^{\gi 2}\ {}_{\gi{22}}B^{\gi 0}+f^{\gi{33}}\ {}_{\gi{30}}B^{\gi 3}\ {}_{\gi{33}}B^{\gi 0}
\\
&=f^{\gi{00}}-f^{\gi{11}}
-f^{\gi{22}}-f^{\gi{33}}
\end{align*}
\begin{align*}
{}_{\gi 0}f^{\gi 1}&=f^{\gi{kr}}\ {}_{\gi{k0}}B^{\gi p}\ {}_{\gi{pr}}B^{\gi 1}
\\
&=f^{\gi{01}}\ {}_{\gi{00}}B^{\gi 0}\ {}_{\gi{01}}B^{\gi 1}+f^{\gi{10}}\ {}_{\gi{10}}B^{\gi 1}\ {}_{\gi{10}}B^{\gi 1}
+f^{\gi{23}}\ {}_{\gi{20}}B^{\gi 2}\ {}_{\gi{23}}B^{\gi 1}+f^{\gi{32}}\ {}_{\gi{30}}B^{\gi 3}\ {}_{\gi{32}}B^{\gi 1}
\\
&=f^{\gi{01}}+f^{\gi{10}}
+f^{\gi{23}}-f^{\gi{32}}
\end{align*}
\begin{align*}
{}_{\gi 0}f^{\gi 2}&=f^{\gi{kr}}\ {}_{\gi{k0}}B^{\gi p}\ {}_{\gi{pr}}B^{\gi 2}
\\
&=f^{\gi{02}}\ {}_{\gi{00}}B^{\gi 0}\ {}_{\gi{02}}B^{\gi 2}+f^{\gi{13}}\ {}_{\gi{10}}B^{\gi 1}\ {}_{\gi{13}}B^{\gi 2}
+f^{\gi{20}}\ {}_{\gi{20}}B^{\gi 2}\ {}_{\gi{20}}B^{\gi 2}+f^{\gi{31}}\ {}_{\gi{30}}B^{\gi 3}\ {}_{\gi{31}}B^{\gi 2}
\\
&=f^{\gi{02}}-f^{\gi{13}}
+f^{\gi{20}}+f^{\gi{31}}
\end{align*}
\begin{align*}
{}_{\gi 0}f^{\gi 3}&=f^{\gi{kr}}\ {}_{\gi{k0}}B^{\gi p}\ {}_{\gi{pr}}B^{\gi 3}
\\
&=f^{\gi{03}}\ {}_{\gi{00}}B^{\gi 0}\ {}_{\gi{03}}B^{\gi 3}+f^{\gi{12}}\ {}_{\gi{10}}B^{\gi 1}\ {}_{\gi{12}}B^{\gi 3}
+f^{\gi{21}}\ {}_{\gi{20}}B^{\gi 2}\ {}_{\gi{21}}B^{\gi 3}+f^{\gi{30}}\ {}_{\gi{30}}B^{\gi 3}\ {}_{\gi{30}}B^{\gi 3}
\\
&=f^{\gi{03}}+f^{\gi{12}}
-f^{\gi{21}}+f^{\gi{30}}
\end{align*}
\begin{align*}
{}_{\gi 1}f^{\gi 0}&=f^{\gi{kr}}\ {}_{\gi{k1}}B^{\gi p}\ {}_{\gi{pr}}B^{\gi 0}
\\
&=f^{\gi{01}}\ {}_{\gi{01}}B^{\gi 1}\ {}_{\gi{11}}B^{\gi 0}+f^{\gi{10}}\ {}_{\gi{11}}B^{\gi 0}\ {}_{\gi{00}}B^{\gi 0}
+f^{\gi{23}}\ {}_{\gi{21}}B^{\gi 3}\ {}_{\gi{33}}B^{\gi 0}+f^{\gi{32}}\ {}_{\gi{31}}B^{\gi 2}\ {}_{\gi{22}}B^{\gi 0}
\\
&=-f^{\gi{01}}-f^{\gi{10}}
+f^{\gi{23}}-f^{\gi{32}}
\end{align*}
\begin{align*}
{}_{\gi 1}f^{\gi 1}&=f^{\gi{kr}}\ {}_{\gi{k1}}B^{\gi p}\ {}_{\gi{pr}}B^{\gi 1}
\\
&=f^{\gi{00}}\ {}_{\gi{01}}B^{\gi 1}\ {}_{\gi{10}}B^{\gi 1}+f^{\gi{11}}\ {}_{\gi{11}}B^{\gi 0}\ {}_{\gi{01}}B^{\gi 1}
+f^{\gi{22}}\ {}_{\gi{21}}B^{\gi 3}\ {}_{\gi{32}}B^{\gi 1}+f^{\gi{33}}\ {}_{\gi{31}}B^{\gi 2}\ {}_{\gi{23}}B^{\gi 1}
\\
&=f^{\gi{00}}-f^{\gi{11}}
+f^{\gi{22}}+f^{\gi{33}}
\end{align*}
\begin{align*}
{}_{\gi 1}f^{\gi 2}&=f^{\gi{kr}}\ {}_{\gi{k1}}B^{\gi p}\ {}_{\gi{pr}}B^{\gi 2}
\\
&=f^{\gi{03}}\ {}_{\gi{01}}B^{\gi 1}\ {}_{\gi{13}}B^{\gi 2}+f^{\gi{12}}\ {}_{\gi{11}}B^{\gi 0}\ {}_{\gi{02}}B^{\gi 2}
+f^{\gi{21}}\ {}_{\gi{21}}B^{\gi 3}\ {}_{\gi{31}}B^{\gi 2}+f^{\gi{30}}\ {}_{\gi{31}}B^{\gi 2}\ {}_{\gi{20}}B^{\gi 2}
\\
&=-f^{\gi{03}}-f^{\gi{12}}
-f^{\gi{21}}+f^{\gi{30}}
\end{align*}
\begin{align*}
{}_{\gi 1}f^{\gi 3}&=f^{\gi{kr}}\ {}_{\gi{k1}}B^{\gi p}\ {}_{\gi{pr}}B^{\gi 3}
\\
&=f^{\gi{02}}\ {}_{\gi{01}}B^{\gi 1}\ {}_{\gi{12}}B^{\gi 3}+f^{\gi{13}}\ {}_{\gi{11}}B^{\gi 0}\ {}_{\gi{03}}B^{\gi 3}
+f^{\gi{20}}\ {}_{\gi{21}}B^{\gi 3}\ {}_{\gi{30}}B^{\gi 3}+f^{\gi{31}}\ {}_{\gi{31}}B^{\gi 2}\ {}_{\gi{21}}B^{\gi 3}
\\
&=f^{\gi{02}}-f^{\gi{13}}
-f^{\gi{20}}-f^{\gi{31}}
\end{align*}
\begin{align*}
{}_{\gi 2}f^{\gi 0}&=f^{\gi{kr}}\ {}_{\gi{k2}}B^{\gi p}\ {}_{\gi{pr}}B^{\gi 0}
\\
&=f^{\gi{02}}\ {}_{\gi{02}}B^{\gi 2}\ {}_{\gi{22}}B^{\gi 0}+f^{\gi{13}}\ {}_{\gi{12}}B^{\gi 3}\ {}_{\gi{33}}B^{\gi 0}
+f^{\gi{20}}\ {}_{\gi{22}}B^{\gi 0}\ {}_{\gi{00}}B^{\gi 0}+f^{\gi{31}}\ {}_{\gi{32}}B^{\gi 1}\ {}_{\gi{11}}B^{\gi 0}
\\
&=-f^{\gi{02}}-f^{\gi{13}}
-f^{\gi{20}}+f^{\gi{31}}
\end{align*}
\begin{align*}
{}_{\gi 2}f^{\gi 1}&=f^{\gi{kr}}\ {}_{\gi{k2}}B^{\gi p}\ {}_{\gi{pr}}B^{\gi 1}
\\
&=f^{\gi{03}}\ {}_{\gi{02}}B^{\gi 2}\ {}_{\gi{23}}B^{\gi 1}+f^{\gi{12}}\ {}_{\gi{12}}B^{\gi 3}\ {}_{\gi{32}}B^{\gi 1}
+f^{\gi{21}}\ {}_{\gi{22}}B^{\gi 0}\ {}_{\gi{01}}B^{\gi 1}+f^{\gi{30}}\ {}_{\gi{32}}B^{\gi 1}\ {}_{\gi{10}}B^{\gi 1}
\\
&=f^{\gi{03}}-f^{\gi{12}}
-f^{\gi{21}}-f^{\gi{30}}
\end{align*}
\begin{align*}
{}_{\gi 2}f^{\gi 2}&=f^{\gi{kr}}\ {}_{\gi{k2}}B^{\gi p}\ {}_{\gi{pr}}B^{\gi 2}
\\
&=f^{\gi{00}}\ {}_{\gi{02}}B^{\gi 2}\ {}_{\gi{20}}B^{\gi 2}+f^{\gi{11}}\ {}_{\gi{12}}B^{\gi 3}\ {}_{\gi{31}}B^{\gi 2}
+f^{\gi{22}}\ {}_{\gi{22}}B^{\gi 0}\ {}_{\gi{02}}B^{\gi 2}+f^{\gi{33}}\ {}_{\gi{32}}B^{\gi 1}\ {}_{\gi{13}}B^{\gi 2}
\\
&=f^{\gi{00}}+f^{\gi{11}}
-f^{\gi{22}}+f^{\gi{33}}
\end{align*}
\begin{align*}
{}_{\gi 2}f^{\gi 3}&=f^{\gi{kr}}\ {}_{\gi{k2}}B^{\gi p}\ {}_{\gi{pr}}B^{\gi 3}
\\
&=f^{\gi{01}}\ {}_{\gi{02}}B^{\gi 2}\ {}_{\gi{21}}B^{\gi 3}+f^{\gi{10}}\ {}_{\gi{12}}B^{\gi 3}\ {}_{\gi{30}}B^{\gi 3}
+f^{\gi{23}}\ {}_{\gi{22}}B^{\gi 0}\ {}_{\gi{03}}B^{\gi 3}+f^{\gi{32}}\ {}_{\gi{32}}B^{\gi 1}\ {}_{\gi{12}}B^{\gi 3}
\\
&=-f^{\gi{01}}+f^{\gi{10}}
-f^{\gi{23}}-f^{\gi{32}}
\end{align*}
\begin{align*}
{}_{\gi 3}f^{\gi 0}&=f^{\gi{kr}}\ {}_{\gi{k3}}B^{\gi p}\ {}_{\gi{pr}}B^{\gi 0}
\\
&=f^{\gi{03}}\ {}_{\gi{03}}B^{\gi 3}\ {}_{\gi{33}}B^{\gi 0}+f^{\gi{12}}\ {}_{\gi{13}}B^{\gi 2}\ {}_{\gi{22}}B^{\gi 0}
+f^{\gi{21}}\ {}_{\gi{23}}B^{\gi 1}\ {}_{\gi{11}}B^{\gi 0}+f^{\gi{30}}\ {}_{\gi{33}}B^{\gi 0}\ {}_{\gi{00}}B^{\gi 0}
\\
&=-f^{\gi{03}}+f^{\gi{12}}
-f^{\gi{21}}-f^{\gi{30}}
\end{align*}
\begin{align*}
{}_{\gi 3}f^{\gi 1}&=f^{\gi{kr}}\ {}_{\gi{k3}}B^{\gi p}\ {}_{\gi{pr}}B^{\gi 1}
\\
&=f^{\gi{02}}\ {}_{\gi{03}}B^{\gi 3}\ {}_{\gi{32}}B^{\gi 1}+f^{\gi{13}}\ {}_{\gi{13}}B^{\gi 2}\ {}_{\gi{23}}B^{\gi 1}
+f^{\gi{20}}\ {}_{\gi{23}}B^{\gi 1}\ {}_{\gi{10}}B^{\gi 1}+f^{\gi{31}}\ {}_{\gi{33}}B^{\gi 0}\ {}_{\gi{01}}B^{\gi 1}
\\
&=-f^{\gi{02}}-f^{\gi{13}}
+f^{\gi{20}}-f^{\gi{31}}
\end{align*}
\begin{align*}
{}_{\gi 3}f^{\gi 2}&=f^{\gi{kr}}\ {}_{\gi{k3}}B^{\gi p}\ {}_{\gi{pr}}B^{\gi 2}
\\
&=f^{\gi{01}}\ {}_{\gi{03}}B^{\gi 3}\ {}_{\gi{31}}B^{\gi 2}+f^{\gi{10}}\ {}_{\gi{13}}B^{\gi 2}\ {}_{\gi{20}}B^{\gi 2}
+f^{\gi{23}}\ {}_{\gi{23}}B^{\gi 1}\ {}_{\gi{13}}B^{\gi 2}+f^{\gi{32}}\ {}_{\gi{33}}B^{\gi 0}\ {}_{\gi{02}}B^{\gi 2}
\\
&=f^{\gi{01}}-f^{\gi{10}}
-f^{\gi{23}}-f^{\gi{32}}
\end{align*}
\begin{align*}
{}_{\gi 3}f^{\gi 3}&=f^{\gi{kr}}\ {}_{\gi{k3}}B^{\gi p}\ {}_{\gi{pr}}B^{\gi 3}
\\
&=f^{\gi{00}}\ {}_{\gi{03}}B^{\gi 3}\ {}_{\gi{30}}B^{\gi 3}+f^{\gi{11}}\ {}_{\gi{13}}B^{\gi 2}\ {}_{\gi{21}}B^{\gi 3}
+f^{\gi{22}}\ {}_{\gi{23}}B^{\gi 1}\ {}_{\gi{12}}B^{\gi 3}+f^{\gi{33}}\ {}_{\gi{33}}B^{\gi 0}\ {}_{\gi{03}}B^{\gi 3}
\\
&=f^{\gi{00}}+f^{\gi{11}}
+f^{\gi{22}}-f^{\gi{33}}
\end{align*}

We group these relationships into systems of linear equations
\eqref{eq: quaternion over real field, 1, 0},
\eqref{eq: quaternion over real field, 1, 1},
\eqref{eq: quaternion over real field, 1, 2},
\eqref{eq: quaternion over real field, 1, 3}.

\eqref{eq: quaternion over real field, 2, 0} is
solution of system of linear equations
\eqref{eq: quaternion over real field, 1, 0}.

\eqref{eq: quaternion over real field, 2, 1} is
solution of system of linear equations
\eqref{eq: quaternion over real field, 1, 1}.

\eqref{eq: quaternion over real field, 2, 2} is
solution of system of linear equations
\eqref{eq: quaternion over real field, 1, 2}.

\eqref{eq: quaternion over real field, 2, 3} is
solution of system of linear equations
\eqref{eq: quaternion over real field, 1, 3}.
\end{proof}

\begin{theorem}
\label{theorem: Quaternion over real field, matrix}
Let us consider division ring of quaternions $E(F)$ as four\Hyph dimensional algebra over field $F$.
Let ${}_{\gi 0}\Vector e=1$, ${}_{\gi 1}\Vector e=i$, ${}_{\gi 2}\Vector e=j$, ${}_{\gi 3}\Vector e=k$ be basis of algebra $E(F)$.
Standard components of additive function over field $F$
and coordinats of this function over field $F$
satisfy relationship
\begin{align}
&
\begin{pmatrix}
{}_{\gi 0}f^{\gi 0}&{}_{\gi 0}f^{\gi 1}&{}_{\gi 0}f^{\gi 2}&{}_{\gi 0}f^{\gi 3}
\\
{}_{\gi 1}f^{\gi 1}&{}_{\gi 1}f^{\gi 0}&{}_{\gi 1}f^{\gi 3}&{}_{\gi 1}f^{\gi 2}
\\
{}_{\gi 2}f^{\gi 2}&{}_{\gi 2}f^{\gi 3}&{}_{\gi 2}f^{\gi 0}&{}_{\gi 2}f^{\gi 1}
\\
{}_{\gi 3}f^{\gi 3}&{}_{\gi 3}f^{\gi 2}&{}_{\gi 3}f^{\gi 1}&{}_{\gi 3}f^{\gi 0}
\end{pmatrix}
\nonumber
\\
=&
\begin{pmatrix}
1&-1&-1&-1
\\
1&-1&1&1
\\
1&1&-1&1
\\
1&1&1&-1
\end{pmatrix}
\begin{pmatrix}
f^{\gi{00}}&-f^{\gi{32}}&-f^{\gi{13}}&-f^{\gi{21}}
\\
f^{\gi{11}}&-f^{\gi{23}}&-f^{\gi{02}}&-f^{\gi{30}}
\\
f^{\gi{22}}&-f^{\gi{10}}&-f^{\gi{31}}&-f^{\gi{03}}
\\
f^{\gi{33}}&-f^{\gi{01}}&-f^{\gi{20}}&-f^{\gi{12}}
\end{pmatrix}
\label{eq: quaternion over real field, 3}
\end{align}
\end{theorem}
\begin{proof}
Let us write equation \eqref{eq: quaternion over real field, 1, 0}
as product of matrices
\begin{equation}
\label{eq: quaternion over real field, 3, 0}
\begin{pmatrix}
{}_{\gi 0}f^{\gi 0}
\\
{}_{\gi 1}f^{\gi 1}
\\
{}_{\gi 2}f^{\gi 2}
\\
{}_{\gi 3}f^{\gi 3}
\end{pmatrix}
=
\begin{pmatrix}
1&-1&-1&-1
\\
1&-1&1&1
\\
1&1&-1&1
\\
1&1&1&-1
\end{pmatrix}
\begin{pmatrix}
f^{\gi{00}}
\\
f^{\gi{11}}
\\
f^{\gi{22}}
\\
f^{\gi{33}}
\end{pmatrix}
\end{equation}
Let us write equation \eqref{eq: quaternion over real field, 1, 1}
as product of matrices
\begin{align}
\begin{pmatrix}
{}_{\gi 0}f^{\gi 1}
\\
{}_{\gi 1}f^{\gi 0}
\\
{}_{\gi 2}f^{\gi 3}
\\
{}_{\gi 3}f^{\gi 2}
\end{pmatrix}
&=
\begin{pmatrix}
1&1&1&-1
\\
-1&-1&1&-1
\\
-1&1&-1&-1
\\
1&-1&-1&-1
\end{pmatrix}
\begin{pmatrix}
f^{\gi{01}}
\\
f^{\gi{10}}
\\
f^{\gi{23}}
\\
f^{\gi{32}}
\end{pmatrix}
\nonumber
\\
&=
\begin{pmatrix}
-1&-1&-1&1
\\
1&1&-1&1
\\
1&-1&1&1
\\
-1&1&1&1
\end{pmatrix}
\begin{pmatrix}
{}-f^{\gi{01}}
\\
-f^{\gi{10}}
\\
-f^{\gi{23}}
\\
-f^{\gi{32}}
\end{pmatrix}
\label{eq: quaternion over real field, 3, 1}
\\
&=
\begin{pmatrix}
1&-1&-1&-1
\\
1&-1&1&1
\\
1&1&-1&1
\\
1&1&1&-1
\end{pmatrix}
\begin{pmatrix}
-f^{\gi{32}}
\\
-f^{\gi{23}}
\\
-f^{\gi{10}}
\\
-f^{\gi{01}}
\end{pmatrix}
\nonumber
\end{align}
Let us write equation \eqref{eq: quaternion over real field, 1, 2}
as product of matrices
\begin{align}
\begin{pmatrix}
{}_{\gi 0}f^{\gi 2}
\\
{}_{\gi 1}f^{\gi 3}
\\
{}_{\gi 2}f^{\gi 0}
\\
{}_{\gi 3}f^{\gi 1}
\end{pmatrix}
&=
\begin{pmatrix}
1&-1&1&1
\\
1&-1&-1&-1
\\
-1&-1&-1&1
\\
-1&-1&1&-1
\end{pmatrix}
\begin{pmatrix}
f^{\gi{02}}
\\
f^{\gi{13}}
\\
f^{\gi{20}}
\\
f^{\gi{31}}
\end{pmatrix}
\nonumber
\\
&=
\begin{pmatrix}
-1&1&-1&-1
\\
-1&1&1&1
\\
1&1&1&-1
\\
1&1&-1&1
\end{pmatrix}
\begin{pmatrix}
-f^{\gi{02}}
\\
-f^{\gi{13}}
\\
-f^{\gi{20}}
\\
-f^{\gi{31}}
\end{pmatrix}
\label{eq: quaternion over real field, 3, 2}
\\
&=
\begin{pmatrix}
1&-1&-1&-1
\\
1&-1&1&1
\\
1&1&-1&1
\\
1&1&1&-1
\end{pmatrix}
\begin{pmatrix}
-f^{\gi{13}}
\\
-f^{\gi{02}}
\\
-f^{\gi{31}}
\\
-f^{\gi{20}}
\end{pmatrix}
\nonumber
\end{align}
Let us write equation \eqref{eq: quaternion over real field, 1, 3}
as product of matrices
\begin{align}
\begin{pmatrix}
{}_{\gi 0}f^{\gi 3}
\\
{}_{\gi 1}f^{\gi 2}
\\
{}_{\gi 2}f^{\gi 1}
\\
{}_{\gi 3}f^{\gi 0}
\end{pmatrix}
&=
\begin{pmatrix}
1&1&-1&1
\\
-1&-1&-1&1
\\
1&-1&-1&-1
\\
-1&1&-1&-1
\end{pmatrix}
\begin{pmatrix}
f^{\gi{03}}
\\
f^{\gi{12}}
\\
f^{\gi{21}}
\\
f^{\gi{30}}
\end{pmatrix}
\nonumber
\\
&=
\begin{pmatrix}
-1&-1&1&-1
\\
1&1&1&-1
\\
-1&1&1&1
\\
1&-1&1&1
\end{pmatrix}
\begin{pmatrix}
-f^{\gi{03}}
\\
-f^{\gi{12}}
\\
-f^{\gi{21}}
\\
-f^{\gi{30}}
\end{pmatrix}
\label{eq: quaternion over real field, 3, 3}
\\
&=
\begin{pmatrix}
1&-1&-1&-1
\\
1&-1&1&1
\\
1&1&-1&1
\\
1&1&1&-1
\end{pmatrix}
\begin{pmatrix}
-f^{\gi{21}}
\\
-f^{\gi{30}}
\\
-f^{\gi{03}}
\\
-f^{\gi{12}}
\end{pmatrix}
\nonumber
\end{align}
We join equations
\eqref{eq: quaternion over real field, 3, 0},
\eqref{eq: quaternion over real field, 3, 1},
\eqref{eq: quaternion over real field, 3, 2},
\eqref{eq: quaternion over real field, 3, 3}
into equation
\eqref{eq: quaternion over real field, 3}.
\end{proof}

\ifx\texFuture\Defined
\section{Maps between different algebras}

We use notation $H=E(R)$ for quaternion algebra over real field.

\begin{theorem}
Линейная функция
\ShowEq{function from R to H}
удовлетворяет равенству
\end{theorem}
\begin{proof}
Над полем действительных чисел
мы можем представить функцию \eqref{eq: function from R to H}
в виде матрицы
\ShowEq{function from R to H, coordinates}
Если в качестве генератора мы выберем отображение
\ShowEq{function from R to H, generator}
то равенство
\eqref{eq: additive multiplicative map over field, G, division ring, relation}
принимает вид
\ShowEq{function from R to H, relation}
\end{proof}
\fi

\ifx\texDifferential\Defined
\section{Differentiable Map of Division Ring of Quaternions}

\begin{theorem}
Since matrix
$\displaystyle
\begin{pmatrix}
\displaystyle\frac{\partial y^{\gi i}}{\partial x^{\gi j}}
\end{pmatrix}
$
is Jacobian of map
$
x\rightarrow y
$
of division ring of quaternions over real field,
then
\ShowEq{quaternion over real field, derivative}
\end{theorem}
\begin{proof}
The statement of theorem is
corollary of theorem
\ifx\IntroCalculusBook\Defined
\ref{theorem: Quaternion over real field}.
\else
\ifx\EprintCalculus\Defined
\ref{theorem: Quaternion over real field}.
\else
\xRef{0701.238}{theorem: Quaternion over real field}.
\fi
\fi
\end{proof}

\begin{theorem}
\label{theorem: quaternion conjugation, derivative}
Quaternionic map
\[
f(x)=\overline x
\]
has the G\^ateaux derivative
\begin{equation}
\label{eq: quaternion conjugation, derivative}
\partial (\overline x)(h)=-\frac 12(h+ihi+jhj+khk)
\end{equation}
\end{theorem}
\begin{proof}
Jacobian of the map $f$ has form
\[
\begin{pmatrix}
1&0&0&0
\\
0&-1&0&0
\\
0&0&-1&0
\\
0&0&0&-1
\end{pmatrix}
\]
From equations
\eqref{eq: quaternion over real field, derivative, 2, 0}
it follows that
\begin{equation}
\label{eq: quaternion conjugation, derivative, 1}
\StandPartial{y}{x}{00}=
\StandPartial{y}{x}{11}=
\StandPartial{y}{x}{22}=
\StandPartial{y}{x}{33}=
-\frac 12
\end{equation}
From equations
\eqref{eq: quaternion over real field, derivative, 2, 1},
\eqref{eq: quaternion over real field, derivative, 2, 2},
\eqref{eq: quaternion over real field, derivative, 2, 3}
it follows that
\begin{equation}
\label{eq: quaternion conjugation, derivative, 2}
\StandPartial{y}{x}{ij}=0
\ \ \ \gi i\ne\gi j
\end{equation}
Equation
\eqref{eq: quaternion conjugation, derivative}
follows from equations
\eqref{eq: Gateaux differential, division ring, standard representation},
\eqref{eq: quaternion conjugation, derivative, 1},
\eqref{eq: quaternion conjugation, derivative, 2}.
\end{proof}

\begin{theorem}
Quaternion conjugation satisfies equation
\[
\overline x=-\frac 12(x+ixi+jxj+kxk)
\]
\end{theorem}
\begin{proof}
The statement of theorem follows from theorem
\ref{theorem: quaternion conjugation, derivative}
and example \ref{example: differential equation, additive function, division ring}.
\end{proof}

\begin{theorem}
Since matrix
$\displaystyle
\begin{pmatrix}
\displaystyle\frac{\partial y^{\gi i}}{\partial x^{\gi j}}
\end{pmatrix}
$
is Jacobian of map
$
x\rightarrow y
$
of division ring of quaternions over real field,
then
\ShowEq{quaternion over real field, derivative, 1}
\end{theorem}
\begin{proof}
The statement of theorem is
corollary of theorem
\ifx\IntroCalculusBook\Defined
\ref{theorem: Quaternion over real field, matrix}.
\else
\ifx\EprintCalculus\Defined
\ref{theorem: Quaternion over real field, matrix}.
\else
\xRef{0701.238}{theorem: Quaternion over real field, matrix}.
\fi
\fi
\end{proof}
\fi

%% file: Quaternion.Eq.tex

\DefEq
{
\begin{equation}
f:R\rightarrow H
\label{eq: function from R to H}
\end{equation}
}
{function from R to H}

\DefEq
{
\begin{equation}
f=
\begin{pmatrix}
f^{\gi 0}&f^{\gi 1}&f^{\gi 2}&f^{\gi 3}
\end{pmatrix}
\label{eq: function from R to H, coordinates}
\end{equation}
}
{function from R to H, coordinates}

\DefEq
{
\begin{equation}
G=
\begin{pmatrix}
G^{\gi 0}&G^{\gi 1}&G^{\gi 2}&G^{\gi 3}
\end{pmatrix}
\label{eq: function from R to H, generator}
\end{equation}
}
{function from R to H, generator}

\DefEq
{
\begin{equation}
f^{\gi j}=&G^{\gi l}\ f^{\gi{kr}}_G\ {}_{\gi{kl}}B^{\gi p}
\ {}_{\gi{pr}}B^{\gi j}
\label{eq: function from R to H, relation}
\end{equation}
}
{function from R to H, relation}

\DefEq
{
\begin{equation}
f(x)=\pC{s}{0}f^{\gi i}\ {}_{\gi i}\Vector e
\ x\ \pC{s}{1}f^{\gi j}\ {}_{\gi j}\Vector e
\label{eq: additive map, division ring, standard representation, 1}
\end{equation}
}
{additive map, division ring, standard representation, 1}

\DefEq
{
\[
f^{\gi i\gi j}=\pC{s}{0}f^{\gi i}\ \pC{s}{1}f^{\gi j}
\]
}
{additive map, division ring, standard representation, 2}

\DefEq
{
\begin{alignat*}{8}
{}_{\gi{00}}B^{\gi 0}&=&1&\ \ &{}_{\gi{01}}B^{\gi 1}&=&1&\ \ &{}_{\gi{02}}B^{\gi 2}&=&1&\ \ &{}_{\gi{03}}B^{\gi 3}&=&1
\\
{}_{\gi{10}}B^{\gi 1}&=&1&&{}_{\gi{11}}B^{\gi 0}&=&-1&&{}_{\gi{12}}B^{\gi 3}&=&1&&{}_{\gi{13}}B^{\gi 2}&=&-1
\\
{}_{\gi{20}}B^{\gi 2}&=&1&&{}_{\gi{21}}B^{\gi 3}&=&-1&&{}_{\gi{22}}B^{\gi 0}&=&-1&&{}_{\gi{23}}B^{\gi 1}&=&1
\\
{}_{\gi{30}}B^{\gi 3}&=&1&&{}_{\gi{31}}B^{\gi 2}&=&1&&{}_{\gi{32}}B^{\gi 1}&=&-1&&{}_{\gi{33}}B^{\gi 0}&=&-1
\end{alignat*}
}
{structural constants, quaternions}

\DefEq
{
\begin{equation}
\label{eq: quaternion over real field, derivative, 1, 0}
\left\{
\begin{matrix}
\displaystyle
\frac{\partial y^{\gi 0}}{\partial x^{\gi 0}}=\StandPartial{y}{x}{00}-\StandPartial{y}{x}{11}-\StandPartial{y}{x}{22}-\StandPartial{y}{x}{33}
\\
\displaystyle
\vphantom{\overset{\rightarrow}{\frac{\partial}{\partial}}}
\frac{\partial y^{\gi 1}}{\partial x^{\gi 1}}=\StandPartial{y}{x}{00}-\StandPartial{y}{x}{11}+\StandPartial{y}{x}{22}+\StandPartial{y}{x}{33}
\\
\displaystyle
\vphantom{\overset{\rightarrow}{\frac{\partial}{\partial}}}
\frac{\partial y^{\gi 2}}{\partial x^{\gi 2}}=\StandPartial{y}{x}{00}+\StandPartial{y}{x}{11}-\StandPartial{y}{x}{22}+\StandPartial{y}{x}{33}
\\
\displaystyle
\vphantom{\overset{\rightarrow}{\frac{\partial}{\partial}}}
\frac{\partial y^{\gi 3}}{\partial x^{\gi 3}}=\StandPartial{y}{x}{00}+\StandPartial{y}{x}{11}+\StandPartial{y}{x}{22}-\StandPartial{y}{x}{33}
\end{matrix}
\right.
\end{equation}
\begin{equation}
\label{eq: quaternion over real field, derivative, 1, 1}
\left\{
\begin{matrix}
\displaystyle
\frac{\partial y^{\gi 1}}{\partial x^{\gi 0}}=\StandPartial{y}{x}{01}+\StandPartial{y}{x}{10}+\StandPartial{y}{x}{23}-\StandPartial{y}{x}{32}
\\
\displaystyle
\vphantom{\overset{\rightarrow}{\frac{\partial}{\partial}}}
\frac{\partial y^{\gi 0}}{\partial x^{\gi 1}}=-\StandPartial{y}{x}{01}-\StandPartial{y}{x}{10}+\StandPartial{y}{x}{23}-\StandPartial{y}{x}{32}
\\
\displaystyle
\frac{\partial y^{\gi 3}}{\partial x^{\gi 2}}=-\StandPartial{y}{x}{01}+\StandPartial{y}{x}{10}-\StandPartial{y}{x}{23}-\StandPartial{y}{x}{32}
\\
\displaystyle
\vphantom{\overset{\rightarrow}{\frac{\partial}{\partial}}}
\frac{\partial y^{\gi 2}}{\partial x^{\gi 3}}=\StandPartial{y}{x}{01}-\StandPartial{y}{x}{10}-\StandPartial{y}{x}{23}-\StandPartial{y}{x}{32}
\end{matrix}
\right.
\end{equation}
\begin{equation}
\label{eq: quaternion over real field, derivative, 1, 2}
\left\{
\begin{matrix}
\displaystyle
\frac{\partial y^{\gi 2}}{\partial x^{\gi 0}}=\StandPartial{y}{x}{02}-\StandPartial{y}{x}{13}+\StandPartial{y}{x}{20}+\StandPartial{y}{x}{31}
\\
\displaystyle
\vphantom{\overset{\rightarrow}{\frac{\partial}{\partial}}}
\frac{\partial y^{\gi 3}}{\partial x^{\gi 1}}=\StandPartial{y}{x}{02}-\StandPartial{y}{x}{13}-\StandPartial{y}{x}{20}-\StandPartial{y}{x}{31}
\\
\displaystyle
\vphantom{\overset{\rightarrow}{\frac{\partial}{\partial}}}
\frac{\partial y^{\gi 0}}{\partial x^{\gi 2}}=-\StandPartial{y}{x}{02}-\StandPartial{y}{x}{13}-\StandPartial{y}{x}{20}+\StandPartial{y}{x}{31}
\\
\displaystyle
\vphantom{\overset{\rightarrow}{\frac{\partial}{\partial}}}
\frac{\partial y^{\gi 1}}{\partial x^{\gi 3}}=-\StandPartial{y}{x}{02}-\StandPartial{y}{x}{13}+\StandPartial{y}{x}{20}-\StandPartial{y}{x}{31}
\end{matrix}
\right.
\end{equation}
\begin{equation}
\label{eq: quaternion over real field, derivative, 1, 3}
\left\{
\begin{matrix}
\displaystyle
\frac{\partial y^{\gi 3}}{\partial x^{\gi 0}}=\StandPartial{y}{x}{03}+\StandPartial{y}{x}{12}-\StandPartial{y}{x}{21}+\StandPartial{y}{x}{30}
\\
\displaystyle
\vphantom{\overset{\rightarrow}{\frac{\partial}{\partial}}}
\frac{\partial y^{\gi 2}}{\partial x^{\gi 1}}=-\StandPartial{y}{x}{03}-\StandPartial{y}{x}{12}-\StandPartial{y}{x}{21}+\StandPartial{y}{x}{30}
\\
\displaystyle
\vphantom{\overset{\rightarrow}{\frac{\partial}{\partial}}}
\frac{\partial y^{\gi 1}}{\partial x^{\gi 2}}=\StandPartial{y}{x}{03}-\StandPartial{y}{x}{12}-\StandPartial{y}{x}{21}-\StandPartial{y}{x}{30}
\\
\displaystyle
\vphantom{\overset{\rightarrow}{\frac{\partial}{\partial}}}
\frac{\partial y^{\gi 0}}{\partial x^{\gi 3}}=-\StandPartial{y}{x}{03}+\StandPartial{y}{x}{12}-\StandPartial{y}{x}{21}-\StandPartial{y}{x}{30}
\end{matrix}
\right.
\end{equation}
\begin{equation}
\label{eq: quaternion over real field, derivative, 2, 0}
\left\{
\begin{matrix}
\displaystyle
4\ \StandPartial{y}{x}{00}=\frac{\partial y^{\gi 0}}{\partial x^{\gi 0}}+\frac{\partial y^{\gi 1}}{\partial x^{\gi 1}}+\frac{\partial y^{\gi 2}}{\partial x^{\gi 2}}+\frac{\partial y^{\gi 3}}{\partial x^{\gi 3}}
\\
\displaystyle
\vphantom{\overset{\rightarrow}{\frac{\partial}{\partial}}}
4\ \StandPartial{y}{x}{11}=-\frac{\partial y^{\gi 0}}{\partial x^{\gi 0}}-\frac{\partial y^{\gi 1}}{\partial x^{\gi 1}}+\frac{\partial y^{\gi 2}}{\partial x^{\gi 2}}+\frac{\partial y^{\gi 3}}{\partial x^{\gi 3}}
\\
\displaystyle
\vphantom{\overset{\rightarrow}{\frac{\partial}{\partial}}}
4\ \StandPartial{y}{x}{22}=-\frac{\partial y^{\gi 0}}{\partial x^{\gi 0}}+\frac{\partial y^{\gi 1}}{\partial x^{\gi 1}}-\frac{\partial y^{\gi 2}}{\partial x^{\gi 2}}+\frac{\partial y^{\gi 3}}{\partial x^{\gi 3}}
\\
\displaystyle
\vphantom{\overset{\rightarrow}{\frac{\partial}{\partial}}}
4\ \StandPartial{y}{x}{33}=-\frac{\partial y^{\gi 0}}{\partial x^{\gi 0}}+\frac{\partial y^{\gi 1}}{\partial x^{\gi 1}}+\frac{\partial y^{\gi 2}}{\partial x^{\gi 2}}-\frac{\partial y^{\gi 3}}{\partial x^{\gi 3}}
\end{matrix}
\right.
\end{equation}
\begin{equation}
\label{eq: quaternion over real field, derivative, 2, 1}
\left\{
\begin{matrix}
\displaystyle
4\ \StandPartial{y}{x}{10}=-\frac{\partial y^{\gi 0}}{\partial x^{\gi 1}}+\frac{\partial y^{\gi 1}}{\partial x^{\gi 0}}-\frac{\partial y^{\gi 2}}{\partial x^{\gi 3}}+\frac{\partial y^{\gi 3}}{\partial x^{\gi 2}}
\\
\displaystyle
\vphantom{\overset{\rightarrow}{\frac{\partial}{\partial}}}
4\ \StandPartial{y}{x}{01}=-\frac{\partial y^{\gi 0}}{\partial x^{\gi 1}}+\frac{\partial y^{\gi 1}}{\partial x^{\gi 0}}+\frac{\partial y^{\gi 2}}{\partial x^{\gi 3}}-\frac{\partial y^{\gi 3}}{\partial x^{\gi 2}}
\\
\displaystyle
\vphantom{\overset{\rightarrow}{\frac{\partial}{\partial}}}
4\ \StandPartial{y}{x}{32}=-\frac{\partial y^{\gi 0}}{\partial x^{\gi 1}}-\frac{\partial y^{\gi 1}}{\partial x^{\gi 0}}-\frac{\partial y^{\gi 2}}{\partial x^{\gi 3}}-\frac{\partial y^{\gi 3}}{\partial x^{\gi 2}}
\\
\displaystyle
\vphantom{\overset{\rightarrow}{\frac{\partial}{\partial}}}
4\ \StandPartial{y}{x}{23}=\frac{\partial y^{\gi 0}}{\partial x^{\gi 1}}+\frac{\partial y^{\gi 1}}{\partial x^{\gi 0}}-\frac{\partial y^{\gi 2}}{\partial x^{\gi 3}}-\frac{\partial y^{\gi 3}}{\partial x^{\gi 2}}
\end{matrix}
\right.
\end{equation}
\begin{equation}
\label{eq: quaternion over real field, derivative, 2, 2}
\left\{
\begin{matrix}
\displaystyle
4\ \StandPartial{y}{x}{20}=-\frac{\partial y^{\gi 0}}{\partial x^{\gi 2}}+\frac{\partial y^{\gi 1}}{\partial x^{\gi 3}}+\frac{\partial y^{\gi 2}}{\partial x^{\gi 0}}-\frac{\partial y^{\gi 3}}{\partial x^{\gi 1}}
\\
\displaystyle
\vphantom{\overset{\rightarrow}{\frac{\partial}{\partial}}}
4\ \StandPartial{y}{x}{31}=+\frac{\partial y^{\gi 0}}{\partial x^{\gi 2}}-\frac{\partial y^{\gi 1}}{\partial x^{\gi 3}}+\frac{\partial y^{\gi 2}}{\partial x^{\gi 0}}-\frac{\partial y^{\gi 3}}{\partial x^{\gi 1}}
\\
\displaystyle
\vphantom{\overset{\rightarrow}{\frac{\partial}{\partial}}}
4\ \StandPartial{y}{x}{02}=-\frac{\partial y^{\gi 0}}{\partial x^{\gi 2}}-\frac{\partial y^{\gi 1}}{\partial x^{\gi 3}}+\frac{\partial y^{\gi 2}}{\partial x^{\gi 0}}+\frac{\partial y^{\gi 3}}{\partial x^{\gi 1}}
\\
\displaystyle
\vphantom{\overset{\rightarrow}{\frac{\partial}{\partial}}}
4\ \StandPartial{y}{x}{13}=-\frac{\partial y^{\gi 0}}{\partial x^{\gi 2}}-\frac{\partial y^{\gi 1}}{\partial x^{\gi 3}}-\frac{\partial y^{\gi 2}}{\partial x^{\gi 0}}-\frac{\partial y^{\gi 3}}{\partial x^{\gi 1}}
\end{matrix}
\right.
\end{equation}
\begin{equation}
\label{eq: quaternion over real field, derivative, 2, 3}
\left\{
\begin{matrix}
\displaystyle
4\ \StandPartial{y}{x}{30}=-\frac{\partial y^{\gi 0}}{\partial x^{\gi 3}}-\frac{\partial y^{\gi 1}}{\partial x^{\gi 2}}+\frac{\partial y^{\gi 2}}{\partial x^{\gi 1}}+\frac{\partial y^{\gi 3}}{\partial x^{\gi 0}}
\\
\displaystyle
\vphantom{\overset{\rightarrow}{\frac{\partial}{\partial}}}
4\ \StandPartial{y}{x}{21}=-\frac{\partial y^{\gi 0}}{\partial x^{\gi 3}}-\frac{\partial y^{\gi 1}}{\partial x^{\gi 2}}-\frac{\partial y^{\gi 2}}{\partial x^{\gi 1}}-\frac{\partial y^{\gi 3}}{\partial x^{\gi 0}}
\\
\displaystyle
\vphantom{\overset{\rightarrow}{\frac{\partial}{\partial}}}
4\ \StandPartial{y}{x}{12}=\frac{\partial y^{\gi 0}}{\partial x^{\gi 3}}-\frac{\partial y^{\gi 1}}{\partial x^{\gi 2}}-\frac{\partial y^{\gi 2}}{\partial x^{\gi 1}}+\frac{\partial y^{\gi 3}}{\partial x^{\gi 0}}
\\
\displaystyle
\vphantom{\overset{\rightarrow}{\frac{\partial}{\partial}}}
4\ \StandPartial{y}{x}{03}=-\frac{\partial y^{\gi 0}}{\partial x^{\gi 3}}+\frac{\partial y^{\gi 1}}{\partial x^{\gi 2}}-\frac{\partial y^{\gi 2}}{\partial x^{\gi 1}}+\frac{\partial y^{\gi 3}}{\partial x^{\gi 0}}
\end{matrix}
\right.
\end{equation}
}
{quaternion over real field, derivative}

\DefEq
{
\begin{align*}
&
\begin{pmatrix}
\displaystyle
\frac{\partial y^{\gi 0}}{\partial x^{\gi 0}}
&
\displaystyle
\frac{\partial y^{\gi 1}}{\partial x^{\gi 0}}
&
\displaystyle
\frac{\partial y^{\gi 2}}{\partial x^{\gi 0}}
&
\displaystyle
\frac{\partial y^{\gi 3}}{\partial x^{\gi 0}}
\\
\displaystyle
\vphantom{\overset{\rightarrow}{\frac{\partial}{\partial}}}
\frac{\partial y^{\gi 1}}{\partial x^{\gi 1}}
&
\displaystyle
\frac{\partial y^{\gi 0}}{\partial x^{\gi 1}}
&
\displaystyle
\frac{\partial y^{\gi 3}}{\partial x^{\gi 1}}
&
\displaystyle
\frac{\partial y^{\gi 2}}{\partial x^{\gi 1}}
\\
\displaystyle
\vphantom{\overset{\rightarrow}{\frac{\partial}{\partial}}}
\frac{\partial y^{\gi 2}}{\partial x^{\gi 2}}
&
\displaystyle
\frac{\partial y^{\gi 3}}{\partial x^{\gi 2}}
&
\displaystyle
\frac{\partial y^{\gi 0}}{\partial x^{\gi 2}}
&
\displaystyle
\frac{\partial y^{\gi 1}}{\partial x^{\gi 2}}
\\
\displaystyle
\vphantom{\overset{\rightarrow}{\frac{\partial}{\partial}}}
\frac{\partial y^{\gi 3}}{\partial x^{\gi 3}}
&
\displaystyle
\frac{\partial y^{\gi 2}}{\partial x^{\gi 3}}
&
\displaystyle
\frac{\partial y^{\gi 1}}{\partial x^{\gi 3}}
&
\displaystyle
\frac{\partial y^{\gi 0}}{\partial x^{\gi 3}}
\end{pmatrix}
\\
=&
\begin{pmatrix}
1&-1&-1&-1
\\
1&-1&1&1
\\
1&1&-1&1
\\
1&1&1&-1
\end{pmatrix}
\begin{pmatrix}
\displaystyle
\StandPartial{y}{x}{00}
&
\displaystyle
-\StandPartial{y}{x}{32}
&
\displaystyle
-\StandPartial{y}{x}{13}
&
\displaystyle
-\StandPartial{y}{x}{21}
\\
\displaystyle
\vphantom{\overset{\rightarrow}{\frac{\partial}{\partial}}}
\StandPartial{y}{x}{11}
&
\displaystyle
-\StandPartial{y}{x}{23}
&
\displaystyle
-\StandPartial{y}{x}{02}
&
\displaystyle
-\StandPartial{y}{x}{30}
\\
\displaystyle
\vphantom{\overset{\rightarrow}{\frac{\partial}{\partial}}}
\StandPartial{y}{x}{22}
&
\displaystyle
-\StandPartial{y}{x}{10}
&
\displaystyle
-\StandPartial{y}{x}{31}
&
\displaystyle
-\StandPartial{y}{x}{03}
\\
\displaystyle
\vphantom{\overset{\rightarrow}{\frac{\partial}{\partial}}}
\StandPartial{y}{x}{33}
&
\displaystyle
-\StandPartial{y}{x}{01}
&
\displaystyle
-\StandPartial{y}{x}{20}
&
\displaystyle
-\StandPartial{y}{x}{12}
\end{pmatrix}
\end{align*}
}
{quaternion over real field, derivative, 1}

%% file: Biblio.English.tex
\OpenBiblio


\BiblioItem{texSpaceTime}{Einstein: Electrodynamics of Moving Bodies}
{
Albert Einstein,
On the Electrodynamics of Moving Bodies, 1905,\\
Hendrik Antoon Lorentz, The Principle of Relativity, 37 - 65,
Translated by W Perrett, G B Jeffery
Courier Dover Publications, 1952
}%

\BiblioItem{texSpaceTime}{Einstein: Foundations of general relativity}
{
Albert Einstein,
Die Grundlage der allgemeinen Relativit\"atstheorie,
Ann. Phys., 1916, {\bf 49}, 769 - 822,\\
Einstein's Annalen Papers: The Complete Collection 1901-1922,
edited by J\"urgen Renn, 517 - 571,\\
Wiley-VCH Verlag GmbH \& Co. KGaA, 2005
}%

\BiblioItem{texSpaceTime}{Einstein: Geometry and Experience}
{
Albert Einstein, Geometry and Experience, (1921)\\
Albert Einstein, Sidelights on Relativity, 25 - 56,\\
Courier Dover Publications, 1983
}%

\BiblioItem{texSpaceTime}{Einstein: Main problems of general relativity}
{
Albert Einstein,
Grundgedanken und Probleme der Relativit\"atstheorie, (1923),\\
Nobelstiftelsen, Les Prix Nobel en 1921 - 1922,
Imprimerie Royale, Stockholm, 1923
}%

\BiblioItem{texSpaceTime}{Einstein: Noneuclidean Geometry and Physics}
{
Albert Einstein,
Nichtenklidische Geometrie in der Physik Neue Rundschan, (1925)
Berlin, S. 16 - 20
}%

\BiblioItem{texSpaceTime}{Einstein: Isaak Newton}
{
Albert Einstein,
Isaak Newton, 1927,
Out of My Later Years, 
Citadel Press, 1995, 219 - 223
}%

\BiblioItem{texSpaceTime}{Einstein: On Science}
{
Albert Einstein,
On Science, 
Cosmic Religion, with Other Opinions and Aphorisms,142 - 146,
New York, 1931, 97 - 103
}%

\BiblioItem{texSpaceTime}{Einstein: Autobiographical Notes}
{
Albert Einstein,
Autobiographical Notes, 1949,\\
Paul A. Schilpp, editor, Albert Einstein: Philosopher-Scientist,
Evanston, 
Illinois, The Library of Living Philosophers, 1949, 1 - 95
}%

\BiblioItem{texSpaceTime}{Cite: 104}
{
Cite 104, Source unknown
}%

\BiblioItem{texGenRelativity}{Ghez}
{
Ghez et al.,
The First Measurement of Spectral Lines in a Short-Period Star Bound to the Galaxy's Central Black Hole: A Paradox of Youth,
\href{http://www.journals.uchicago.edu/ApJ/journal/issues/ApJL/v586n2/16990/brief/16990.abstract.html}{ApJL, 586, L127} (2003),
eprint \href{http://arxiv.org/abs/astro-ph/0302299}{arXiv:astro-ph/0302299} (2003)
}%

\BiblioItem{texGenRelativity}{Schodel}
{
R. Sch\"odel et al.,
A star in a 15.2-year orbit around the supermassive black hole at the centre of the Milky Way,
\href{http://www.nature.com/cgi-taf/DynaPage.taf?file=/nature/journal/v419/n6908/abs/nature01121_fs.html}{Nature 419, 694} (2002)
}%

\BiblioItem{texAffine,texGeomObject}{Mielke}
{
Eckehard W. Mielke, Affine generalization of the Komar complex of general relativity,
\href{http://prola.aps.org/searchabstract/PRD/v63/i4/e044018}{Phys. Rev. D 63, 044018} (2001)
}%

\BiblioItem{texAffine}{Obukhov}
{
Yu. N. Obukhov and J. G. Pereira, Metric\hyph affine approach to teleparallel gravity,
\href{http://scitation.aip.org/getabs/servlet/GetabsServlet?prog=normal&id=PRVDAQ000067000004044016000001&idtype=cvips&gifs=Yes}
{Phys. Rev. D 67, 044016} (2003),
eprint \href{http://arxiv.org/abs/gr-qc/0212080}{arXiv:gr-qc/0212080} (2002)
}%

\BiblioItem{texAffine}{Sardanashvily}
{
Giovanni Giachetta, Gennadi Sardanashvily, Dirac Equation in Gauge and Affine-Metric Gravitation Theories,
eprint \href{http://arxiv.org/abs/gr-qc/9511035}{arXiv:gr-qc/9511035} (1995)
}%

\BiblioItem{texAffine}{Gauge}
{
Frank Gronwald and Friedrich W. Hehl, On the Gauge Aspects of Gravity, eprint
\href{http://arxiv.org/abs/gr-qc/9602013}{arXiv:gr-qc/9602013} (1996)
}%

\BiblioItem{texAffine}{Neeman}
{
Yuval Neeman, Friedrich W. Hehl, Test Matter in a Spacetime with Nonmetricity, eprint
\href{http://arxiv.org/abs/gr-qc/9604047}{arXiv:gr-qc/9604047} (1996)
}%

\BiblioItem{texTidal,texAffine,texGeomObject}{torsion}
{
F. W. Hehl, P. von der Heyde, G. D. Kerlick, and J. M. Nester,
General relativity with spin and torsion: Foundations and prospects,\\
\href{http://prola.aps.org/abstract/RMP/v48/i3/p393_1}{Rev. Mod. Phys. 48, 393} (1976)
}%

\BiblioItem{texTidal,texNewton}{Megged}
{
O. Megged, Post-Riemannian Merger of Yang-Mills Interactions with Gravity,
eprint \href{http://arxiv.org/abs/hep-th/0008135}{arXiv:hep-th/0008135} (2001)
}%


\BiblioItem{texNewton}{gr-qc-9604027}
{
Yu.N. Obukhov, E.J. Vlachynsky, W. Esser, R. Tresguerres and F.W. Hehl,
An exact solution of the metric\hyph affine gauge theory with dilation, shear, and spin charges,
eprint \href{http://arxiv.org/abs/gr-qc/9604027}{arXiv:gr-qc/9604027} (1996)
}%

\BiblioItem{texLagrange}{Weinberg}
{
Steven Weinberg. The Quantum Theory of Fields. Cambridge university press.
}%

\BiblioItem{texLagrange}{Reinhardt}
{
Greiner Reinhardt. Field Quantization. Springer.
}%

\BiblioItem{texLagrange}{Landau}
{
L. D. Landau, E. M. Lifshich, The classical theory of fields.
Oxford, New York, Pergamon Press
}%

\BiblioItem{texTidal}{Wheeler}
{
Ignazio Ciufolini, John Wheeler. Gravitation and Inertia.
Princeton university press.
}%

\BiblioItem{texPrefaceTidal,texPrefaceRefernceFrame}{Anderson02}
{
J. D. Anderson, P. A. Laing, E. L. Lau, A. S. Liu, M. M. Nieto, and S. G. Turyshev,
Study of the anomalous acceleration of Pioneer 10 and 11,
\href{http://prola.aps.org/searchabstract/PRD/v65/i8/e082004}{Phys. Rev. D 65, 082004, 50 pp.}, (2002),
eprint \href{http://arxiv.org/abs/gr-qc/0104064}{arXiv:gr-qc/0104064} (2001)
}%

\BiblioItem{texTidal}{Anderson98}
{
J. D. Anderson, P. A. Laing, E. L. Lau, A. S. Liu, M. M. Nieto, and S. G. Turyshev,
Indication, from Pioneer 10/11, Galileo, and Ulysses Data, of an Apparent Anomalous, Weak, Long-Range Acceleration,
\href{http://prola.aps.org/abstract/PRL/v81/i14/p2858_1}{Phys. Rev. Lett. 81, 2858}, (1998),
eprint \href{http://arxiv.org/abs/gr-qc/9808081}{arXiv:gr-qc/9808081} (1998)
}%


\BiblioItem{texDifferential,texReferenceFrame,texFiberedAlgebra,texVocab,texAdditiveMap,texDivisionRing}
{Serge Lang}
{
Serge Lang,
Algebra, Springer, 2002
}%

\BiblioItem{texFiberedAlgebra,texTstarMorphism}{Burris Sankappanavar}
{
S. Burris, H.P. Sankappanavar,
A Course in Universal Algebra, Springer-Verlag (March, 1982),
\\eprint
\href{http://www.math.uwaterloo.ca/~snburris/htdocs/ualg.html}
{http://www.math.uwaterloo.ca/~snburris/htdocs/ualg.html}
\\(The Millennium Edition)
}%

\BiblioItem{texGeomObject}{Shilov}
{
G. E. Shilov,
Calculus, Multivariable Functions,
Moscow, Nauka, 1972
}%

\BiblioItem{texBundle}{Kolmogorov Fomin}
{
A. N. Kolmogorov and S. V. Fomin,
Elements of the Theory of Functions and Functional Analysis,
Courier Dover Publication, 1999
}%

\BiblioItem{}{Lebedev Vorovich}
{
L. P. Lebedev, I. I. Vorovich,
Functional Analysis in Mechanics,
Springer, 2002
}%

\BiblioItem{texAffine,texRepresentation,texTstarRepresentation,texRingAdditiveMap,texPolyvector,texAffineSpace}
{Rashevsky}
{
P. K. Rashevsky, Riemann Geometry and Tensor Calculus,\\
Moscow, Nauka, 1967
}%

\BiblioItem{texBiadditiveMap}
{Kurosh: High Algebra}
{
A. G. Kurosh, High Algebra,
Moscow, Nauka, 1968
}%

\BiblioItem{texPolyvector}{Dubrovin Fomenko Novikov part 1}
{
B. A. Dubrovin, A. T. Fomenko, S. P. Novikov,
Modern Geometry - Methods and Applications,\\
Part 1, The Geometry of Surfaces, Transformation Groups, and Fields,\\
Translated by Robert G. Burns,\\
Springer - New York, 1992
}%

\BiblioItem{texDrcBasis,texBasis}{Korn}
{
Granino A. Korn, Theresa M. Korn,
Mathematical Handbook for Scientists and Engineer,
McGraw-Hill Book Company, New York, San Francisco,
Toronto, London, Sydney, 1968
}%

\BiblioItem{texBundle}{Hocking Young Topology}
{
John G. Hocking, Gail S. Young,
Topology,\\
Courier Dover Publications, 1988
}%


\BiblioItem{texPrefaceRefernceFrame}{Tartaglia}
{
Angelo Tartaglia and Matteo Luca Ruggiero,
Angular Momentum Effects in Michelson\Hyph Morley Type Experiments,
Gen.Rel.Grav. 34, 1371-1382 (2002),\\
eprint \href{http://arxiv.org/abs/gr-qc/0110015}{arXiv:gr-qc/0110015} (2001)
}%

\BiblioItem{texPrefaceRefernceFrame}{Tomozawa}
{
Yukio Tomozawa, Speed of Light in Gravitational Fields, eprint
\href{http://arxiv.org/abs/astro-ph/0303047}{arXiv:astro-ph/0303047} (2004)
}%

\BiblioItem{texPrefaceRefernceFrame}{Magueijo}
{
Joao Magueijo,
Covariant and locally Lorentz-invariant varying speed of light theories,
\href{http://prola.aps.org/abstract/PRD/v62/i10/e103521}{Phys. Rev. D 62, 103521} (2000),
eprint \href{http://arxiv.org/abs/gr-qc/0007036}{arXiv:gr-qc/0007036} (2000)
}%

\BiblioItem{texPrefaceRefernceFrame}{Bassett}
{
Bruce A. Bassett, Stefano Liberati, Carmen Molina-Paris, and Matt Visser,
Geometrodynamics of variable-speed-of-light cosmologies,
\href{http://prola.aps.org/abstract/PRD/v62/i10/e103518}{Phys. Rev. D 62}, 103518 (2000),
eprint \href{http://arxiv.org/abs/astro-ph/0001441}{arXiv:astro-ph/0001441} (2000)
}%

\BiblioItem{texPrefaceRefernceFrame}{Straumann}
{
Lochlainn O'Raifeartaigh and Norbert Straumann,
Gauge theory: Historical origins and some modern developments,
\href{http://prola.aps.org/abstract/RMP/v72/i1/p1_1}{Rev. Mod. Phys. 72, 1} (2000)
}%

\BiblioItem{texPrefaceRefernceFrame}{Lammerzahl}
{
Claus L\"ammerzahl, Mark P. Haugan,
On the interpretation of Michelson\Hyph Morley experiments,
{Phys. Lett. A282 223-229} (2001),\\
eprint \href{http://arxiv.org/abs/gr-qc/0103052}{arXiv:gr-qc/0103052} (2001)
}%

\BiblioItem{texPrefaceRefernceFrame}{Muller}
{
Holger Muller et al.,
Modern Michelson-Morley Experiment using Cryogenic Optical Resonators,
\href{http://prola.aps.org/searchabstract/PRL/v91/i2/e020401}{Phys. Rev. Lett. 91, 020401} (2003),
eprint \href{http://arxiv.org/abs/physics/0305117}{arXiv:physics/0305117} (2000)
}%

\BiblioItem{texPrefaceRefernceFrame,texPrefaceTidal}{Ranada}
{
Antonio F. Ranada,
Pioneer acceleration and variation of light speed: experimental situation,
eprint \href{http://arxiv.org/abs/gr-qc/0402120}{arXiv:gr-qc/0402120} (2004)
}%

\BiblioItem{texBiring,texVectorSpace}{math.QA-0208146}
{
I. Gelfand, S. Gelfand, V. Retakh, R. Wilson,
Quasideterminants,\\
eprint \href{http://arxiv.org/abs/math.QA/0208146}{arXiv:math.QA/0208146} (2002)
}%

\BiblioItem{texBiring,texVectorSpace}{q-alg-9705026}
{
I.Gelfand, V.Retakh,
Quasideterminants, I,\\
eprint \href{http://arxiv.org/abs/q-alg/9705026}{arXiv:q-alg/9705026} (1997)
}%

\BiblioItem{texVectorSpace}{Gelfand Retakh 1991}
{
I. Gelfand and V. Retakh, Determinants of Matrices over Noncommutative Rings, Funct.
Anal. Appl. 25 (1991), no. 2, 91-102
}%

\BiblioItem{texVectorSpace}{Gelfand Retakh 1992}
{
I. Gelfand and V. Retakh, A Theory of Noncommutative Determinants and Characteristic
Functions of Graphs, Funct. Anal. Appl. 26 (1992), no. 4, 1-20
}%

\BiblioItem{texVectorSpace}{hep-th-9407124}
{
I. M. Gelfand, D. Krob, A. Lascoux, B. Leclerc, V.S. Retakh and J.-Y. Thibon,
Noncommutative symmetric functions,\\
eprint \href{http://arxiv.org/abs/hep-th/9407124}{arXiv:hep-th/9407124} (1994)
}%

\BiblioItem{texVectorSpace}{Carl Faith 1}
{
Carl Faith, Algebra: Rings, Modules and Categories I,
Springer - Verlag, Berlin - Heidelberg - New York, 1973
}%



\ifx\PrintBook\Defined
\else
\BiblioItem{texDrcReferenceFrame,texRefernceFrame,texLie,texLieRepresentation}{0412.391}
{
Aleks Kleyn,
Basis Manifold,
eprint \href{http://arxiv.org/abs/math.DG/0412391}{arXiv:math.DG/0412391} (2004)
}%
\fi

\BiblioItem{}{0405.027}
{
Aleks Kleyn,
Reference Frame in General Relativity,\\
eprint \href{http://arxiv.org/abs/gr-qc/0405027}{arXiv:gr-qc/0405027} (2004)
}%

\BiblioItem{}{0405.028}
{
Aleks Kleyn, Metric\hyph Affine Manifold,\\
eprint \href{http://arxiv.org/abs/gr-qc/0405028}{arXiv:gr-qc/0405028} (2004)
}%

\BiblioItem{texFiberedAlgebra,texBundleRelation,texDrcLie}{0701.238}
{
Aleks Kleyn,
Lectures on Linear Algebra over Division Ring,\\
eprint \href{http://arxiv.org/abs/math.GM/0701238}{arXiv:math.GM/0701238} (2007)
}%

\BiblioItem{texBundleRelation,texPrefaceRelation,IntroCalculusBook}{0702.561}
{
Aleks Kleyn,
Algebra Bundle,\\
eprint \href{http://arxiv.org/abs/math.DG/0702561}{arXiv:math.DG/0702561} (2007)
}%


\BiblioItem{texPolymodule}{math.RA-0501237v1}
{
Aleks Kleyn,
Module Over Division Ring, version 1,\\
eprint \href{http://arxiv.org/abs/math/0501237v1}{arXiv:math.RA/0501237v1} (2005)
}%

\ifx\texBiring\Defined
\else
\BiblioItem{texVectorSpace,texFiberedAlgebra}{0612.111}
{
Aleks Kleyn,
Biring of Matrices,\\
eprint \href{http://arxiv.org/abs/math.OA/0612111}{arXiv:math.OA/0612111} (2006)
}%
\fi

\ifx\texBundleRelation\Defined
\else
\BiblioItem{texFiberedMorphism}{0707.2246}
{
Aleks Kleyn,
Fibered Correspondence,\\
eprint \href{http://arxiv.org/abs/0707.2246}{arXiv:0707.2246} (2007)
}%
\fi

\ifx\PrintBook\Defined
\BiblioItem{texPrefaceRelation}{0707.2246}
{
Aleks Kleyn,
Fibered Correspondence,\\
eprint \href{http://arxiv.org/abs/0707.2246}{arXiv:0707.2246} (2007)
}%
\fi

\BiblioItem{texLie}{0803.3276}
{
Aleks Kleyn,
Lorentz Transformation and General Covariance Principle,\\
eprint \href{http://arxiv.org/abs/0803.3276}{arXiv:0803.3276} (2008)
}%

\BiblioItem{EprintCalculus}{0812.4763}
{
Aleks Kleyn,
Introduction into Calculus over Division Ring,\\
eprint \href{http://arxiv.org/abs/0812.4763}{arXiv:0812.4763} (2008)
}%

\BiblioItem{texHomotopy,texTowerRepresentation,texPrefaceIntroGeometry}{9705.009}
{
John C. Baez,
An Introduction to n-Categories,\\
eprint \href{http://arxiv.org/abs/q-alg/9705009}{arXiv:q-alg/9705009} (1997)
}%

\BiblioItem{texPrefaceRelation}{Tolstoi about Anna Karenina}
{
Tolstoi about Anna Karenina,
in book A Karenina Companion, by C. J. G. Turner,
published by Wilfrid Laurier University Press (August 1993)
}%

\BiblioItem{texBundleRelation,texPrefaceRelation,texTstarMorphism,texBundle}
{Cohn: Universal Algebra}
{
Paul M. Cohn,
Universal Algebra,
Springer, 1981
}%

\BiblioItem{texBundle}
{Maunder: Algebraic Topology}
{
C. R. F. Maunder,
Algebraic Topology,
Dover Publications, Inc, Mineola, New York, 1996
}%

\BiblioItem{texFiberedAlgebra}{Pommaret: Partial Differential Equations}
{
J.-F. Pommaret,
Partial Differential Equations and Group Theory,
Springer, 1994
}%

\BiblioItem{texBundleRelation}{Bourbaki: Set Theory}
{
N. Bourbaki,
Theory of sets,
Springer, 2004
}%

\BiblioItem{texBundle,texCartesian,texFiberedAlgebra,texBundleRelation,texFiberedMorphism}
{Bourbaki: General Topology 1}
{
N. Bourbaki,
General Topology, Chapters 1 - 4,
Springer, 1989
}

\BiblioItem{texDifferential,texSecondDifferential}{Bourbaki: General Topology: Chapter 5 - 10}
{
N. Bourbaki,
General Topology, Chapters 5 - 10,
Springer, 1989
}

\BiblioItem{texSecondDifferential}{Bourbaki: Topological Vector Space}
{
N. Bourbaki,
Topological Vector Spaces, Chapters 1 - 5,
Springer, 1987
}

\BiblioItem{texDifferential}{Pontryagin: Topological Group}
{
L. S. Pontryagin,
Selected Works, Volume Two, Topological Groups,
Gordon and Breach Science Publishers, 1986
}

\BiblioItem{texVocab,texDiffEq,texDrcDiffEq,texPrefaceLie}
{Eisenhart: Continuous Groups of Transformations}
{
Eisenhart,
Continuous Groups of Transformations,
Dover Publications, New York, 1961
}

\BiblioItem{texVocab}
{Condon Odabasi}
{
Edward Uhler Condon, Halis Odabasi,
Atomic Structure,
CUP Archive, 1980
}

\BiblioItem{texFiberedMorphism}{Postnikov: Differential Geometry}
{
Postnikov M. M.,
Geometry IV: Differential geometry,
Moscow, Nauka, 1983
}

\BiblioItem{texSecondDifferential,texDifferential}{Fihtengolts: Calculus volume 1}
{
Fihtengolts G. M.,
Differential and Integral Calculus Course, volume 1,
Moscow, Nauka, 1969
}

\BiblioItem{texFiberedAlgebra,texFiberedMorphism}{Hatcher: Algebraic Topology}
{
Allen Hatcher,
Algebraic Topology,
Cambridge University Press, 2002
}

\BiblioItem{texFiberedMorphism}{geometry of differential equations}
{
Vinogradov, A. M., Krasil'shchik, I. S., and Lychagin, V. V.,
Introduction to geometry of nonlinear differential equations,
Nauka, Moscow, 1986
}

\BiblioItem{texFiberedMorphism}{cohomological analysis}
{
A. M. Vinogradov,
Cohomological Analysis of Partial Differential Equations
and Secondary Calculus,
American Mathematical Society, 2001
}

\BiblioItem{texPolyvector}{0801.1734}
{
Brandon S. DiNunno, Richard A. Matzner,
The Volume Inside a Black Hole,\\
eprint \href{http://arxiv.org/abs/0801.1734v1}{arXiv:0801.1734v1} (2008)
}

\BiblioItem{texQuaternion}{Izrail M. Gelfand: Quaternion Groups}
{
I. M. Gelfand, M. I. Graev,
Representation of Quaternion Groups over Localy Compact and
Functional Fields,\\
Funct. Anal. Appl. {\bf 2} (1968) 19 - 33;\\
Izrai.l Moiseevich Gelfand, Semen Grigorevich Gindikin,\\
Izrail M. Gelfand: Collected Papers, volume II, 435 - 449,\\
Springer, 1989
}

\BiblioItem{texAffineSpace}{Bamberg Sternberg}
{
Paul Bamberg, Shlomo Sternberg,
A course in mathematics for students of physics,
Cambridge University Press, 1991
}

\BiblioItem{texPrefaceIntroCalculus}{Vadim Komkov}
{
Vadim Komkov,
Variational Principles of Continuum Mechanics with Engineering Applications: Critical Points Theory,\\
Springer, 1986
}

\CloseBiblio

%% file: Index.English.tex
\OpenIndex
\SetIndexSpace%
\Index{texLinearMap}
   {$1$-\drc form}%
   {1-drc form, vector spaces}%
\SetIndexSpace%
\Index{texBundleRelation}
   {$2$\Hyph ary fibered relation}%
   {2 ary fibered relation}%
\SetIndexSpace%
\Index{texSecondDifferential}
   {$A$\Hyph valued function}%
   {A valued function}%
\Index{texBiring}
   {$(^a_b)$\hyph \CR quasideterminant}%
   {a b cr-quasideterminant}%
\Index{texBiring}
   {$(^a_b)$\hyph \RC quasideterminant}%
   {a b RC-quasideterminant}%
\Index{texDifferential}
   {absolute value on division ring}%
   {absolute value on division ring}%
\Index{texBasis}
   {active representation}%
   {active representation}%
\Index{texDrcBasis}
   {active \sT representation}%
   {active representation, vector space}%
\Index{texBasis}
   {active transformation on basis manifold}%
   {active transformation}%
\Index{texDrcBasis}
   {active transformation on the set of \drc bases}%
   {active transformation, vector space}%
\Index{texAdditiveMap}
   {additive mapping of $D$\Hyph vector spaces}%
   {additive map of D vector spaces}%
\Index{texRingAdditiveMap}
   {additive map of division ring generated by map $G$}%
   {additive map generated by map, division ring}%
\Index{texAdditiveMap}
   {additive map of ring}%
   {Additive map of Ring}%
\Index{texBasis}
   {affine basis}%
   {Affine Basis}%
\Index{texAffineSpace}
   {affine basis}%
   {Affine Basis}%
\Index{texAffineSpace}
   {affine transformation group}%
   {AffineTransformationGroup}%
\Index{texBasis}
   {affine transformation group}%
   {AffineTransformationGroup}%
\Index{texBasis}
   {affine transformation on basis manifold}%
   {affine transformation}%
\Index{texBiring}
   {alternative representation of matrix}%
   {Alternative representation}%
\Index{texRefernceFrame}
   {anholonomic coordinate}%
   {anholonomic coordinate}%
\Index{texRefernceFrame}
   {anholonomic coordinates of connection}%
   {anholonomic coordinates of connection}%
\Index{texRefernceFrame}
   {anholonomic coordinates of vector}%
   {vector anholonomic coordinates}%
\Index{texRefernceFrame}
   {anholonomic coordinates on manifold}%
   {anholonomic coordinates on manifold}%
\Index{texRefernceFrame}
   {anholonomity object}%
   {anholonomity object}%
\Index{texFiberedGroup}
   {antihomomorphism of fibered groups}%
   {antihomomorphism of fibered groups}%
\Index{texBundleRelation}
   {antisymmetric $2$\Hyph ary fibered relation}%
   {antisymmetric 2 ary fibered relation}%
\Index{texFiberedAlgebra}
   {arity of operation}%
   {arity of operation}%
\Index{texTstarRepresentation}
   {associative law for covariant \sT representation}%
   {associative law for covariant starT representation}%
\Index{texTstarRepresentation}
   {associative law for covariant \Ts representation}%
   {associative law for covariant Tstar representation}%
\Index{texVectorField}
   {associative law for \Drc linear maps of vector bundles}%
   {associative law for drc linear maps of vector bundles}%
\Index{texDrcMorphism}
   {associative law for \drc linear maps of vector spaces}%
   {associative law for drc linear maps of vector spaces}%
\Index{texVectorField}
   {associative law for $\mathcal D\star$\Hyph vector fields}%
   {associative law, Dstar vector fields}%
\Index{texVectorSpace}
   {associative law for $D\star$\Hyph vector space}%
   {associative law, Dstar vector space}%
\Index{texLinearMap}
   {associative law for twin representations}%
   {associative law for twin representations}%
\Index{texBundleRelation}
   {associative law of composition of fibered correspondences}%
   {associative law, composition of fibered correspondences}%
\Index{texAffine}
   {auto parallel line}%
   {auto parallel line}%
\Index{texTstarMorphism}
   {automorphism of representation $f$}%
   {automorphism of representation f}%
\SetIndexSpace%
\Index{texPolymodule}
   {($S\star$, $\star T$)\hyph bimodule}%
   {(Sstar,starT)-bimodule}%
\Index{texBundleRelation}
   {base of fibered correspondence}%
   {base of fibered correspondence}%
\Index{texBundle}
   {base of map}%
   {base of map}%
\Index{texBasis}
   {basis manifold of affine space}%
   {Basis Manifold, Affine Space}%
\Index{texAffineSpace}
   {basis manifold of affine space}%
   {Basis Manifold, Affine Space}%
\Index{texAffineSpace}
   {basis manifold of central affine space}%
   {Basis Manifold, Central Affine Space}%
\Index{texBasis}
   {basis manifold of central affine space}%
   {Basis Manifold, Central Affine Space}%
\Index{texDrcBasis}
   {basis manifold of \drc vector space}%
   {basis manifold of vector space}%
\Index{texAffineSpace}
   {basis manifold of Euclid space}%
   {Basis Manifold, Euclid Space}%
\Index{texBasis}
   {basis manifold of Euclid space}%
   {Basis Manifold, Euclid Space}%
\Index{texBasis}
   {basis manifold of vector space}%
   {basis manifold of vector space}%
\Index{texBasis}
   {basis of vector space}%
   {Basis}%
\Index{texLieRepresentation}
   {basis vector of \sT representation}%
   {basis vector of starT representation}%
\Index{texLieRepresentation}
   {basis vector of \Ts representation}%
   {basis vector of Tstar representation}%
\Index{texBiring}
   {biring}%
   {biring}%
\Index{texFiberedMorphism}
   {bundle of level $2$}%
   {bundle of level 2}%
\Index{texFiberedMorphism}
   {bundle of level $n$}%
   {bundle of level n}%
\SetIndexSpace%
\Index{texVectorSpace}
   {\subs rows \dcr vector space}%
   {subs rows dcr vector space}%
\Index{texBiring}
   {\subs row of matrix}%
   {c row}%
\Index{texBiring}
   {$c$\hyph row of matrix}%
   {c-row}%
\Index{texAffine}
   {Cartan connection}%
   {Cartan connection}%
\Index{texAffine}
   {Cartan curvature}%
   {Cartan curvature}%
\Index{texAffine}
   {Cartan derivative}%
   {Cartan derivative}%
\Index{texAffine}
   {Cartan symbol}%
   {Cartan symbol}%
\Index{texAffine}
   {Cartan transport}%
   {Cartan transport}%
\Index{texBundle}
   {Cartesian power $\mathcal{A}$ of bundle $\mathcal{B}$}%
   {Cartesian power A of bundle B}%
\Index{texBundle}
   {Cartesian power $A$ of set $B$}%
   {Cartesian power of set}%
\Index{texCartesian}
   {Cartesian power $n$ of bundle $\mathcal{E}$}%
   {Cartesian power n of bundle E}%
\Index{texCartesian}
   {direct product of bundles}%
   {Cartesian product of bundles}%
\Index{texCartesian}
   {direct product of total spaces}%
   {Cartesian product of total spaces}%
\Index{texHomotopy}
   {category of \drc vector spaces}%
   {category of drc vector spaces}%
\Index{texBundleRelation}
   {category of fibered correspondences over diagonal}%
   {category of fibered correspondences over diagonal}%
\Index{texBundleRelation}
   {category of reduced fibered correspondences}%
   {category of reduced fibered correspondences}%
\Index{texTstarMorphism}
   {category of \Ts representations of $\mathfrak{F}$\Hyph algebra $A$}%
   {category of Tstar representations of F algebra}%
\Index{texTstarMorphism}
   {category of \Ts representations of $\mathfrak{F}$\Hyph algebra from category $\mathcal A$}%
   {category of Tstar representations of F algebra from category}%
\Index{texDifferential}
   {Cauchy sequence in valued division ring}%
   {Cauchy sequence, valued division ring}%
\Index{texDivisionRing}
   {center of ring $D$}%
   {center of ring}%
\Index{texAffineSpace}
   {central affine basis}%
   {Central Affine Basis}%
\Index{texBasis}
   {central affine basis}%
   {Central Affine Basis}%
\Index{texRepresentation}
   {column vector}%
   {column vector}%
\Index{texBundleRelation}
   {commutative diagram of correspondences}%
   {commutative diagram of correspondences}%
\Index{texBundle}
   {compact\hyph open topology}%
   {compact open topology}%
\Index{texDifferential}
   {complete division ring}%
   {complete division ring}%
\Index{texDiffEq}
   {complete system of linear partial differential equations}%
   {Complete System of Linear Partial Differential Equations}%
\Index{texDiffEq}
   {completely integrable system}%
   {completely integrable system}%
\Index{texDifferential}
   {component of the G\^ateaux derivative of map $f(x)$}%
   {component of Gateaux derivative of map, division ring}%
\Index{texSecondDifferential}
   {component of the G\^ateaux derivative of map $\Vector f(\Vector x)$}%
   {component of Gateaux derivative of map, D vector space}%
\Index{texSecondDifferential}
   {component of the G\^ateaux derivative of second order of map  of division ring}%
   {component of Gateaux derivative of Second Order, division ring}%
\Index{texAdditiveMap}
   {component of linear map of $D$\Hyph vector space}%
   {component of linear map, D vector space}%
\Index{texRingAdditiveMap}
   {component of linear map $f$ of division ring}%
   {component of linear map, division ring}%
\Index{texAdditiveMap}
   {component of polyadditive map $\Vector A$}%
   {component of polyadditive map, D vector space}%
\Index{texAdditiveMap}
   {component of polylinear map of division ring}%
   {component of polylinear map, division ring}%
\Index{texSecondDifferential}
   {component of the G\^ateaux derivative of second order of map $\Vector f(\Vector x)$}%
   {component of Gateaux derivative of Second Order, D vector space}%
\Index{texBundleRelation}
   {composition of fibered correspondences}%
   {composition of fibered correspondences}%
\Index{texBundleRelation}
   {composition of reduced fibered correspondences}%
   {composition of reduced fibered correspondences}%
\Index{texBiring}
   {condition of reducibility of products}%
   {condition of reducibility of products}%
\Index{texBundleRelation}
   {continuous correspondence}%
   {continuous correspondence}%
\Index{texDifferential}
   {continuous function of division ring}%
   {continuous function, division ring}%
\Index{texFiberedGroup}
   {contravariant \sT representation of fibered group}%
   {contravariant starT representation of fibered group}%
\Index{texTstarRepresentation}
   {contravariant \sT representation of group}%
   {contravariant starT representation of group}%
\Index{texTstarRepresentation}
   {contravariant \Ts representation of group}%
   {contravariant Tstar representation of group}%
\Index{texFiberedGroup}
   {contravariant \Ts representation of fibered group}%
   {contravariant Tstar representation of fibered group}%
\Index{texVectorSpace}
   {coordinate \drc isomorphism}%
   {coordinate drc isomorphism}%
\Index{texVectorField}
   {coordinate \Drc vector bundle}%
   {coordinate drc vector bundle}%
\Index{texVectorSpace}
   {coordinate \drc vector space}%
   {coordinate drc vector space}%
\Index{texBasis}
   {coordinate isomorphism}%
   {coordinate isomorphism}%
\Index{texVectorSpace}
   {coordinate matrix of set of vectors in \dcr rows vector space}%
   {coordinate matrix of set of vectors, dcr vector space}%
\Index{texVectorSpace}
   {coordinate matrix of set of vectors in \drc rows vector space}%
   {coordinate matrix of set of vectors, drc vector space}%
\Index{texVectorField}
   {coordinate matrix of vector field in \Drc basis}%
   {coordinate matrix of vector field in drc basis}%
\Index{texVectorSpace}
   {coordinate matrix of vector in \drc basis}%
   {coordinate matrix of vector in drc basis}%
\Index{texRefernceFrame}
   {coordinate reference frame}%
   {coordinate reference frame}%
\Index{texDrcBasis}
   {coordinate representation in \drc vector space}%
   {coordinate representation, vector space}%
\Index{texBasis}
   {coordinate representation of group in vector space}%
   {coordinate representation, vector space}%
\Index{texBasis}
   {coordinate vector space}%
   {coordinate vector space}%
\Index{texBasis}
   {coordinates of geometrical object}%
   {coordinates of geometrical object, vector space}%
\Index{texDrcBasis}
   {coordinates of geometrical object in coordinate \drc vector space}%
   {coordinates of geometrical object, coordinate vector space}%
\Index{texBasis}
   {coordinates of geometrical object in coordinate representation}%
   {coordinates of geometrical object, coordinate vector space}%
\Index{texDrcBasis}
   {coordinates of geometrical object in \drc vector space}%
   {coordinates of geometrical object, vector space}%
\Index{texDrcBasis}
   {coordinates of representation}%
   {coordinates of representation}%
\Index{texBasis}
   {coordinates of representation}%
   {coordinates of representation}%
\Index{texVectorSpace}
   {coordinates of set of vectors in \dcr vector space}%
   {coordinates of set of vectors, dcr vector space}%
\Index{texVectorSpace}
   {coordinates of set of vectors in \drc vector space}%
   {coordinates of set of vectors, drc vector space}%
\Index{texVectorField}
   {coordinates of vector field in \Drc basis}%
   {coordinates of vector field in drc basis}%
\Index{texVectorSpace}
   {coordinates of vector in \drc basis}%
   {coordinates of vector in drc basis}%
\Index{texBundleRelation}
   {correspondence continuous on the set}%
   {correspondence continuous on the set}%
\Index{texBundleRelation}
   {correspondence of homomorphism}%
   {correspondence of homomorphism}%
\Index{texFiberedGroup}
   {covariant \sT representation of fibered group}%
   {covariant starT representation of fibered group}%
\Index{texTstarRepresentation}
   {covariant \sT representation of group}%
   {covariant starT representation of group}%
\Index{texTstarRepresentation}
   {covariant \Ts representation of group}%
   {covariant Tstar representation of group}%
\Index{texFiberedGroup}
   {covariant \Ts representation of fibered group}%
   {covariant Tstar representation of fibered group}%
\Index{texBiring}
   {\CR inverse element of biring}%
   {cr-inverse element}%
\Index{texVectorSpace}
   {\CR matrix group}%
   {cr-matrix group}%
\Index{texBiring}
   {\CR power}%
   {cr power}%
\Index{texBiring}
   {\CR product of matrices}%
   {cr-product of matrices}%
\Index{texVectorSpace}
   {\crd vector space}%
   {crd vector space}%
\SetIndexSpace%
\Index{texLinearMap}
   {\rcd basis dual to \drc basis of vector space}%
   {basis dual to basis, drc vector space}%
\Index{texSecondDifferential}
   {$D$\Hyph vector function}%
   {d vector function}%
\Index{texLinearMap}
   {$D$\hyph vector space}%
   {D vector space}%
\Index{texSecondDifferential}
   {$D$\Hyph valued variable}%
   {D valued variable}%
\Index{texVectorSpace}
   {\dcr basis of \subs rows vector space}%
   {dcr basis, c rows vector space}%
\Index{texVectorSpace}
   {\dcr vector}%
   {dcr vector}%
\Index{texVectorSpace}
   {\dcr vector space}%
   {dcr vector space}%
\Index{texBiring}
   {determinant of matrix}%
   {determinant}%
\Index{texTidal}
   {deviation of trajectories}%
   {deviation of trajectories}%
\Index{texBundleRelation}
   {diagonal in bundle}%
   {diagonal in bundle}%
\Index{texBundleRelation}
   {diagram of correspondences}%
   {diagram of correspondences}%
\Index{texVectorSpace}
   {dimension of \drc vector space}%
   {dimension of vector space}%
\Index{texPolymodule}
   {direct product of $D$\Hyph vector spaces}%
   {direct product of D vector spaces}%
\Index{texPolymodule}
   {direct product of division rings}%
   {direct product of division rings}%
\Index{texPolymodule}
   {direct product of \drc vector spaces}%
   {direct product, drc vector space}%
\Index{texFiberedGroup}
   {direct product of representations of fibered group}%
   {direct product of representations of fibered group}%
\Index{texRepresentation}
   {direct product of representations of group}%
   {direct product of representations of group}%
\Index{texTstarRepresentation}
   {direct product of \Ts representations of group}%
   {direct product of representations of group}%
\Index{texLieRepresentation}
   {direct sum of representations}%
   {direct sum of representations}%
\Index{texAdditiveMap}
   {direction in $D$\Hyph vector space $\Vector V$ over field $P$}%
   {direction over field, D vector space}%
\Index{texAdditiveMap}
   {direction in ring $D$ over commutative ring $P$}%
   {direction over commutative ring, ring}%
\Index{texVectorField}
   {distributive law for $\mathcal D\star$\Hyph vector fields}%
   {distributive law, Dstar vector fields}%
\Index{texVectorSpace}
   {distributive law for $D\star$\Hyph vector space}%
   {distributive law, Dstar vector space}%
\Index{texAffineSpace}
   {\drc affine space}%
   {drc affine space}%
\Index{texVectorSpace}
   {\drc automorphism of vector space}%
   {automorphism of vector space}%
\Index{texVectorSpace}
   {\drc basis for \sups rows vector space}%
   {drc basis, r rows vector space}%
\Index{texVectorField}
   {\Drc basis for vector  bundle}%
   {drc basis, vector bundle}%
\Index{texVectorSpace}
   {\drc basis for vector space}%
   {drc basis, vector space}%
\Index{texVectorSpace}
   {\drc isomorphism of vector spaces}%
   {isomorphism of vector spaces}%
\Index{texVectorField}
   {\Drc linear map of vector bundles}%
   {drc linear map of vector bundles}%
\Index{texDrcMorphism}
   {\drc linear map of vector spaces}%
   {drc linear map of vector spaces}%
\Index{texVectorSpace}
   {\drc linear span in vector space}%
   {linear span, vector space}%
\Index{texDrcBasis}
   {\drc linear \sT representation of group}%
   {drc linear starT representation of group}%
\Index{texVectorField}
   {\Drc linearly dependent vector fields}%
   {linearly dependent vector fields}%
\Index{texVectorSpace}
   {\drc linearly dependent vectors}%
   {linearly dependent, vector space}%
\Index{texVectorField}
   {\Drc linearly independent vector fields}%
   {linearly independent vector fields}%
\Index{texVectorSpace}
   {\drc linearly independent vectors}%
   {linearly independent, vector space}%
\Index{texVectorSpace}
   {\drc system of linear equations}%
   {system of linear equations}%
\Index{texVectorSpace}
   {\drc vector}%
   {drc vector}%
\Index{texSecondDifferential}
   {\drc vector function}%
   {drc vector function}%
\Index{texVectorSpace}
   {\drc vector space}%
   {drc vector space}%
\Index{texLinearMap}
   {\Ds component of coordinates of vector $\Vector r$}%
   {Dstar component of coordinates of vector, D vector space}%
\Index{texVectorField}
   {$\mathcal D\star$\hyph product of vector field over scalar}%
   {Dstar product of vector field over scalar, vector space}%
\Index{texVectorField}
   {$\mathcal D\star$\Hyph vector bundle}%
   {Dstar vector bundle}%
\Index{texVectorField}
   {$\mathcal D\star$\Hyph vector field}%
   {Dstar vector field}%
\Index{texVectorSpace}
   {$D\star$\hyph  vector space}%
   {Dstar vector space}%
\Index{texVectorField}
   {$\mathcal D\star$\hyph linear composition of vector fields}%
   {linear composition of vector fields}%
\Index{texVectorSpace}
   {$D\star$\hyph product of vector over scalar}%
   {Dstar product of vector over scalar, vector space}%
\Index{texLinearMap}
   {dual space of \drc vector space}%
   {dual space of drc vector space}%
\Index{texBiring}
   {duality principle for biring}%
   {duality principle for biring}%
\Index{texBiring}
   {duality principle for biring of matrices}%
   {duality principle for biring of matrices}%
\SetIndexSpace%
\Index{texTstarMorphism}
   {effective representation of $\mathfrak{F}$\Hyph algebra $A$}%
   {effective representation of algebra}%
\Index{texVectorSpace}
   {effective representation of division ring}%
   {effective representation of division ring}%
\Index{texVectorField}
   {effective \Ts representation of fibered division ring}%
   {effective representation of fibered division ring}%
\Index{texFiberedAlgebra}
   {effective representation of fibered $\mathfrak{F}$\Hyph algebra}%
   {effective representation of fibered F-algebra}%
\Index{texFiberedGroup}
   {effective \Ts representation of fibered group}%
   {effective representation of fibered group}%
\Index{texTstarRepresentation}
   {effective representation of group}%
   {effective representation of group}%
\Index{texRepresentation}
   {effective representation of group}%
   {effective representation of group}%
\Index{texELie}
   {enhanced Lie group}%
   {enhanced Lie group}%
\Index{texDiffEq}
   {essential parameters}%
   {essential parameters}%
\Index{texVectorSpace}
   {extended matrix of \drc linear equations}%
   {extended matrix, system of drc linear equations}%
\Index{texVectorSpace}
   {extended matrix of \rcd linear equations}%
   {extended matrix, system of rcd linear equations}%
\Index{texBundleRelation}
   {extension of correspondence}%
   {extension of correspondence}%
\Index{texAffine}
   {extreme line}%
   {extreme line}%
\SetIndexSpace%
\Index{texVectorField}
   {fibered coordinate \Drc isomorphism}%
   {fibered coordinate drc isomorphism}%
\Index{texBundleRelation}
   {fibered correspondence from $\mathcal{A}$ to $\mathcal{B}$}%
   {fibered correspondence from A to B}%
\Index{texBundleRelation}
   {fibered correspondence in $\mathcal{A}$}%
   {fibered correspondence in A}%
\Index{texBundleRelation}
   {fibered correspondence of homomorphism}%
   {fibered correspondence of homomorphism}%
\Index{texBundleRelation}
   {fibered equivalence}%
   {fibered equivalence}%
\Index{texFiberedAlgebra}
   {fibered $\mathfrak{F}$\Hyph algebra}%
   {fibered F-algebra}%
\Index{texFiberedAlgebra}
   {fibered $\mathfrak{F}$\Hyph subalgebra}%
   {fibered F-subalgebra}%
\Index{texFiberedAlgebra}
   {fibered group}%
   {fibered group}%
\Index{texFiberedMorphism}
   {fibered identification morphism}%
   {fibered identification morphism}%
\Index{texFiberedMorphism}
   {fibered little group}%
   {fibered little group}%
\Index{texBundle}
   {fibered morphism from bundle $\mathcal{A}$ into $\mathcal{B}$}%
   {fibered morphism from A into B}%
\Index{texFiberedMorphism}
   {fibered natural morphism}%
   {fibered natural morphism}%
\Index{texBundleRelation}
   {fibered ordering}%
   {fibered ordering}%
\Index{texBundleRelation}
   {fibered preordering}%
   {fibered preordering}%
\Index{texFiberedAlgebra}
   {fibered ring}%
   {fibered ring}%
\Index{texFiberedMorphism}
   {fibered stability group}%
   {fibered stability group}%
\Index{texBundle}
   {fibered subset}%
   {fibered subset}%
\Index{texNewton}
   {field-strength tensor}%
   {field-strength tensor}%
\Index{texBundleRelation}
   {filter $\mathfrak{F}$ converges to $A$}%
   {filter converges}%
\Index{texNewton}
   {first Newton law}%
   {First Newton law}%
\Index{texFrechet}
   {the Fr\'echet \Ds derivative of map $f$ of division ring $D$ at point $x$}%
   {Frechet Dstar derivative of map, division ring}%
\Index{texFiberedMorphism}
   {free \Ts representation of fibered group}%
   {free representation of fibered group}%
\Index{texTstarRepresentation}
   {free \Ts representation of group}%
   {free representation of group}%
\Index{texAffine}
   {Frenet transport}%
   {Frenet transport}%
\Index{texSecondDifferential}
   {function continuous with respect to set of arguments}%
   {function continuous with respect to set of arguments}%
\Index{texDifferential}
   {function is projective over field $R$ and continuous in direction over field $R$}%
   {projective function is continuous in direction over field R, division ring}%
\Index{texSecondDifferential}
   {function of $D$\Hyph vector space $\Vector{V}$ to $D$|Hyph vector space $\Vector W$ differentiable in the G\^ateaux sense}%
   {function differentiable in Gateaux sense, D vector space}%
\Index{texDifferential}
   {function of division ring differentiable in the G\^ateaux sense}%
   {function differentiable in Gateaux sense, division ring}%
\Index{texFrechet}
   {function of division ring \Ds differentiable in the Fr\'echet sense}%
   {function Dstar differentiable in Frechet sense, division ring}%
\Index{texSecondDifferential}
   {function of $\gi n$ $D$\Hyph valued variables}%
   {function of n D valued variables}%
\Index{texDifferential}
   {fundamental sequence in valued division ring}%
   {fundamental sequence, valued division ring}%
\SetIndexSpace%
\Index{texRefernceFrame}
   {$G$\Hyph reference frame}%
   {G reference frame}%
\Index{texBasis}
   {\Gbasis\ of vector space}%
   {G-basis}%
\Index{texBasis}
   {\Gcoords\ of basis}%
   {G-coordinates}%
\Index{texSecondDifferential}
   {the G\^ateaux \crd derivative of map $\Vector f$ of $D$\hyph vector space $\Vector V$ to $D$\hyph vector space $\Vector W$}%
   {Gateaux crd derivative of map, D vector space}%
\Index{texDifferential}
   {the G\^ateaux derivative of map $f$}%
   {Gateaux derivative of map, division ring}%
\Index{texSecondDifferential}
   {the G\^ateaux derivative of order $n$ of map $f$ of division ring}%
   {Gateaux derivative of Order n, division ring}%
\Index{texSecondDifferential}
   {the G\^ateaux derivative of order $n$ of map $\Vector f$}%
   {Gateaux derivative of Order n, D vector space}%
\Index{texSecondDifferential}
   {the G\^ateaux derivative of order $n$ of map $\Vector f$}%
   {Gateaux derivative of Order n, D vector space}%
\Index{texDifferential}
   {the G\^ateaux derivative of second order of map of division ring}%
   {Gateaux derivative of Second Order, division ring}%
\Index{texDifferential}
   {the G\^ateaux differential of map $f$}%
   {Gateaux differential of map, division ring}%
\Index{texSecondDifferential}
   {the G\^ateaux differential of second order of map of division ring}%
   {Gateaux differential of Second Order, division ring}%
\Index{texSecondDifferential}
   {the G\^ateaux differential of second order of mapping $\Vector f$}%
   {Gateaux differential of Second Order, D vector space}%
\Index{texSecondDifferential}
   {the G\^ateaux \drc derivative of map $\Vector f$ of $D$\Hyph vector space $\Vector V$ to $D$\Hyph vector space $\Vector W$}%
   {Gateaux drc derivative of map, D vector space}%
\Index{texDifferential}
   {the G\^ateaux \Ds derivative of map $f$ of division ring $D$}%
   {Gateaux Dstar derivative of map, division ring}%
\Index{texSecondDifferential}
   {the G\^ateaux Jacobian of map of $D$\Hyph vector space}%
   {Gateaux Jacobian of map, D vector space}%
\Index{texSecondDifferential}
   {the G\^ateaux partial \drc derivative of map $f^b$ with respect to variable $x^a$}%
   {Gateaux partial drc derivative of map with respect to variable, D vector space}%
\Index{texSecondDifferential}
   {the G\^ateaux partial \drc derivative of map $f^b$ with respect to variable $x^a$}%
   {Gateaux partial crd derivative of map with respect to variable, D vector space}%
\Index{texDifferential}
   {the G\^ateaux \sD derivative of map $f$ of division ring $D$}%
   {Gateaux starD derivative of map, division ring}%
\Index{texRingAdditiveMap}
   {generator of additive map}%
   {generator of additive map, division ring}%
\Index{texDrcBasis}
   {geometrical object defined in \drc vector space}%
   {geometrical object, vector space}%
\Index{texBasis}
   {geometrical object in coordinate representation}%
   {geometrical object, coordinate vector space}%
\Index{texDrcBasis}
   {geometrical object in coordinate representation defined in \drc vector space}%
   {geometrical object, coordinate vector space}%
\Index{texBasis}
   {geometrical object in vector space}%
   {geometrical object, vector space}%
\Index{texDrcBasis}
   {geometrical object of type $A$}%
   {geometrical object of type A, vector space}%
\Index{texBasis}
   {geometrical object of type $A$ in vector space}%
   {geometrical object of type A, vector space}%
\Index{texGroupRing}
   {group algebra}%
   {group algebra}%
\Index{texBasis}
   {\Gspace}%
   {GSpace}%
\Index{texSecondDifferential}
   {the G\^ateaux derivative of map $\Vector f$ of normed $D$\Hyph vector space $\Vector{V}$ to normed $D$\Hyph vector space $\Vector{W}$}%
   {Gateaux derivative of map, D vector space}%
\Index{texSecondDifferential}
   {the G\^ateaux differential of map $\Vector f$ of normed $D$\Hyph vector space $\Vector{V}$ to normed $D$\Hyph vector space $\Vector{W}$}%
   {Gateaux differential of map, D vector space}%
\Index{texSecondDifferential}
   {the G\^ateaux mixed partial derivative of map $f^j$ with respect to variables $v^i$, $v^j$}%
   {Gateaux partial derivative of Second Order, D vector space}%
\Index{texSecondDifferential}
   {the G\^ateaux partial derivative of map $f^j$ with respect to variable $v^i$}%
   {Gateaux partial derivative, D vector space}%
\SetIndexSpace%
\Index{texBiring}
   {Hadamard inverse of matrix}%
   {Hadamard inverse of matrix}%
\Index{texRefernceFrame}
   {holonomic coordinates of connection}%
   {holonomic coordinates of connection}%
\Index{texRefernceFrame}
   {holonomic coordinates of vector}%
   {vector holonomic coordinates}%
\Index{texFiberedGroup}
   {homogeneous bundle of fibered group}%
   {homogeneous bundle of fibered group}%
\Index{texSecondDifferential}
   {homogeneous map of degree $k$ over field $F$}%
   {homogeneous map of degree over field, D vector space}%
\Index{texTstarRepresentation}
   {homogeneous space of group}%
   {homogeneous space of group}%
\Index{texRepresentation}
   {homogeneous space of group}%
   {homogeneous space of group}%
\Index{texFiberedAlgebra}
   {homomorphism of fibered $\mathfrak{F}$\Hyph algebras}%
   {homomorphism of fibered F-algebras}%
\Index{texFiberedGroup}
   {homomorphism of fibered groups}%
   {homomorphism of fibered groups}%
\SetIndexSpace%
\Index{texLieRepresentation}
   {infinitesimal generator}%
   {infinitesimal generator}%
\Index{texLinearLie}
   {infinitesimal generators of group Lie}%
   {infinitesimal generators of group Lie}%
\Index{texDrcBasis}
   {invariance principle}%
   {invariance principle}%
\Index{texBasis}
   {invariance principle in vector space}%
   {invariance principle, vector space}%
\Index{texBundleRelation}
   {inverse fibered correspondence}%
   {inverse fibered correspondence}%
\Index{texBundleRelation}
   {inverse reduced fibered correspondence}%
   {inverse reduced fibered correspondence}%
\Index{texFiberedAlgebra}
   {isomorphism of fibered $\mathfrak{F}$\Hyph algebras}%
   {isomorphism of fibered F-algebras}%
\SetIndexSpace%
\Index{texAdditiveMap}
   {kernel of additive map of division ring}%
   {kernel of additive map, division ring}%
\Index{texAdditiveMap}
   {kernel of additive map of $D$\Hyph vector space}%
   {kernel of additive map, D vector space}%
\Index{texFiberedGroup}
   {kernel of inefficiency of representation of fibered group}%
   {kernel of inefficiency of representation of fibered group}%
\Index{texRepresentation}
   {kernel of inefficiency of representation of group}%
   {kernel of inefficiency of representation of group}%
\Index{texTstarRepresentation}
   {kernel of inefficiency of \Ts representation of group $G$}%
   {kernel of inefficiency of representation of group}%
\Index{texDiffProperty}
   {Killing equation}%
   {Killing equation}%
\Index{texDiffProperty}
   {Killing equation of second type}%
   {Killing equation second type}%
\Index{texDiffProperty}
   {Killing vector of second type}%
   {Killing vector second type}%
\Index{texBiring}
   {Kronecker symbol}%
   {Kronecker symbol}%
\SetIndexSpace%
\Index{texLie}
   {left defined Lie algebra of Lie group}%
   {left defined Lie algebra}%
\Index{texLie}
   {left invariant vector field}%
   {left invariant vector}%
\Index{texVectorSpace}
   {left module}%
   {left module}%
\Index{texTstarRepresentation}
   {left shift on group}%
   {left shift}%
\Index{texFiberedGroup}
   {left shift on fibered group}%
   {Tstar shift, fibered group}%
\Index{texRepresentation}
   {left shift on group}%
   {left shift, group}%
\Index{texLie}
   {left structural constant of Lie algebra}%
   {left structural constant of Lie algebra}%
\Index{texVectorSpace}
   {left vector space}%
   {left vector space}%
\Index{texRepresentation}
   {left-side contravariant representation of group}%
   {left-side contravariant representation of group}%
\Index{texRepresentation}
   {left-side covariant representation of group}%
   {left-side covariant representation of group}%
\Index{texTstarMorphism}
   {left-side representation of $\mathfrak{F}$\Hyph algebra $A$ in set $M$}%
   {left-side representation of algebra}%
\Index{texFiberedAlgebra}
   {left-side representation of fibered $\mathfrak{F}$\Hyph algebra}%
   {left-side representation of fibered F-algebra}%
\Index{texTstarMorphism}
   {left-side transformation}%
   {left-side transformation}%
\Index{texFiberedAlgebra}
   {left-side transformation on bundle}%
   {left-side transformation of bundle}%
\Index{texLie}
   {Lie algebra of Lie group}%
   {algebra Lie group Lie}%
\Index{texDiffProperty}
   {Lie derivative}%
   {Lie derivative}%
\Index{texAffine}
   {Lie derivative of connection}%
   {Lie derivative of connection}%
\Index{texAffine}
   {Lie derivative of metric}%
   {Lie derivative of metric}%
\Index{texLie}
   {Lie group basic operators}%
   {Lie group basic operators}%
\Index{texBundleRelation}
   {lift of correspondence}%
   {lift of correspondence}%
\Index{texBundle}
   {lift of map}%
   {lift of map}%
\Index{texBundleRelation}
   {limit of correspondence with respect to the filter}%
   {limit of correspondence with respect to the filter}%
\Index{texBundleRelation}
   {limit of filter}%
   {limit of filter}%
\Index{texDifferential}
   {limit of sequence in valued division ring}%
   {limit of sequence, valued division ring}%
\Index{texBundleRelation}
   {limit set of filter}%
   {limit set of filter}%
\Index{texAdditiveMap}
   {linear map of $D$\Hyph vector spaces}%
   {linear map, vector space}%
\Index{texAdditiveMap}
   {linear map of $D$\Hyph vector spaces over field $F$}%
   {linear map over field, vector space}%
\Index{texRingAdditiveMap}
   { linear map of division ring}%
   {linear map, division ring}%
\Index{texRepresentation}
   {linear representation of group}%
   {linear representation of group}%
\Index{texTstarRepresentation}
   {little group}%
   {little group}%
\Index{texRefernceFrame}
   {local reference frame}%
   {local reference frame}%
\Index{texBundle}
   {locally compact at point $p$ space}%
   {locally compact at point space}%
\Index{texBundle}
   {locally compact space}%
   {locally compact space}%
\Index{texRefernceFrame}
   {Lorentz transformation}%
   {Lorentz transformation}%
\SetIndexSpace%
\Index{texPolyvector}
   {$m$\Hyph dimensional parallelepiped}%
   {m dimensional parallelepiped}%
\Index{texPolyvector}
   {$m$\Hyph vector}%
   {m-vector}%
\Index{texAdditiveMap}
   {map of $D$\Hyph vector space, multiplicative over field}%
   {map multiplicative over field, D vector space}%
\Index{texAdditiveMap}
   {map of $D$\Hyph vector space, projective over field}%
   {map projective over field, D vector space}%
\Index{texAdditiveMap}
   {map of division ring, projective over commutative ring}%
   {map projective over commutative ring, ring}%
\Index{texAdditiveMap}
   {map of ring $D$, linear over commutative ring $F$}%
   {linear map over commutative ring, ring}%
\Index{texAdditiveMap}
   {map of ring, multiplicative over commutative ring}%
   {map multiplicative over commutative ring, ring}%
\Index{texAdditiveMap}
   {map of rings $R_1$, ..., $R_n$ polylinear over commutative ring $P$}%
   {map polylinear over commutative ring, ring}%
\Index{texDrcReferenceFrame}
   {map of type $G$ on manifold}%
   {map of type G on manifold}%
\Index{texAdditiveMap}
   {map polylinear over commutative ring}%
   {map polylinear over commutative ring, ring}%
\Index{texCartesian}
   {mapping space}%
   {mapping space}%
\Index{texDrcMorphism}
   {matrix of \drc linear map}%
   {matrix of drc linear map}%
\Index{texVectorField}
   {matrix of fibered \Drc linear map}%
   {matrix of fibered drc linear map}%
\Index{texGeomObject}
   {metric-affine manifold}%
   {metric-affine manifold}%
\Index{texTowerRepresentation}
   {morphism from tower of \Ts representations into tower of \Ts representations}%
   {morphism from tower of representations into tower of representations}%
\Index{texFiberedMorphism}
   {morphism of fibered \Ts representations from $\mathcal{F}$ into $\mathcal{G}$}%
   {morphism of fibered representations from f into g}%
\Index{texTstarMorphism}
   {morphism of representations from $f$ into $g$}%
   {morphism of representations from f into g}%
\Index{texTstarMorphism}
   {morphism of representations of $\mathfrak{F}$ algebra}%
   {morphism of representations of F algebra}%
\Index{texTstarMorphism}
   {morphism of representations of $\mathfrak{F}$\Hyph algebra in $\mathfrak{H}$\Hyph algebra}%
   {morphism of representations of F algebra in H algebra}%
\Index{texFiberedMorphism}
   {morphism of \Ts representations of fibered $\mathfrak{F}$\Hyph algebra}%
   {morphism of representations of fibered F algebra}%
\Index{texAffineSpace}
   {movement on basis manifold}%
   {movement transformation}%
\SetIndexSpace%
\Index{texBundleRelation}
   {$n$\Hyph ary fibered relation}%
   {fibered relation}%
\Index{texGeomObject}
   {nonmetricity}%
   {nonmetricity}%
\Index{texVectorSpace}
   {nonsingular \drc system of linear equations}%
   {nonsingular system of linear equations}%
\Index{texRepresentation}
   {nonsingular transformation}%
   {nonsingular transformation}%
\Index{texDifferential}
   {norm of mapping of division ring}%
   {norm of map, division ring}%
\Index{texSecondDifferential}
   {norm of map $\Vector A$ of normed $D$\hyph vector space}%
   {norm of map, D vector space}%
\Index{texAdditiveMap}
   {norm of quaternion}%
   {norm of quaternion}%
\Index{texSecondDifferential}
   {norm on $D$\Hyph vector space}%
   {norm on D vector space}%
\Index{texSecondDifferential}
   {normed $D$\Hyph vector space}%
   {normed D vector space}%
\SetIndexSpace%
\Index{texFiberedAlgebra}
   {operation on bundle}%
   {operation on bundle}%
\Index{texBundleRelation}
   {opposite fibered preordering}%
   {opposite fibered preordering}%
\Index{texFiberedGroup}
   {orbit of representation of fibered group}%
   {orbit of representation of fibered group}%
\Index{texRepresentation}
   {orbit of representation of group}%
   {orbit of representation of group}%
\Index{texTstarRepresentation}
   {orbit of \Ts representation of group}%
   {orbit of representation of group}%
\Index{texBasis}
   {orthonornal basis}%
   {Orthonornal Basis}%
\Index{texAffineSpace}
   {orthonornal basis}%
   {Orthonornal Basis}%
\SetIndexSpace%
\Index{texGeomObject}
   {parallelogram}%
   {parallelogram}%
\Index{texAdditiveMap}
   {partial additive map of variable $v^i$}%
   {partial additive map of variable}%
\Index{texBasis}
   {passive representation}%
   {passive representation}%
\Index{texBasis}
   {passive transformation on basis manifold}%
   {passive transformation}%
\Index{texDrcBasis}
   {passive transformation on the set of \drc bases}%
   {passive transformation, vector space}%
\Index{texDrcBasis}
   {passive \Ts representation}%
   {passive representation}%
\Index{texRefernceFrame}
   {pfaffian derivative}%
   {pfaffian derivative}%
\Index{texAdditiveMap}
   {polyadditive mapping of $(n)$\hyph $D$\hyph vector spaces}%
   {polyadditive map of D vector spaces}%
\Index{texRingAdditiveMap}
   {polyadditive map of rings}%
   {polyadditive map of rings}%
\Index{texRingAdditiveMap}
   {polylinear map of rings}%
   {polylinear map of rings}%
\Index{texRingAdditiveMap}
   {polylinear skew symmetric map}%
   {polylinear map skew symmetric, division ring}%
\Index{texRingAdditiveMap}
   {polylinear symmetric map}%
   {polylinear map symmetric, division ring}%
\Index{texPolyvector}
   {polyvector}%
   {polyvector}%
\Index{texEuclideanSpace}
   {positive definite biadditive map}%
   {positive definite biadditive map}%
\Index{texNewton}
   {potential energy}%
   {potential energy}%
\Index{texDrcBasis}
   {product of geometrical object and constant}%
   {product of geometrical object and constant}%
\Index{texBasis}
   {product of geometrical object and constant in vector space}%
   {product of geometrical object and constant, vector space}%
\Index{texPolymodule}
   {product of groups}%
   {product of groups}%
\Index{texTstarMorphism}
   {product of morphisms of representations of $\mathfrak{F}$\Hyph algebra}%
   {product of morphisms of representations of F algebra}%
\Index{texVectorField}
   {product of morphisms of \Ts representations of fibered $\mathfrak{F}$\Hyph algebra}%
   {product of morphisms of representations of fibered F algebra}%
\Index{texAdditiveMap}
   {product of objects in category}%
   {product of objects in category}%
\Index{texBundle}
   {projection of bundle $\mathcal{E}$ along fiber $E$}%
   {projection of bundle along fiber}%
\SetIndexSpace%
\Index{texAffineSpace}
   {quasi affine transformation on basis manifold}%
   {quasi affine transformation}%
\Index{texBasis}
   {quasi affine transformation on basis manifold}%
   {quasi affine transformation}%
\Index{texAffineSpace}
   {quasi movement on basis manifold}%
   {quasi movement}%
\Index{texBasis}
   {quasi movement on basis manifold}%
   {quasi movement}%
\Index{texAdditiveMap}
   {quaternion algebra $E$ over the field $F$}%
   {quaternion algebra over the field}%
\Index{texFiberedMorphism}
   {quotient bundle}%
   {quotient bundle}%
\SetIndexSpace%
\Index{texBiring}
   {\sups row of matrix}%
   {r row}%
\Index{texVectorSpace}
   {\sups rows \drc vector space}%
   {sups rows drc vector space}%
\Index{texBiring}
   {$r$\hyph row of matrix}%
   {r-row}%
\Index{texBiring}
   {\RC inverse element of biring}%
   {rc-inverse element}%
\Index{texVectorSpace}
   {\RC major minor}%
   {RC-major minor}%
\Index{texVectorSpace}
   {\RC matrix group}%
   {rc-matrix group}%
\Index{texVectorSpace}
   {\RC nonsingular matrix}%
   {RC nonsingular matrix}%
\Index{texBiring}
   {\RC power}%
   {rc power}%
\Index{texBiring}
   {\RC product of matrices}%
   {rc-product of matrices}%
\Index{texBiring}
   {\RC quasideterminant}%
   {RC-quasideterminant}%
\Index{texVectorSpace}
   {\RC rank of matrix}%
   {rc-rank of matrix}%
\Index{texVectorSpace}
   {\RC singular matrix}%
   {RC singular matrix}%
\Index{texDrcBasis}
   {\rcd representation of group}%
   {rcd linear representation of group}%
\Index{texVectorSpace}
   {\rcd vector space}%
   {rcd vector space}%
\Index{texCartesian}
   {reduced Cartesian product of bundles}%
   {reduced Cartesian product of bundles}%
\Index{texCartesian}
   {reduced Cartesian product of total spaces}%
   {reduced Cartesian product of total spaces}%
\Index{texBundleRelation}
   {reduced fibered correspondence from $\mathcal{A}$ to $\mathcal{B}$}%
   {reduced fibered correspondence from A to B}%
\Index{texBundleRelation}
   {reduced fibered correspondence in $\mathcal{A}$}%
   {reduced fibered correspondence in A}%
\Index{texBiring}
   {reducible biring}%
   {reducible biring}%
\Index{texRefernceFrame}
   {reference frame in event space}%
   {reference frame in event space}%
\Index{texDrcReferenceFrame}
   {reference frame manifold}%
   {reference frame manifold}%
\Index{texBundleRelation}
   {reflexive $2$\Hyph ary fibered relation}%
   {reflexive 2 ary fibered relation}%
\Index{texTowerRepresentation}
   {representation of $\mathfrak{F}$\Hyph algebra $A$ in category $\mathcal B$}%
   {representation of F algebra in category}%
\Index{texTstarMorphism}
   {representation of $\mathfrak{F}$\Hyph algebra $A$ in set $M$}%
   {representation of algebra}%
\Index{texRepresentation}
   {representation of group}%
   {representation of group}%
\Index{texDrcBasis}
   {representative of geometrical object in vector space}%
   {representative of geometrical object, vector space}%
\Index{texBasis}
   {representative of geometrical object in vector space}%
   {representative of geometrical object, vector space}%
\Index{texBundleRelation}
   {restriction of correspondence $\Phi$ to set $C$}%
   {restriction of correspondence}%
\Index{texLie}
   {right defined Lie algebra of Lie group}%
   {right defined Lie algebra}%
\Index{texLie}
   {right invariant vector field}%
   {right invariant vector}%
\Index{texVectorSpace}
   {right module}%
   {right module}%
\Index{texTstarRepresentation}
   {right shift on group}%
   {right shift}%
\Index{texRepresentation}
   {right shift on group}%
   {right shift, group}%
\Index{texLie}
   {right structural constant of Lie algebra}%
   {right structural constant of Lie algebra}%
\Index{texVectorSpace}
   {right vector space}%
   {right vector space}%
\Index{texRepresentation}
   {right-side contravariant representation of group}%
   {right-side contravariant representation of group}%
\Index{texRepresentation}
   {right-side covariant representation of group}%
   {right-side covariant representation of group}%
\Index{texTstarMorphism}
   {right-side representation of $\mathfrak{F}$\Hyph algebra $A$ in set $M$}%
   {right-side representation of algebra}%
\Index{texFiberedAlgebra}
   {right-side representation of fibered $\mathfrak{F}$\Hyph algebra}%
   {right-side representation of fibered F-algebra}%
\Index{texTstarMorphism}
   {right-side transformation}%
   {right-side transformation}%
\Index{texRepresentation}
   {right-side transformation}%
   {right-side transformation}%
\Index{texDivisionRing}
   {ring has characteristic $0$}%
   {ring has characteristic 0}%
\Index{texDivisionRing}
   {ring has characteristic $p$}%
   {ring has characteristic p}%
\Index{texRepresentation}
   {row vector}%
   {row vector}%
\Index{texVectorSpace}
   {$R\star$\Hyph module}%
   {Rstar-module}%
\SetIndexSpace%
\Index{texNewton}
   {scalar potential}%
   {scalar potential}%
\Index{texNewton}
   {second Newton law}%
   {Second Newton law}%
\Index{texPolyvector}
   {simple polyvector}%
   {simple polyvector}%
\Index{texTstarMorphism}
   {single transitive representation of $\mathfrak{F}$\Hyph algebra $A$}%
   {single transitive representation of algebra}%
\Index{texFiberedAlgebra}
   {single transitive representation of fibered $\mathfrak{F}$\Hyph algebra}%
   {single transitive representation of fibered F-algebra}%
\Index{texRepresentation}
   {single transitive representation of group}%
   {single transitive representation of group}%
\Index{texAdditiveMap}
   {singular linear mapp of $D$\Hyph vector space}%
   {singular additive map, D vector space}%
\Index{texAdditiveMap}
   {singular linear map of division ring}%
   {singular additive map, division ring}%
\Index{texPolyvector}
   {skew product of vectors}%
   {skew product of vectors}%
\Index{texTstarRepresentation}
   {space of orbits of \Ts representation}%
   {space of orbits of Ts representation}%
\Index{texTidal}
   {speed of deviation}%
   {speed of deviation}%
\Index{texVectorField}
   {$(\mathcal S\RCstar,\mathcal T\RCstar)$\Hyph linear map of vector bundles}%
   {src trc linear map of vector bundles}%
\Index{texDrcMorphism}
   {$(S\RCstar,T\RCstar)$\Hyph linear map of vector spaces}%
   {src trc linear map of vector spaces}%
\Index{texTstarRepresentation}
   {stability group}%
   {stability group}%
\Index{texAdditiveMap}
   {standard component of additive map $f$ over field $F$}%
   {standard component of additive map, division ring}%
\Index{texBiadditiveMap}
   {standard component of biadditive map $f$ over field $F$}%
   {standard component of biadditive map, division ring}%
\Index{texDifferential}
   {standard component of the G\^ateaux differential of map $f$}%
   {standard component of Gateaux differential, division ring}%
\Index{texRingAdditiveMap}
   {standard component of linear map of division ring}%
   {standard component of linear map, division ring}%
\Index{texRingAdditiveMap}
   {standard component of polylinear map $f$ of division ring}%
   {standard component of polylinear map, division ring}%
\Index{texBiadditiveMap}
   {standard component of quadratic map $f$ over field $F$}%
   {standard component of quadratic map, division ring}%
\Index{texTensorProduct}
   {standard component of tensor}%
   {standard component of tensor, division ring}%
\Index{texDrcBasis}
   {standard coordinates of basis}%
   {standard coordinates of basis}%
\Index{texBasis}
   {standard coordinates of basis}%
   {standard coordinates of basis}%
\Index{texAdditiveMap}
   {standard representation of additive map of division ring over field $F$}%
   {additive map, standard representation, division ring}%
\Index{texBiadditiveMap}
   {standard representation of biadditive map of division ring over field $F$}%
   {biadditive map, standard representation, division ring}%
\Index{texBiadditiveMap}
   {standard representation of biadditive map of division ring over field}%
   {biadditive map, standard representation, division ring}%
\Index{texDifferential}
   {standard representation of the G\^ateaux differential of map of division ring over field $F$}%
   {Gateaux differential, standard representation, division ring}%
\Index{texRingAdditiveMap}
   {standard representation of linear map of division ring}%
   {linear map, standard representation, division ring}%
\Index{texBiring}
   {standard representation of matrix}%
   {Standard representation}%
\Index{texRingAdditiveMap}
   {standard representation of polylinear map of division ring}%
   {polylinear map, standard representation, division ring}%
\Index{texBiadditiveMap}
   {standard representation of quadratic map of division ring over field $F$}%
   {quadratic map, standard representation, division ring}%
\Index{texLinearMap}
   {\sD component of coordinates of vector $\Vector r$}%
   {starD component of coordinates of vector, D vector space}%
\Index{texVectorSpace}
   {$\star D$\hyph vector space}%
   {starD-vector space}%
\Index{texLinearMap}
   {$\star D$\Hyph product of \drc linear map $A$ over scalar}%
   {starD product of drc linear map over scalar}%
\Index{texVectorSpace}
   {$\star R$\hyph module}%
   {starR-module}%
\Index{texTstarRepresentation}
   {\sT shift}%
   {starT shift}%
\Index{texFiberedGroup}
   {\sT shift on fibered group}%
   {starT shift, fibered group}%
\Index{texTstarMorphism}
   {\sT representation of $\mathfrak{F}$\Hyph algebra $A$ in set $M$}%
   {starT representation of algebra}%
\Index{texFiberedAlgebra}
   {\sT representation of fibered $\mathfrak{F}$\Hyph algebra}%
   {starT representation of fibered F-algebra}%
\Index{texFiberedGroup}
   {\sT representation of fibered group}%
   {starT representation of fibered group}%
\Index{texTstarMorphism}
   {\sT transformation}%
   {starT transformation}%
\Index{texFiberedAlgebra}
   {\sT transformation on bundle}%
   {starT transformation of bundle}%
\Index{texDivisionRing}
   {structural constants of division ring $D$ over field $F$}%
   {structural constants of division ring over field}%
\Index{texBundle}
   {subbundle}%
   {subbundle}%
\Index{texVectorField}
   {subbundle of $\mathcal D\star$\hyph vector space}%
   {subbundle of Dstar vector bundle}%
\Index{texLinearMap}
   {sum of \drc linear maps}%
   {sum of drc linear maps, drc vector spaces}%
\Index{texDrcBasis}
   {sum of geometrical objects}%
   {sum of geometrical objects}%
\Index{texBasis}
   {sum of geometrical objects in vector space}%
   {sum of geometrical objects, vector space}%
\Index{texBundleRelation}
   {symmetric $2$\Hyph ary fibered relation}%
   {symmetric 2 ary fibered relation}%
\Index{texEuclideanSpace}
   {symmetric biadditive map}%
   {symmetric biadditive map}%
\Index{texDrcBasis}
   {symmetry group}%
   {SymmetryGroup}%
\Index{texBasis}
   {symmetry group}%
   {symmetry group}%
\Index{texGenRelativity}
   {synchronization of reference frame}%
   {synchronization of reference frame}%
\SetIndexSpace%
\Index{texSecondDifferential}
   {Taylor polynomial}%
   {Taylor polynomial, division ring}%
\Index{texSecondDifferential}
   {Taylor series}%
   {Taylor series, division ring}%
\Index{texTensorProduct}
   {tensor product of $D$\Hyph vector spaces}%
   {tensor product of D vector spaces}%
\Index{texTensorProduct}
   {tensor product of division rings}%
   {tensor product of division rings}%
\Index{texTensorProduct}
   {tensor product of \Ds vector spaces}%
   {tensor product of Dstar vector spaces}%
\Index{texLie}
   {tensor product of representations}%
   {tensor product of representations}%
\Index{texTensorProduct}
   {tensor product of rings over commutative ring}%
   {tensor product of rings}%
\Index{texSecondDifferential}
   {topological $D$\Hyph vector space}%
   {topological D vector space}%
\Index{texDifferential}
   {topological division ring}%
   {topological division ring}%
\Index{texSecondDifferential}
   {topological \drc vector space}%
   {topological drc vector space}%
\Index{texGeomObject}
   {torsion form}%
   {torsion form}%
\Index{texGeomObject}
   {torsion tensor}%
   {torsion tensor}%
\Index{texFiberedMorphism}
   {tower of bundles}%
   {tower of bundles}%
\Index{texTowerRepresentation}
   {tower of representations of $\mathfrak{F}$\Hyph algebras}%
   {tower of representations of algebras}%
\Index{texTstarMorphism}
   {transformation coordinated with equivalence}%
   {transformation coordinated with equivalence}%
\Index{texTstarMorphism}
   {transformation of set}%
   {transformation of set}%
\Index{texFiberedAlgebra}
   {transformation on bundle}%
   {transformation of bundle}%
\Index{texBundleRelation}
   {transitive $2$\Hyph ary fibered relation}%
   {transitive 2 ary fibered relation}%
\Index{texTstarMorphism}
   {transitive representation of $\mathfrak{F}$\Hyph algebra $A$}%
   {transitive representation of algebra}%
\Index{texFiberedAlgebra}
   {transitive representation of fibered $\mathfrak{F}$\Hyph algebra}%
   {transitive representation of fibered F-algebra}%
\Index{texRepresentation}
   {transitive representation of group}%
   {transitive representation of group}%
\Index{texVectorSpace}
   {\Ts linear composition of  vectors}%
   {linear composition of  vectors}%
\Index{texVectorSpace}
   {\Ts matrices vector space}%
   {matrices vector space}%
\Index{texTstarRepresentation}
   {\Ts shift}%
   {Tstar shift}%
\Index{texTstarMorphism}
   {\Ts representation of $\mathfrak{F}$\Hyph algebra $A$ in set $M$}%
   {Tstar representation of algebra}%
\Index{texFiberedAlgebra}
   {\Ts representation of fibered $\mathfrak{F}$\Hyph algebra}%
   {Tstar representation of fibered F-algebra}%
\Index{texTstarMorphism}
   {\Ts transformation}%
   {Tstar transformation}%
\Index{texFiberedAlgebra}
   {\Ts transformation on bundle}%
   {Tstar transformation of bundle}%
\Index{texLinearMap}
   {twin representations of division ring}%
   {twin representations of division ring}%
\Index{texFiberedGroup}
   {twin representations of fibered group}%
   {twin representations of fibered group}%
\Index{texTstarRepresentation}
   {twin representations of group}%
   {twin representations of group}%
\SetIndexSpace%
\Index{texDifferential}
   {unit sphere in division ring}%
   {unit sphere in division ring}%
\Index{texVectorField}
   {unitarity law for  $\mathcal D\star$\Hyph vector fields}%
   {unitarity law, Dstar vector fields}%
\Index{texVectorSpace}
   {unitarity law for $D\star$\Hyph vector space}%
   {unitarity law, Dstar vector space}%
\SetIndexSpace%
\Index{texDifferential}
   {valued division ring}%
   {valued division ring}%
\Index{texFiberedAlgebra}
   {vector bundle}%
   {vector bundle}%
\Index{texNewton}
   {vector potential}%
   {vector potential}%
\Index{texVectorSpace}
   {vector space type}%
   {vector space type}%

\CloseIndex

%% file: Symbol.English.tex
\def\indexname{Special Symbols and Notations}
\OpenIndex

\SetIndexSpace
\Symb{texBiring}
   {$(^a_b)$\hyph\CR quasideterminant}%
   {a b CR quasideterminant definition}%
\Symb{texBiring}
   {$(^a_b)$\hyph \RC quasideterminant}%
   {a b RC-quasideterminant definition}%
\Symb{texBiring}
   {minor}%
   {A from b a}%
\Symb{texBiring}
   {minor}%
   {A from columns T}%
\Symb{texBiring}
   {minor}%
   {A from rows S}%
\Symb{texBiring}
   {minor}%
   {A without column a}%
\Symb{texBiring}
   {minor}%
   {A without columns T}%
\Symb{texBiring}
   {minor}%
   {A without row b}%
\Symb{texBiring}
   {minor}%
   {A without rows S}%
\Symb{texAffineSpace}
   {affine space}%
   {affine space}%
\Symb{texBasis}
   {affine space}%
   {An}%
\Symb{texBiring}
   {\subs row ($c$\hyph row) of matrix}%
   {c row}%
\Symb{texAdditiveMap}
   {component of linear map $\Vector{A}$ of $D$\Hyph vector space}%
   {component of linear map, D vector space}%
\Symb{texAdditiveMap}
   {component $p$ of polyadditive map $\Vector A$}%
   {component of polyadditive map, D vector space}%
\Symb{texBiring}
   {\CR power of element $A$ of biring}%
   {cr power}%
\Symb{texBiring}
   {\CR inverse element of biring}%
   {cr-inverse element}%
\Symb{texBiring}
   {\CR product of matrices}%
   {cr-product of matrices}%
\Symb{texVectorSpace}
   {\dcr vector}%
   {dcr vector}%
\Symb{texLie}
   {derivative of left shift}%
   {derivative of left shift}%
\Symb{texLie}
   {derivative of left shift}%
   {derivative of left shift, 1-Parameter Group}%
\Symb{texLie}
   {derivative of right shift}%
   {derivative of right shift}%
\Symb{texLie}
   {}%
   {derivative of right shift}%
\Symb{texLie}
   {derivative of right shift}%
   {derivative of right shift, 1-Parameter Group}%
\Symb{texLie}
   {derivative of left shift}%
   {derivative of Tstar shift}%
\Symb{texVectorSpace}
   {\drc vector}%
   {drc vector}%
\Symb{texSecondDifferential}
   {norm of map $\Vector A$ of normed $D$\hyph vector space}%
   {norm of map, D vector space}%
\Symb{texAffine}
   {derivative}%
   {overline nabla_l, definition 2}%
\Symb{texAdditiveMap}
   {partial additive map of variable $v^i$}%
   {partial additive map of variable}%
\Symb{texBiring}
   {\sups row ($r$\hyph row) of matrix}%
   {r row}%
\Symb{texBiring}
   {\RC power of element $A$ of biring}%
   {rc power}%
\Symb{texBiring}
   {\RC inverse element of biring}%
   {rc-inverse element}%
\Symb{texBiring}
   {\RC product of matrices}%
   {rc-product of matrices}%
\Symb{texBiring}
   {\RC quasideterminant}%
   {RC-quasideterminant definition}%
\Symb{texAdditiveMap}
   {set of additive maps of $D$\Hyph vector space $\Vector{V}$ into $D$\Hyph vector space $\Vector{W}$}%
   {set additive maps, D vector space}%
\Symb{texRingAdditiveMap}
   {set of additive maps of ring $R_1$ into ring $R_2$}%
   {set additive maps, ring}%
\Symb{texAdditiveMap}
   {set of polyadditive mappings}%
   {set polyadditive maps, D vector space}%
\Symb{texRingAdditiveMap}
   {set of polyadditive maps of rings $R_1$, ..., $R_n$ into module $S$}%
   {set polyadditive maps, ring}%
\Symb{texRingAdditiveMap}
   {set of polylinear maps of rings $R_1$, ..., $R_n$ into module $S$}%
   {set polylinear maps, ring}%
\Symb{texPolyvector}
   {skew product of vectors $\Vector a_1$, ..., $\Vector a_m$}%
   {skew product of vectors}%
\Symb{texTstarRepresentation}
   {right shift}%
   {starT shift}%
\Symb{texFiberedGroup}
   {\sT shift}%
   {starT shift, fibered group}%
\Symb{texTstarRepresentation}
   {left shift}%
   {Tstar shift}%
\Symb{texFiberedGroup}
   {\Ts shift}%
   {Tstar shift, fibered group}%
\Symb{texRefernceFrame}
   {anholonomic coordinates of vector}%
   {vector anholonomic coordinates}%
\Symb{texRefernceFrame}
   {holonomic coordinates of vector}%
   {vector holonomic coordinates}%

\SetIndexSpace
\Symb{texBasis}
   {basis manifold of affine space}%
   {BAn}%
\Symb{texBasis}
   {basis manifold of vector space}%
   {basis manifold of vector space}%
\Symb{texBasis}
   {basis manifold of vector space $\mathcal{V}$}%
   {basis manifold of vector space}%
\Symb{texBasis}
   {basis manifold of central affine space}%
   {BCAn}%
\Symb{texBasis}
   {basis manifold of Euclid space}%
   {BEn}%
\Symb{texBundle}
   {Cartesian power $\mathcal{A}$ of bundle $\mathcal{B}$}%
   {Cartesian power of bundle}%
\Symb{texBundle}
   {Cartesian power $A$ of set $B$}%
   {Cartesian power of set}%
\Symb{texAffineSpace}
   {basis manifold of affine space}%
   {FAn}%
\Symb{texAffineSpace}
   {basis manifold of central affine space}%
   {FCAn}%
\Symb{texAffineSpace}
   {basis manifold of Euclid space}%
   {FEn}%
\Symb{texPolymodule}
   {product of objects $B_1$, ..., $B_n$ in category $\mathcal A$}%
   {product of objects in category, 1 n}%
\Symb{texDivisionRing}
   {structural constants of division ring $D$ over field $F$}%
   {structural constants of division ring over field}%

\SetIndexSpace
\Symb{texBasis}
   {central affine space}%
   {CAn}%
\Symb{texAffineSpace}
   {central affine space}%
   {central affine space}%
\Symb{texLie}
   {left structural constant of Lie algebra}%
   {left structural constant of Lie algebra}%
\Symb{texLie}
   {right structural constant of Lie algebra}%
   {right structural constant of Lie algebra}%

\SetIndexSpace
\Symb{texLieRepresentation}
   {basis vector of \sT representation}%
   {basis vector of starT representation}%
\Symb{texLieRepresentation}
   {basis vector of \sT representation}%
   {basis vector of starT representation, coordinates}%
\Symb{texLieRepresentation}
   {basis vector of \Ts representation}%
   {basis vector of Tstar representation}%
\Symb{texLieRepresentation}
   {basis vector of \Ts representation}%
   {basis vector of Tstar representation, coordinates}%
\Symb{texVectorSpace}
   {\subs rows \dcr vector space}%
   {c rows dcr vector space}%
\Symb{texSecondDifferential}
   {component of the G\^ateaux derivative of map $\Vector f(\Vector x)$}%
   {component of Gateaux derivative of map, D vector space}%
\Symb{texSecondDifferential}
   {component of the G\^ateaux derivative of map $\Vector f(\Vector x)$}%
   {component of Gateaux derivative of map, D vector space, short}%
\Symb{texDifferential}
   {component of the G\^ateaux differential of map $f(x)$}%
   {component of Gateaux derivative of map, division ring}%
\Symb{texSecondDifferential}
   {component of the G\^ateaux derivative of second order of map $\Vector f(\Vector x)$}%
   {component of Gateaux derivative of Second Order, D vector space}%
\Symb{texDifferential}
   {component of the G\^ateaux derivative of second order of map $f(x)$ of division ring}%
   {component of Gateaux derivative of Second Order, division ring}%
\Symb{texVectorField}
   {coordinate \Drc vector bundle}%
   {coordinate drc vector bundle}%
\Symb{texVectorSpace}
   {coordinate \drc vector space}%
   {coordinate drc vector space}%
\Symb{texRefernceFrame}
   {coordinate reference frame}%
   {coordinate reference frame, extensive definition}%
\Symb{texBundleRelation}
   {diagonal in bundle $\mathcal{A}$}%
   {diagonal in bundle, 1}%
\Symb{texPolymodule}
   {direct product of division rings $D_1$, ..., $D_n$}%
   {direct product of division rings, 1 n}%
\Symb{texFrechet}
   {the Fr\'echet \Ds derivative of map $f$ of division ring}%
   {Frechet Dstar derivative of map, division ring}%
\Symb{texSecondDifferential}
   {the G\^ateaux \crd derivative of map $\Vector f$ of $D$\hyph vector space $\Vector V$ to $D$\hyph vector space $\Vector W$}%
   {Gateaux crd derivative of map, D vector space}%
\Symb{texSecondDifferential}
   {the G\^ateaux derivative of map $\Vector f$ of normed $D$\Hyph vector space $\Vector{V}$ to normed $D$\Hyph vector space $\Vector{W}$}%
   {Gateaux derivative of map, D vector space}%
\Symb{texDifferential}
   {the G\^ateaux derivative of map $f$}%
   {Gateaux derivative of map, division ring}%
\Symb{texDifferential}
   {the G\^ateaux derivative of map $f$}%
   {Gateaux derivative of map, fraction, division ring}%
\Symb{texSecondDifferential}
   {the G\^ateaux derivative of order $n$ of map $\Vector f$}%
   {Gateaux derivative of Order n, D vector space}%
\Symb{texSecondDifferential}
   {the G\^ateaux derivative of order $n$ of map $f$ of division ring}%
   {Gateaux derivative of Order n, division ring}%
\Symb{texSecondDifferential}
   {the G\^ateaux derivative of order $n$ of map $f$ of division ring}%
   {Gateaux derivative of Order n, fraction, division ring}%
\Symb{texSecondDifferential}
   {the G\^ateaux derivative of second order of map $\Vector f$}%
   {Gateaux derivative of Second Order, D vector space}%
\Symb{texDifferential}
   {the G\^ateaux derivative of second order of map $f$ of division ring}%
   {Gateaux derivative of Second Order, division ring}%
\Symb{texDifferential}
   {the G\^ateaux derivative of second order of map $f$ of division ring}%
   {Gateaux derivative of Second Order, fraction, division ring}%
\Symb{texSecondDifferential}
   {the G\^ateaux differential of map $\Vector f$ of normed $D$\Hyph vector space $\Vector{V}$ to normed $D$\Hyph vector space $\Vector{W}$}%
   {Gateaux differential of map, D vector space}%
\Symb{texDifferential}
   {the G\^ateaux differential of map $f$}%
   {Gateaux differential of map, division ring}%
\Symb{texSecondDifferential}
   {the G\^ateaux differential of second order of mapping $\Vector f$}%
   {Gateaux differential of Second Order, D vector space}%
\Symb{texDifferential}
   {the G\^ateaux differential of second order of mapping $f$ of division ring}%
   {Gateaux differential of Second Order, division ring}%
\Symb{texSecondDifferential}
   {the G\^ateaux \drc derivative of map $\Vector f$ of $D$\Hyph vector space $\Vector V$ to $D$\Hyph vector space $\Vector W$}%
   {Gateaux drc derivative of map, D vector space}%
\Symb{texDifferential}
   {the G\^ateaux \Ds derivative of map $f$ of division ring $D$}%
   {Gateaux Dstar derivative of map, division ring}%
\Symb{texSecondDifferential}
   {the G\^ateaux Jacobian of map of $D$\Hyph vector space}%
   {Gateaux Jacobian of map, D vector space}%
\Symb{texSecondDifferential}
   {the G\^ateaux partial \drc derivative of map $f^b$ with respect to variable $v^a$}%
   {Gateaux partial crd derivative of map, 1, D vector space}%
\Symb{texSecondDifferential}
   {the G\^ateaux partial \drc derivative of map $f^b$ with respect to variable $v^a$}%
   {Gateaux partial crd derivative of map, 2, D vector space}%
\Symb{texSecondDifferential}
   {the G\^ateaux partial \drc derivative of map $f^b$ with respect to variable $v^a$}%
   {Gateaux partial crd derivative of map, 3, D vector space}%
\Symb{texSecondDifferential}
   {the G\^ateaux mixed partial derivative of map $f^j$ with respect to variables $v^i$, $v^j$}%
   {Gateaux partial derivative of Second Order, D vector space}%
\Symb{texSecondDifferential}
   {the G\^ateaux partial derivative of map $f^j$ with respect to variable $v^i$}%
   {Gateaux partial derivative, D vector space}%
\Symb{texSecondDifferential}
   {the G\^ateaux partial \drc derivative of map $f^b$ with respect to variable $v^a$}%
   {Gateaux partial drc derivative of map, 1, D vector space}%
\Symb{texSecondDifferential}
   {the G\^ateaux partial \drc derivative of map $f^b$ with respect to variable $v^a$}%
   {Gateaux partial drc derivative of map, 2, D vector space}%
\Symb{texSecondDifferential}
   {the G\^ateaux partial \drc derivative of map $f^b$ with respect to variable $v^a$}%
   {Gateaux partial drc derivative of map, 3, D vector space}%
\Symb{texDifferential}
   {the G\^ateaux \sD derivative of map $f$ of division ring $D$}%
   {Gateaux starD derivative of map, division ring}%
\Symb{texVectorSpace}
   {matrices vector space}%
   {matrices vector space}%
\Symb{texAffine}
   {Cartan derivative}%
   {overbrace D}%
\Symb{texAffine}
   {derivative}%
   {overline D}%
\Symb{texRefernceFrame}
   {derivative $e_{(k)}$}%
   {partial(k)}%
\Symb{texVectorSpace}
   {\sups rows \drc vector space}%
   {r rows drc vector space}%
\Symb{texTidal}
   {speed of deviation}%
   {speed of deviation}%
\Symb{texDifferential}
   {standard component of the G\^ateaux differential of map $f$}%
   {standard component of Gateaux differential, division ring}%
\Symb{texTensorProduct}
   {tensor product of division rings}%
   {tensor product of division rings}%
\Symb{texTensorProduct}
   {tensor product of rings}%
   {tensor product of rings}%
\Symb{texVectorSpace}
   {vector space type}%
   {vector space type}%

\SetIndexSpace
\Symb{texAffineSpace}
   {affine basis}%
   {Affine Basis}%
\Symb{texBasis}
   {affine basis}%
   {Affine Basis}%
\Symb{texAffineSpace}
   {basis}%
   {basis}%
\Symb{texBasis}
   {basis of vector space}%
   {Basis e}%
\Symb{texBasis}
   {basis in vector space $\mathcal{V}$}%
   {basis in V}%
\Symb{texVectorSpace}
   {basis of vector space}%
   {basis, vector space}%
\Symb{texPolymodule}
   {basis of $(n)$\hyph vector space}%
   {basis,n vector space}%
\Symb{texCartesian}
   {Cartesian power of total spaces}%
   {Cartesian power of total spaces}%
\Symb{texCartesian}
   {Cartesian product of total spaces}%
   {Cartesian product of total spaces, definition 1}%
\Symb{texBasis}
   {central affine basis}%
   {Central Affine Basis}%
\Symb{texVectorField}
   {basis for \Drc vector bundle}%
   {drc basis, vector bundle}%
\Symb{texRefernceFrame}
   {form of reference frame}%
   {dual forms, reference frame}%
\Symb{texBasis}
   {Euclid space}%
   {En}%
\Symb{texAffineSpace}
   {Euclid space}%
   {En}%
\Symb{texAffineSpace}
   {pseudo Euclid space}%
   {Enm}%
\Symb{texBasis}
   {pseudo Euclid space}%
   {Enm}%
\Symb{texFiberedAlgebra}
   {identical transformation of bundle}%
   {identical transformation of bundle}%
\Symb{texBasis}
   {orthonornal basis}%
   {Orthonornal Basis}%
\Symb{texAdditiveMap}
   {quaternion algebra over the field $F$}%
   {quaternion algebra over the field}%
\Symb{texCartesian}
   {reduced Cartesian product of total spaces}%
   {reduced Cartesian product of total spaces, definition 1}%
\Symb{texFiberedAlgebra}
   {set of nonsingular \sT transformations of bundle $\mathcal{E}$}%
   {set of starT nonsingular transformations of bundle}%
\Symb{texFiberedAlgebra}
   {set of nonsingular \Ts transformations of bundle $\mathcal{E}$}%
   {set of Tstar nonsingular transformations of bundle}%
\Symb{texBasis}
   {standard coordinates of basis}%
   {standard coordinates of basis}%
\Symb{texRefernceFrame}
   {standard coordinates of reference frame}%
   {standard coordinates of reference frame}%
\Symb{texRefernceFrame}
   {vector field of reference frame}%
   {vector field of reference frame}%
\Symb{texBasis}
   {vector of basis}%
   {vector of basis}%

\SetIndexSpace
\Symb{texVectorSpace}
   {coordinates of basis in \subs rows \dcr vector space}%
   {basis coordinates, c rows dcr vector space}%
\Symb{texVectorSpace}
   {coordinates of basis in \sups rows \drc vector space}%
   {basis coordinates, r rows drc vector space}%
\Symb{texVectorSpace}
   {basis for \subs rows \dcr vector space}%
   {basis, c rows dcr vector space}%
\Symb{texVectorSpace}
   {basis for \sups rows \drc vector space}%
   {basis, r rows drc vector space}%
\Symb{texDiffEq}
   {central affine basis}%
   {Central Affine Basis}%
\Symb{texRingAdditiveMap}
   {component of linear map $f$ of division ring}%
   {component of linear map, division ring}%
\Symb{texRingAdditiveMap}
   {component of polylinear map of division ring}%
   {component of polylinear map, division ring}%
\Symb{texBundle}
   {fibered morphism from bundle $\mathcal{A}$ into $\mathcal{B}$}%
   {fibered morphism from A into B}%
\Symb{texBundleRelation}
   {filter $\mathfrak{F}$ converges to set $A$}%
   {filter converges}%
\Symb{texFiberedAlgebra}
   {homomorphism of fibered $\mathfrak{F}$\Hyph algebras}%
   {homomorphism of fibered F-algebras}%
\Symb{texBundleRelation}
   {inverse fibered correspondence}%
   {inverse fibered correspondence, 1}%
\Symb{texBundleRelation}
   {inverse reduced fibered correspondence}%
   {inverse reduced fibered correspondence, 1}%
\Symb{texCartesian}
   {map to Cartesian product}%
   {map to Cartesian product}%
\Symb{texDifferential}
   {norm of map $f$ of division ring}%
   {norm of map, division ring}%
\Symb{texTstarRepresentation}
   {representation orbit of group $G$}%
   {orbit of representation of group}%
\Symb{texAffineSpace}
   {orthonornal basis}%
   {Orthonornal Basis}%
\Symb{texRefernceFrame}
   {reference frame}%
   {reference frame}%
\Symb{texRefernceFrame}
   {reference frame, extensive definition}%
   {reference frame, extensive definition}%
\Symb{texRingAdditiveMap}
   {standard component of additive map $f$ over field $F$}%
   {standard component of additive map, division ring}%
\Symb{texBiadditiveMap}
   {standard component of biadditive map $f$ over field $F$}%
   {standard component of biadditive map, division ring}%
\Symb{texRingAdditiveMap}
   {standard component of linear map $f$ of division ring}%
   {standard component of linear map, division ring}%
\Symb{texRingAdditiveMap}
   {standard component of polylinear map $f$ of division ring}%
   {standard component of polylinear map, division ring}%
\Symb{texBiadditiveMap}
   {standard component of quadratic map $f$ over field $F$}%
   {standard component of quadratic map, division ring}%
\Symb{texTensorProduct}
   {standard component of tensor}%
   {standard component of tensor, division ring}%

\SetIndexSpace
\Symb{texVectorSpace}
   {\CR matrix group}%
   {cr-matrix group}%
\Symb{texFiberedMorphism}
   {fibered little group of section $h$}%
   {fibered little group}%
\Symb{texFiberedMorphism}
   {fibered stability group of section $h$}%
   {fibered stability group}%
\Symb{texLie}
   {algebra Lie of group Lie}%
   {g}%
\Symb{texLie}
   {left defined algebra Lie of group Lie}%
   {gl}%
\Symb{texAffineSpace}
   {affine transformation group}%
   {GLAn}%
\Symb{texBasis}
   {affine transformation group}%
   {GLAn}%
\Symb{texLie}
   {right defined algebra Lie of group Lie}%
   {gr}%
\Symb{texBasis}
   {group of homomorphisms of vector space $\mathcal{V}$}%
   {GV}%
\Symb{texTstarRepresentation}
   {little group of $x$}%
   {little group}%
\Symb{texFiberedGroup}
   {orbit of effective covariant \sT representation of fibered group}%
   {orbit of effective covariant starT representation of fibered group}%
\Symb{texTstarRepresentation}
   {orbit of effective covariant \sT representation of group}%
   {orbit of effective covariant starT representation of group}%
\Symb{texFiberedGroup}
   {orbit of effective covariant \Ts representation of fibered group}%
   {orbit of effective covariant Tstar representation of fibered group}%
\Symb{texTstarRepresentation}
   {orbit of effective covariant \Ts representation of group}%
   {orbit of effective covariant Tstar representation of group}%
\Symb{texPolymodule}
   {product of groups $G_1$, ..., $G_n$}%
   {product of groups, 1 n}%
\Symb{texVectorSpace}
   {\RC matrix group}%
   {rc-matrix group}%
\Symb{texTstarRepresentation}
   {stability group of $x$}%
   {stability group}%

\SetIndexSpace
\Symb{texBiring}
   {Hadamard inverse of matrix}%
   {Hadamard inverse of matrix}%

\SetIndexSpace
\Symb{texLieRepresentation}
   {infinitesimal generator of representation}%
   {infinitesimal generator of representation}%
\Symb{texLinearLie}
   {Lie group infinitesimal generators}%
   {Lie group infinitesimal generators}%

\SetIndexSpace
\Symb{texAdditiveMap}
   {kernel of additive map of $D$\Hyph vector space}%
   {kernel of additive map, D vector space}%
\Symb{texAdditiveMap}
   {kernel of additive map of division ring}%
   {kernel of additive map, division ring}%

\SetIndexSpace
\Symb{texRepresentation}
   {left shift}%
   {left shift}%
\Symb{texDiffProperty}
   {Lie derivative of connection}%
   {Lie derivative of connection}%
\Symb{texDiffProperty}
   {Lie derivative of metric}%
   {Lie derivative of metric}%
\Symb{texBundleRelation}
   {limit of correspondence $\Phi$ with respect to the filter $\mathfrak{F}$}%
   {limit of correspondence with respect to the filter}%
\Symb{texDifferential}
   {limit of sequence in valued division ring}%
   {limit of sequence, valued division ring}%
\Symb{texBasis}
   {passive transformation}%
   {passive transformation}%
\Symb{texLinearMap}
   {\rcd vector space of \drc linear maps of \drc vector space $\Vector{V}$ into \drc vector space $\Vector{W}$}%
   {set drc linear maps, drc vector space}%
\Symb{texAdditiveMap}
   {set of linear maps of $D$\Hyph vector space $\Vector{V}$ into $D$\Hyph vector space $\Vector{W}$}%
   {set linear maps, D vector space}%
\Symb{texRepresentation}
   {set of left-side nonsingular transformations of set $M$}%
   {set of left-side nonsingular transformations}%
\Symb{texLinearMap}
   {\drc vector space of \rcd linear maps of \rcd vector space $\Vector{V}$ into \rcd vector space $\Vector{W}$}%
   {set rcd linear maps, rcd vector space}%
\Symb{texLinearMap}
   {set of \sT representations of division ring $S$ in additive group of division ring $R$}%
   {set sT representations, division ring}%
\Symb{texLinearMap}
   {set of \Ts representations of division ring $S$ in additive group of division ring $R$}%
   {set Ts representations, division ring}%

\SetIndexSpace
\Symb{texTstarMorphism}
   {set of \sT transformations of set $M$}%
   {set of starT transformations}%
\Symb{texTstarMorphism}
   {set of transformations of set $M$}%
   {set of transformations}%
\Symb{texTstarMorphism}
   {set of \Ts transformations of set $M$}%
   {set of Tstar transformations}%
\Symb{texTstarRepresentation}
   {space of orbits of effective \sT covariant representation of the group}%
   {space of orbits of effective sT representation}%
\Symb{texTstarRepresentation}
   {space of orbits of effective \Ts covariant representation of the group}%
   {space of orbits of effective Ts representation}%
\Symb{texTstarRepresentation}
   {space of orbits of \Ts representation $f$ of group $G$ in set $M$}%
   {space of orbits of Ts representation}%

\SetIndexSpace
\Symb{texBasis}
   {geometrical object in coordinate representation}%
   {geometrical object, coordinate vector space}%
\Symb{texBasis}
   {geometrical object}%
   {geometrical object, vector space}%
\Symb{texFiberedGroup}
   {orbit of representation of fibered group $\mathcal{G}$}%
   {orbit of representation of fibered group}%
\Symb{texRepresentation}
   {orbit of representation of the group $G$}%
   {orbit of representation of group}%

\SetIndexSpace
\Symb{texBundle}
   {bundle}%
   {bundle}%
\Symb{texFiberedMorphism}
   {bundle of level $2$}%
   {bundle of level 2}%
\Symb{texFiberedMorphism}
   {bundle of level $n$}%
   {bundle of level n}%
\Symb{texCartesian}
   {Cartesian power of bundle}%
   {Cartesian power of bundle}%
\Symb{texCartesian}
   {Cartesian product of bundles}%
   {Cartesian product of bundles, definition 1}%
\Symb{texCartesian}
   {reduced Cartesian product of bundles}%
   {reduced Cartesian product of bundles, definition 1}%
\Symb{texFiberedAlgebra}
   {set of nonsingular \sT transformations of bundle $\bundle{}pE{}$}%
   {set of starT nonsingular transformations of bundle, projection}%
\Symb{texFiberedAlgebra}
   {set of nonsingular \Ts transformations of bundle $\bundle{}pE{}$}%
   {set of Tstar nonsingular transformations of bundle, projection}%

\SetIndexSpace
\Symb{texBasis}
   {active transformation}%
   {active transformation}%
\Symb{texAffine}
   {Cartan curvature}%
   {Cartan curvature}%
\Symb{texVectorSpace}
   {\CR rank of matrix}%
   {cr-rank of matrix}%
\Symb{texBundleRelation}
   {diagonal in bundle  $\bundle{}pA{}$}%
   {diagonal in bundle, 2}%
\Symb{texBundleRelation}
   {diagonal in bundle $\mathcal{A}$}%
   {diagonal in reduced bundle, 2}%
\Symb{texLinearMap}
   {\Ds component of coordinates of vector $\Vector r$}%
   {Dstar component of coordinates of vector, D vector space}%
\Symb{texAffine}
   {curvature}%
   {GLn curvature_overline}%
\Symb{texVectorSpace}
   {\RC rank of matrix}%
   {rc-rank of matrix}%
\Symb{texRepresentation}
   {right shift}%
   {right shift}%
\Symb{texRepresentation}
   {set of right-side nonsingular transformations of set $M$}%
   {set of right-side nonsingular transformations}%
\Symb{texLinearMap}
   {\sD component of coordinates of vector $\Vector r$}%
   {starD component of coordinates of vector, D vector space}%

\SetIndexSpace
\Symb{texBundleRelation}
   {composition of fibered correspondences}%
   {composition of fibered correspondences}%
\Symb{texBundleRelation}
   {inverse fibered correspondence}%
   {inverse fibered correspondence, 2}%
\Symb{texBundleRelation}
   {inverse reduced fibered correspondence}%
   {inverse reduced fibered correspondence, 2}%
\Symb{texVectorSpace}
   {linear span in vector space}%
   {linear span, vector space}%

\SetIndexSpace
\Symb{texTstarMorphism}
   {category of \Ts representations of $\mathfrak{F}$\Hyph algebra $A$}%
   {category of Tstar representations of F algebra}%
\Symb{texTstarMorphism}
   {category of \Ts representations of $\mathfrak{F}$\Hyph algebra from category $\mathcal A$}%
   {category of Tstar representations of F algebra from category}%
\Symb{texLie}
   {tangent plane to group $G$}%
   {TaG}%

\SetIndexSpace
\Symb{texBasis}
   {coordinate vector space}%
   {coordinate vector space}%
\Symb{texBasis}
   {coordinates in vector space}%
   {coordinates in vector space}%
\Symb{texPolymodule}
   {direct product of $D_i\RCstar$\hyph vector spaces $\Vector V_1$, ..., $\Vector V_n$}%
   {direct product, drc vector space, 1 n}%
\Symb{texLinearMap}
   {dual space of \drc vector space $\Vector V$}%
   {dual space of drc vector space}%
\Symb{texVectorSpace}
   {\dcr vector space}%
   {left CR vector space}%
\Symb{texVectorSpace}
   {\drc vector space}%
   {left RC vector space}%
\Symb{texVectorSpace}
   {\crd vector space}%
   {right CR vector space}%
\Symb{texVectorSpace}
   {\rcd vector space}%
   {right RC vector space}%
\Symb{texTensorProduct}
   {tensor product of $D$\Hyph vector spaces}%
   {tensor product of D vector spaces}%
\Symb{texTensorProduct}
   {tensor product of \Ds vector spaces}%
   {tensor product of Dstar vector spaces}%
\Symb{texBasis}
   {vector space}%
   {V}%

\SetIndexSpace
\Symb{texBasis}
   {geometrical object in vector space}%
   {geometrical object, vector space}%

\SetIndexSpace
\Symb{texRefernceFrame}
   {anholonomic coordinate}%
   {x(k)}%

\SetIndexSpace
\Symb{texDivisionRing}
   {center of ring $D$}%
   {center of ring}%

\SetIndexSpace
\Symb{texTidal}
   {deviation of trajectories}%
   {deviation of trajectories}%
\Symb{texRepresentation}
   {identical transformation}%
   {identical transformation}%
\Symb{texTstarMorphism}
   {identical transformation}%
   {identical transformation}%
\Symb{texBasis}
   {image of vector $\Vector e_k\in\Basis e$ under isomorphism to coordinate vector space}%
   {image of vector e_k, coordinate vector space}%
\Symb{texBiring}
   {Kronecker symbol}%
   {Kronecker symbol}%

\SetIndexSpace
\Symb{texRefernceFrame}
   {anholonomic coordinates of connection}%
   {anholonomic coordinates of connection}%
\Symb{texAffine}
   {Cartan symbol}%
   {Cartan symbol}%
\Symb{texAffine}
   {connection}%
   {conection overline}%
\Symb{texRefernceFrame}
   {holonomic coordinates of connection}%
   {holonomic coordinates of connection}%
\Symb{texAffine}
   {Cartan connection}%
   {overbrace Gamma i kl}%
\Symb{texBundle}
   {set of sections of bundle}%
   {set of sections of bundle}%

\SetIndexSpace
\Symb{texLie}
   {inverse operator to operator $\psi_l$}%
   {inverse operator to operator psi l}%
\Symb{texLie}
   {inverse operator to operator $\psi_r$}%
   {inverse operator to operator psi r}%

\SetIndexSpace
\Symb{texRefernceFrame}
   {anholonomity object}%
   {anholonomity object}%

\SetIndexSpace
\Symb{texLie}
   {basic operator of group Lie}%
   {Lie Basic Operator L}%
\Symb{texLie}
   {}%
   {Lie Basic Operator L}%
\Symb{texLie}
   {basic operator of group Lie}%
   {Lie Basic Operator L, 1-Parameter Group}%
\Symb{texLie}
   {basic operator of group Lie}%
   {Lie Basic Operator R}%
\Symb{texLie}
   {}%
   {Lie Basic Operator R}%
\Symb{texLie}
   {basic operator of group Lie}%
   {Lie Basic Operator R, 1-Parameter Group}%

\SetIndexSpace
\Symb{texLie}
   {Lie group composition law}%
   {Lie group composition law}%

\SetIndexSpace
\Symb{texAffine}
   {Cartan derivative}%
   {overbrace nabla_l}%
\Symb{texAffine}
   {derivative}%
   {overline nabla_l, definition 1}%

\SetIndexSpace
\Symb{texBundleRelation}
   {restriction of correspondence $\Phi$ to set $C$}%
   {restriction of correspondence}%

\SetIndexSpace
\Symb{texCartesian}
   {Cartesian product of bundles}%
   {Cartesian product of bundles, definition 2}%
\Symb{texCartesian}
   {Cartesian product of total spaces}%
   {Cartesian product of total spaces, definition 2}%
\Symb{texPolymodule}
   {direct product of division rings $D_i$, $i\in I$}%
   {direct product of division rings}%
\Symb{texPolymodule}
   {direct product of division rings $D_1$, ..., $D_n$}%
   {direct product of division rings, i 1 n}%
\Symb{texPolymodule}
   {direct product of $D_i\RCstar$\hyph vector spaces $\Vector V_i$, $i\in I$}%
   {direct product, drc vector space}%
\Symb{texPolymodule}
   {direct product of $D_i\RCstar$\hyph vector spaces}%
   {direct product, drc vector space, i 1 n}%
\Symb{texPolymodule}
   {product of groups $G_i$, $i\in I$}%
   {product of groups}%
\Symb{texPolymodule}
   {product of groups $G_1$, ..., $G_n$}%
   {product of groups, i 1 n}%
\Symb{texPolymodule}
   {product of objects $\{B_i,i\in I\}$ in category $\mathcal A$}%
   {product of objects in category}%
\Symb{texPolymodule}
   {product of objects $B_1$, ..., $B_n$ in category $\mathcal A$}%
   {product of objects in category, i 1 n}%
\Symb{texCartesian}
   {reduced Cartesian product of bundles}%
   {reduced Cartesian product of bundles, definition 2}%
\Symb{texCartesian}
   {reduced Cartesian product of total spaces}%
   {reduced Cartesian product of total spaces, definition 2}%

\SetIndexSpace
\Symb{texBundle}
   {fibered subset}%
   {fibered subset}%
\Symb{texBundle}
   {subbundle}%
   {subbundle}%

\CloseIndex